\documentclass[10pt]{amsart}
\usepackage{amsfonts}
\usepackage{mathrsfs}
\usepackage{amscd}
\usepackage{amsmath}
\usepackage{amssymb}
\usepackage{latexsym}
\usepackage{lscape}
\usepackage{xypic}
\usepackage{comment}
\usepackage{amscd}
\usepackage{wasysym}  
\usepackage{tikz} 
\usetikzlibrary{matrix,arrows}
\usepackage{appendix}
\usepackage{geometry}
\usepackage[T1]{fontenc}
\usepackage{lmodern}
\usepackage[nice]{nicefrac}
\usepackage{dsfont, stmaryrd}
\usepackage{varwidth} 
\usepackage{dsfont}

\usepackage[pagebackref,hyperindex,linktocpage=true]{hyperref}
\hypersetup{
    colorlinks,
    linkcolor={red!50!black},
    citecolor={blue!50!black},
    urlcolor={blue!80!black}
}

\usepackage{cleveref}

\newtheorem{thm}{Theorem}[section]

\newtheorem{propo}[thm]{Proposition}
\newtheorem{lem}[thm]{Lemma}
\newtheorem{sublem}[thm]{Sublemma}
\newtheorem{lem-def}[thm]{Lemma-Definition}
\newtheorem{cor}[thm]{Corollary}
\newtheorem{conject}[thm]{Conjecture}
\newtheorem{propert}[thm]{Properties}
\newtheorem{observ}[thm]{Observation}

\theoremstyle{definition}
\newtheorem*{ack}{Acknowledgement}

\newtheorem{ex}[thm]{Example}

\newtheorem{rmk}[thm]{Remark}
\newtheorem{dfn}[thm]{Definition}
\newtheorem{quest}[thm]{Question}

\numberwithin{equation}{section}

\newcommand{\nc}{\newcommand}

\nc{\theo}{\begin{thm}} \nc{\xtheo}{\end{thm}}
\nc{\prop}{\begin{propo}} \nc{\xprop}{\end{propo}}
\nc{\lemm}{\begin{lem}} \nc{\xlemm}{\end{lem}}
\nc{\sublemm}{\begin{sublem}} \nc{\xsublemm}{\end{sublem}}
\nc{\lemmdefi}{\begin{lem-def}} \nc{\xlemmdefi}{\end{lem-def}}
\nc{\coro}{\begin{cor}} \nc{\xcoro}{\end{cor}}
\nc{\conj}{\begin{conject}} \nc{\xconj}{\end{conject}}
\nc{\proper}{\begin{propert}} \nc{\xproper}{\end{propert}}
\nc{\obse}{\begin{observ}} \nc{\xobse}{\end{observ}}
\nc{\ques}{\begin{quest}} \nc{\xques}{\end{quest}}

\nc{\ackn}{\begin{ack}} \nc{\xackn}{\end{ack}}
\nc{\exam}{\begin{ex}} \nc{\xexam}{\end{ex}}
\nc{\rema}{\begin{rmk}} \nc{\xrema}{\end{rmk}}
\nc{\defi}{\begin{dfn}} \nc{\xdefi}{\end{dfn}}

\nc{\pf}{\begin{proof}} \nc{\xpf}{\end{proof}}

\nc{\on}{\operatorname}
\nc{\fraka}{{\mathfrak a}} \nc{\bba}{{\mathbf a}}
\nc{\frakb}{{\mathfrak b}}
\nc{\frakc}{{\mathfrak c}}
\nc{\frakd}{{\mathfrak d}}
\nc{\frake}{{\mathfrak e}}
\nc{\frakf}{{\mathfrak f}}
\nc{\frakg}{{\mathfrak g}}
\nc{\frakh}{{\mathfrak h}}
\nc{\fraki}{{\mathfrak i}}
\nc{\frakj}{{\mathfrak j}}
\nc{\frakk}{{\mathfrak k}}
\nc{\frakl}{{\mathfrak l}}
\nc{\frakm}{{\mathfrak m}}
\nc{\frakn}{{\mathfrak n}}
\nc{\frako}{{\mathfrak o}}
\nc{\frakp}{{\mathfrak p}}
\nc{\frakq}{{\mathfrak q}}
\nc{\frakr}{{\mathfrak r}}
\nc{\fraks}{{\mathfrak s}}
\nc{\frakt}{{\mathfrak t}}
\nc{\fraku}{{\mathfrak u}}
\nc{\frakv}{{\mathfrak v}}
\nc{\frakw}{{\mathfrak w}}
\nc{\frakx}{{\mathfrak x}}
\nc{\fraky}{{\mathfrak y}}
\nc{\frakz}{{\mathfrak z}}
\nc{\frakA}{{\mathfrak A}}
\nc{\frakB}{{\mathfrak B}}
\nc{\frakC}{{\mathfrak C}}
\nc{\frakD}{{\mathfrak D}}
\nc{\frakE}{{\mathfrak E}}
\nc{\frakF}{{\mathfrak F}}
\nc{\frakG}{{\mathfrak G}}
\nc{\frakH}{{\mathfrak H}}
\nc{\frakI}{{\mathfrak I}}
\nc{\frakJ}{{\mathfrak J}}
\nc{\frakK}{{\mathfrak K}}
\nc{\frakL}{{\mathfrak L}}
\nc{\frakM}{{\mathfrak M}}
\nc{\frakN}{{\mathfrak N}}
\nc{\frakO}{{\mathfrak O}}
\nc{\frakP}{{\mathfrak P}}
\nc{\frakQ}{{\mathfrak Q}}
\nc{\frakR}{{\mathfrak R}}
\nc{\frakS}{{\mathfrak S}}
\nc{\frakT}{{\mathfrak T}}
\nc{\frakU}{{\mathfrak U}}
\nc{\frakV}{{\mathfrak V}}
\nc{\frakW}{{\mathfrak W}}
\nc{\frakX}{{\mathfrak X}}
\nc{\frakY}{{\mathfrak Y}}
\nc{\frakZ}{{\mathfrak Z}}
\nc{\bbA}{{\mathbb A}}
\nc{\bbB}{{\mathbb B}}
\nc{\bbC}{{\mathbb C}}
\nc{\bbD}{{\mathbb D}}
\nc{\bbE}{{\mathbb E}}
\nc{\bbF}{{\mathbb F}} \nc{\bbf}{{\mathbf f}}
\nc{\bbG}{{\mathbb G}}
\nc{\bbH}{{\mathbb H}}
\nc{\bbI}{{\mathbb I}}
\nc{\bbJ}{{\mathbb J}}
\nc{\bbK}{{\mathbb K}}
\nc{\bbL}{{\mathbb L}}
\nc{\bbM}{{\mathbb M}}
\nc{\bbN}{{\mathbb N}}
\nc{\bbO}{{\mathbb O}}
\nc{\bbP}{{\mathbb P}}
\nc{\bbQ}{{\mathbb Q}}
\nc{\bbR}{{\mathbb R}}
\nc{\bbS}{{\mathbb S}}
\nc{\bbT}{{\mathbb T}}
\nc{\bbU}{{\mathbb U}}
\nc{\bbV}{{\mathbb V}}
\nc{\bbW}{{\mathbb W}}
\nc{\bbX}{{\mathbb X}}
\nc{\bbY}{{\mathbb Y}}
\nc{\bbZ}{{\mathbb Z}}
\nc{\calA}{{\mathcal A}}
\nc{\calB}{{\mathcal B}}
\nc{\calC}{{\mathcal C}}
\nc{\calD}{{\mathcal D}}
\nc{\calE}{{\mathcal E}}
\nc{\calF}{{\mathcal F}}
\nc{\calG}{{\mathcal G}}
\nc{\calH}{{\mathcal H}}
\nc{\calI}{{\mathcal I}}
\nc{\calJ}{{\mathcal J}}
\nc{\calK}{{\mathcal K}}
\nc{\calL}{{\mathcal L}}
\nc{\calM}{{\mathcal M}}
\nc{\calN}{{\mathcal N}}
\nc{\calO}{{\mathcal O}}
\nc{\calP}{{\mathcal P}}
\nc{\calQ}{{\mathcal Q}}
\nc{\calR}{{\mathcal R}}
\nc{\calS}{{\mathcal S}}
\nc{\calT}{{\mathcal T}}
\nc{\calU}{{\mathcal U}}
\nc{\calV}{{\mathcal V}}
\nc{\calW}{{\mathcal W}}
\nc{\calX}{{\mathcal X}}
\nc{\calY}{{\mathcal Y}}
\nc{\calZ}{{\mathcal Z}}

\nc{\scrA}{{\mathscr A}}
\nc{\scrB}{{\mathscr B}}
\nc{\scrR}{{\mathscr R}}
\nc{\scrO}{{\mathscr O}}

\nc{\Bmu}{\mbox{$\raisebox{-0.59ex}{$l$}\hspace{-0.18em}\mu\hspace{-0.88em}\raisebox{-0.98ex}{\scalebox{2}{$\color{white}.$}}\hspace{-0.416em}\raisebox{+0.88ex}{$\color{white}.$}\hspace{0.46em}$}{}}

\nc{\bnu}{{\bar{ \nu}}}

\nc{\olO}{\bar{\calO}}

\nc{\al}{{\alpha}} 
\nc{\be}{{\beta}}
\nc{\ga}{{\gamma}} \nc{\Ga}{{\Gamma}}
 \nc{\hGa}{\hat{\Gamma}}
\nc{\ve}{{\varepsilon}} 
\nc{\la}{{\lambda}} \nc{\La}{{\Lambda}}
\nc{\om}{\omega} \nc{\Om}{\Omega} 
\nc{\sig}{{\sigma}} \nc{\Sig}{{\Sigma}}

\nc{\tnb}{\psi_{\rm tame}}
\nc{\oM}{\overline{{M}}}
\nc{\op}{{\on{op}}}
\nc{\ad}{{\on{ad}}}
\nc{\alg}{{\on{alg}}}
\nc{\Ad}{{\on{Ad}}}
\nc{\Adm}{{\on{Adm}}} \nc{\aff}{{\on{aff}}}
\nc{\Aut}{{\on{Aut}}}
\nc{\Bun}{{\on{Bun}}}
\nc{\cha}{{\on{char}}}
\nc{\der}{{\on{der}}}
\nc{\Der}{{\on{Der}}}
\nc{\Dist}{{\on{Dist}}}
\nc{\diag}{{\on{diag}}}
\nc{\End}{{\on{End}}}
\nc{\Fl}{{\on{Fl}}}
\nc{\Tr}{{\on{Transp}}}
\nc{\TR}{{\calT\!\calR}}
\nc{\Gal}{{\on{Gal}}}
\nc{\Gr}{{\on{Gr}}}
\nc{\rH}{{\on{H}}}
\nc{\Hom}{{\on{Hom}}}
\nc{\IC}{{\on{IC}}}
\nc{\id}{{\on{id}}}
\nc{\Id}{{\on{Id}}}
\nc{\ind}{{\on{ind}}}
\nc{\Ind}{{\on{Ind}}}
\nc{\Lie}{{\on{Lie}}}
\nc{\Pic}{{\on{Pic}}}
\nc{\pr}{{\on{pr}}}
\nc{\Res}{{\on{Res}}}
\nc{\res}{{\on{res}}} \nc{\Sat}{{\on{Sat}}}
\nc{\s}{{\on{sc}}}
\nc{\drv}{{\on{der}}}
\nc{\sgn}{{\on{sgn}}}
\nc{\Spec}{{\on{Spec}}}
\nc{\Spf}{\on{Spf}} 
\nc{\Sph}{\on{Sph}}
\nc{\St}{{\on{St}}}
\nc{\tr}{{\on{tr}}}
\nc{\Mod}{{\mathrm{-Mod}}}
\nc{\Hilb}{{\on{Hilb}}} 
\nc{\Ext}{{\on{Ext}}} 
\nc{\vs}{{\on{Vec}}}
\nc{\ev}{{\on{ev}}}
\nc{\nO}{{\breve{\calO}}}
\nc{\tS}{{\tilde{S}}}
\nc{\spe}{{\on{sp}}}
\nc{\loc}{{\on{loc}}}
\nc{\gr}{{\on{gr}}}

\nc{\nscrR}{{\mathscr{R}^{\on{nr}}}}
\nc{\GL}{{\on{GL}}}
\nc{\Gl}{\on{Gl}} 
\nc{\GSp}{{\on{GSp}}}
\nc{\gl}{{\frakg\frakl}}
\nc{\SL}{{\on{SL}}} 
\nc{\SU}{{\on{SU}}} 
\nc{\SO}{{\on{SO}}}
\nc{\PGL}{{\on{PGL}}}

\nc{\Conv}{{\on{Conv}}}
\nc{\Rep}{{\on{Rep}}}
\nc{\Dom}{{\on{Dom}}}
\nc{\red}{{\on{red}}}
\nc{\act}{{\on{act}}}
\nc{\nr}{{\on{nr}}}
\nc{\ctf}{{\on{ctf}}}

\nc{\str}{{\on{-}}} 
\nc{\os}{{\bar{s}}}
\nc{\oeta}{{\bar{\eta}}}

\nc{\hookto}{\hookrightarrow}
\nc{\longto}{\longrightarrow}
\nc{\leftto}{\leftarrow}
\nc{\onto}{\twoheadrightarrow}
\nc{\lonto}{\twoheadleftarrow}

\nc{\uG}{{\underline{G}}}
\nc{\uH}{{\underline{H}}}
\nc{\uA}{{\underline{A}}}
\nc{\uS}{{\underline{S}}}
\nc{\uT}{{\underline{T}}}
\nc{\uM}{{\underline{M}}}
\nc{\uP}{{\underline{P}}}
\nc{\uB}{{\underline{B}}}
\nc{\uN}{{\underline{N}}}
\nc{\uD}{{\underline{D}}}
\nc{\uU}{{\underline{U}}}

\nc{\ucG}{{\underline{\calG}}}
\nc{\ucH}{{\underline{\calH}}}
\nc{\ucA}{{\underline{\calA}}}
\nc{\ucS}{{\underline{\calS}}}
\nc{\ucT}{{\underline{\calT}}}
\nc{\ucM}{{\underline{\calM}}}
\nc{\ucP}{{\underline{\calP}}}
\nc{\ucN}{{\underline{\calN}}}
\nc{\ucU}{{\underline{\calU}}}
\nc{\ucO}{{\underline{\calO}}}
\nc{\ucC}{{\underline{\calC}}}

\nc{\bF}{{\breve{F}}}

\nc{\uFl}{{\underline{\Fl}}} 
\nc{\oFl}{{\overline{\Fl}}} 
\nc{\bU}{{\overline{U}}}
\nc{\tGr}{{\tilde{\Gr}}}
\nc{\cGr}{\calG\! r}
\nc{\oGr}{\overline{\on{Gr}}} 
\nc{\ocGr}{\overline{\calG\! r}}
\nc{\co}{{\colon}}
\nc{\sch}[1]{(Sch/{#1})}
\nc{\HypLoc}[1]{HypLoc({#1})}

\nc{\ohtimes}{\stackrel{!}{\otimes}}
\nc{\boxtilde}{\widetilde{\boxtimes}}
\nc{\vstar}{{\varhexstar}}

\nc{\Div}{\on{Div}}
\nc{\Sht}{\on{Sht}}
\nc{\Frob}{\on{Frob}}

\nc{\x}{\times}
\nc{\bsl}{\backslash}
\nc{\algQl}{{\bar{\bbQ}_\ell}}
\nc{\sF}{{\bar{F}}}
\nc{\nF}{{\breve{F}}}
\nc{\nW}{{W^{\on{nr}}}}
\nc{\sk}{{\bar{k}}}
\nc{\cont}{\on{c}}
\nc{\Supp}{\on{Supp}}
\nc{\blt}{\bullet}  
\nc{\dom}{\on{dom}}
\nc{\scon}{{\on{sc}}} 
\nc{\Affine}{\on{Aff}} 
\nc{\nscrA}{\mathscr{A}^{\on{nr}}} 
\nc{\nfraka}{{\bbf^{\on{nr}}}}
\nc{\ran}{{\rangle}}
\nc{\lan}{{\langle}}
\nc{\bk}{{\bar{k}}}
\nc{\tF}{{\tilde{F}}}
\nc{\sS}{{\bar{S}}}
\nc{\LG}{{^\text{L}\hspace{-0.04cm}G}}
\nc{\LL}{{^\text{L}\hspace{-0.07cm}L}}
\nc{\et}{{\text{\rm \'et}}}
\nc{\inv}{{\on{inv}}}
\nc{\Hecke}{{\on{Hecke}}}
\nc{\Isom}{{\on{Isom}}}
\nc{\oSht}{{\overline{\on{Sht}}}}
\nc{\umu}{{\underline \mu}}

\nc{\LRS}{{\on{LRS}}}
\nc{\IndSch}{{\on{IndSch}}}
\nc{\Pos}{{\on{Pos}}}
\nc{\Sets}{{\on{Sets}}}
\nc{\AffSch}{{\on{AffSch}}}
\nc{\Groups}{{\on{Groups}}}
\nc{\Rings}{{\on{Rings}}}
\nc{\Gpds}{{\on{Gpds}}}
\nc{\Sch}{{\on{Sch}}}
\nc{\fl}{{\on{flat}}}
\nc{\opp}{{\on{op}}}

\nc{\test}[1]{\framebox{\begin{varwidth}{0.9\textwidth} Test: #1 \end{varwidth}}}

\nc{\pot}[1]{ [\hspace{-0,5mm}[ {#1} ]\hspace{-0,5mm}] }
\nc{\rpot}[1]{ (\hspace{-0,7mm}( {#1} )\hspace{-0,7mm}) }

\nc{\defined}{\hspace{0.1cm}\stackrel{\text{\tiny \rm def}}{=}\hspace{0.1cm}}

\begin{document}

\title[Normality of Schubert varieties]{On the normality of Schubert varieties: \\ Remaining cases in positive characteristic\\ \hfill \\ Sur la normalit\'{e} des vari\'{e}t\'{e}s de Schubert:\\ les cas restants en caract\'{e}ristique positive}
\author[T.~J.~Haines, J.~Louren\c{c}o, T.~Richarz]{Thomas J.~Haines, Jo\~ao Louren\c{c}o, Timo Richarz}

\address{Department of Mathematics, University of Maryland, College Park, MD 20742-4015, DC, USA}
\email{tjh@umd.edu}

\address{Mathematisches Institut, Universität Münster, Einsteinstrasse 62, Münster, Germany}
\email{j.lourenco@uni-muenster.de}

\address{Technical University of Darmstadt, Department of Mathematics, 64289 Darmstadt, Germany}
\email{richarz@mathematik.tu-darmstadt.de}

\thanks{Research of T.H.~partially supported by NSF DMS-1801352, research of J.L.~funded by the Deutsche Forschungsgemeinschaft and research of T.R.~funded by the European Research Council (ERC) under Horizon Europe (grant agreement nº 101040935) and by the Deutsche Forschungsgemeinschaft (DFG, German Research Foundation) TRR 326 \textit{Geometry and Arithmetic of Uniformized Structures}, project number 444845124.}

\maketitle

\begin{abstract} 
We study the geometry of equicharacteristic partial affine flag varieties associated to tamely ramified groups $G$, with particular attention to the characteristic $p>0$ setting. We prove that when $p$ divides the order of the fundamental group $\pi_1(G_{\text{der}})$, most Schubert varieties attached to $G$ are not normal, and we provide a criterion for when normality holds. 
Apart from this, we show, on the one hand, that loop groups of semisimple groups satisfying $p \mid \# \pi_1(G_{\text{der}})$ are not reduced, and on the other hand, that their integral realizations are ind-flat.
Our methods allow us to classify all tamely ramified Pappas-Zhu local models of Hodge type which are normal.

\medskip

Nous \'{e}tudions la g\'{e}om\'{e}trie des vari\'{e}t\'{e}s de drapeaux affines partielles associ\'{e}es \`{a} des groupes $G$ mod\'{e}r\'{e}ment ramifi\'{e}s, avec un accent particulier sur le cadre de la caractéristique $p>0$. On démontre que, lorsque $p$ divise l'ordre du groupe fondamental $\pi_1(G_{\rm der})$, la plupart des vari\'{e}t\'{e}s de Schubert ne sont pas normales et nous fournissons une condition n\'{e}cessaire et suffisante pour que cela se produise. De plus, nous montrons, d'une part, que les groupes de lacets de groupes semisimples satisfaisant $p\,\mid\,\#\pi_1(G_{\rm der})$ ne sont pas r\'{e}duits, et d'autre part, que leurs r\'{e}alisations int\'{e}grales sont ind-plates. Nos m\'{e}thodes nous permettent de classifier tous les mod\`{e}les locaux de type Hodge au sens de Pappas-Zhu qui sont normaux.

\end{abstract}

\tableofcontents



\section{Introduction}

 Partial affine flag varieties are important objects in arithmetic algebraic geometry for their intimate relation to local models of Shimura varieties and moduli stacks of shtukas. 
 They first appeared extensively in the realm of Kac-Moody theory by means of (integral) representation theory of Kac-Moody algebras. 
 They were later reinterpreted via the theory of affine Grassmannians as parametrizing torsors under parahoric group schemes over the formal disk equipped with a trivialization over the punctured one.

In the works of Faltings \cite{Fal03}, Pappas-Rapoport \cite{PR08}, Zhu \cite{Zhu14} and Pappas-Zhu \cite{PZ13}, the authors establish several geometric properties of affine flag varieties, such as normality of Schubert varieties or reducedness of the special fiber of local models, under the following working hypothesis: the reductive group $G$ over the non-archimedean local base field is {\em tamely ramified}, and the residue characteristic $p>0$ {\em does not divide} the order of $\pi_1(G_{\text{der}})$, that is, the simply connected cover $G_\scon\to G_\der$ of the derived subgroup of $G$ is an \'etale isogeny.  
The first type of restriction has been substantially lifted in the work of Levin \cite{Lev16} for Weil-restricted groups, and in \cite{Lou19} for absolutely almost simple, wildly ramified groups. The second type of restriction is dealt with in this paper, whose main finding can be summarized as follows:
Let $F=k\rpot{t}$ be the Laurent series field in the formal variable $t$ with algebraically closed residue field $k$ of characteristic $p>0$.
Let $G$ be a tamely ramified connected reductive $F$-group, $\bbf$ a facet of its Bruhat-Tits building and $\bba$ an alcove containing $\bbf$ in its closure. 
For each class $w\in W/W_\bbf$ in the Iwahori-Weyl group quotient, let $S_w=S_w(\bba,\bbf)$ be the associated Schubert variety in the partial affine flag variety $\Fl_{G,\bbf}$. 
Note that $W/W_\bbf$ is always a (countable) infinite set when $G$ is nontrivial.

\theo[Prop.~\ref{normality.prop}, Thm.~\ref{fin.many.normal.sch.vars}, Prop.~\ref{general.facets.prop}, App.\,\ref{geom_props}] \label{main.thm.intro}
Assume $G$ is absolutely almost simple \textup{(}in particular, semisimple\textup{)}. If $p$ divides $\#\pi_1(G)$, then only finitely many Schubert varieties $S_w$, $w\in W/W_\bbf$ in the partial affine flag variety $\Fl_{G,\bbf}$ are normal.
The non-normal Schubert varieties are geometrically unibranch and regular in codimension $1$, but do not satisfy the $\textup{(}S2\textup{)}$ property, do not have rational singularities, and are neither Cohen-Macaulay, nor weakly normal, nor Frobenius split.

\xtheo

The existence of non-normal Schubert varieties in bad residue characteristics was first observed by the second named author. 
This came as a total surprise to us as these seem to be the very first examples of non-normal Schubert varieties in the literature.
The easiest such example occurs for the quasi-minuscule Schubert variety inside the affine Grassmannian for $G=\PGL_2$ in residue characteristic $p=2$:
the complete local ring at the singular point is isomorphic to the $k$-algebra
\[
k\pot{x,y,v,w}/(vw+x^2y^2, v^2+x^3y, w^2+xy^3, xw+yv).
\]
This is a surface singularity which is not weakly normal.
Its (weak) normalization morphism identifies with the inclusion map of the subalgebra of $k\pot{x,y,z}/(z^2+xy)$ generated by $x, y, v=xz, w=yz$ (see Appendix \ref{appendix_PGL2}).

The reason why non-normal Schubert varieties must exist can be summarized in a few lines. Up to translation by a suitable element in $G(F)$ which stabilizes $\bba$, we may assume that $S_w$ lies in the neutral component of $\Fl_{G,\bbf}$. 
By \cite[Prop. 6.6]{PR08}, the reduction of the neutral component identifies with that of $G_{\on{der}}$, so for this discussion we may assume that $G=G_{\on{der}}$ is semisimple.
Then one has a map 
\begin{equation}\label{intro.eq1}  
S_{\scon,w}=S_{\scon,w}(\bba,\bbf) \longto S_{w}(\bba,\bbf)=S_{w}
\end{equation}
where $S_{\scon,w}$ is the Schubert variety for $w$ inside $\Fl_{G_\scon,\bbf}$ and $G_\scon\to G$ is the simply connected cover.
The Schubert variety $S_{\scon,w}$ is known to be normal by \cite[Thm.~8.4]{PR08}, and the map \eqref{intro.eq1} can be shown to be finite, birational and a universal homeomorphism by using Demazure resolutions. 
In other words, the map \eqref{intro.eq1} is the (weak) normalization morphism of $S_w$, just as in the example of the quasi-minuscule Schubert variety above. 
On the other hand, the affine flag variety $\Fl_{G_\scon,\bbf}$ is reduced as an ind-scheme by \cite[Thm.~6.1]{PR08}, that is, equals the colimit of its Schubert varieties. 
If all Schubert varieties in $\Fl_{G,\bbf}$ were normal, then these two facts would imply the map $\Fl_{G_\scon,\bbf} \rightarrow \Fl_{G,\bbf}$ is a monomorphism. 
By looking at tangent spaces, this is clearly not true as soon as the kernel of $G_\scon\to G$ is non-\'etale, or equivalently, as soon as $p$ divides $\# \pi_1(G)$. 
We should however stress that the above reasoning only shows that there are infinitely many non-normal Schubert varieties in $\Fl_{G,\bbf}$ when $p$ divides $\#\pi_1(G)$.

For the rest of this introduction, we let $G$ denote an arbitrary tamely ramified connected reductive $F$-group. Exploiting tangent spaces a bit further, we show that the normality of $S_w$ is equivalent to the injectivity of the induced map $T_eS_{\scon,w}\to  T_eS_w$ on tangent spaces, which yields the following key observation:

\lemm[Cor.~\ref{normality.cor}] \label{normality.lemm.intro}
Let $w \in W/W_\bbf$. If $S_w$ is normal, then $S_v$ is normal for all $v\leq w$.
\xlemm

In order to give an effective normality criterion, we are led to a deeper study of tangent spaces of Schubert varieties for simply connected groups. 
In this, we recast in Section \ref{Kac-Moody-Interlude} results of Kumar \cite{Kum96}, Mathieu \cite{Mat89}, Ramanathan \cite{Ram87} and Polo \cite{Pol94} in the following fashion. 

We lift our whole setting to the Witt vectors ${\mathbb W}(k)$ as in \cite[\S\S7--9]{PR08}, and denote by $\uS_{\scon,w}\subset \uFl_{G_\scon,\bbf}$ the lift to ${\mathbb W}(k)$ of $S_{\scon,w}\subset \Fl_{G_\scon,\bbf}$ which comes equipped with a section $e\co \Spec\,{\mathbb W}(k)\to \uS_{\scon,w}$ given by the base point. 
Given any equivariant ample line bundle $\calL$ on $\uFl_{G_\scon,\bbf}$, we obtain the Kac-Moody action of $T_e\uFl_{G_\scon,\bbf}$ on the vector space $\Gamma(\uFl_{G_\scon,\bbf},\calL)^\vee$ dual to $\Gamma(\uFl_{G_\scon,\bbf},\calL)$; see Section \ref{schubert.tangent.sec}.

\theo[Lem.~\ref{tangent.schubert.lemm}] \label{tangent.thm.intro}
Assume $w \in W_{\rm aff}$. For any ${\mathbb W}(k)$-algebra $R$, the $R$-valued tangent space 
\[
T_e\uS_{\scon,w}(R) \;=\; \Hom_{{\mathbb W}(k)}(e^*\Om_{\uS_{\scon,w}/{\mathbb W}(k)},R)
\]
identifies with the submodule of $T_e\uFl_{G_\scon,\bbf}(R)$ consisting of those $X$ such that $X \Theta^\vee_\calL$ lies in the subspace $\Gamma(\uS_{\scon,w},\calL)^\vee\otimes R$, where $\Theta_\calL \in \Gamma(\uFl_{G_\scon,\bbf},\calL)$ is the usual theta divisor attached to $\calL$, and $\Theta^\vee_\calL \in \Gamma(\uFl_{G_\scon,\bbf},\calL)^\vee$ is the vector that sends $\Theta_\calL$ to $1$ and the remaining weight spaces to $0$.
\xtheo

This formula can in principle be used to determine whether a given Schubert variety is normal or not (see Corollary~\ref{normality.criterion.cor}).
We also think that it is of independent interest to have a good source for this material (some of which was known before in related contexts), and that having a Witt-vector framework which links to characteristic zero settings would potentially help in a future classification of all normal Schubert varieties when $p \mid \# \pi_1(G_{\text{der}})$.

It is not clear to us whether tangent spaces of Schubert varieties can be computed in a characteristic-independent way determined by the characteristic $0$ description, see Remark \ref{polo.remark} which comments on the argument in \cite[Cor.~4.1]{Pol94}.

The key to Theorem \ref{main.thm.intro} is to show that the tangent spaces of quasi-minuscule Schubert varieties in twisted affine Grassmannians for absolutely special vertices in characteristic $p>0$ are big enough, see Proposition \ref{quasi.minuscule}.
Thanks to some elementary observations (see Lemmas \ref{base_change_tangent_lem} and \ref{flat.coker.tangent}) the calculation can be reduced to characteristic $0$ where we identify the tangent spaces with tangent spaces at minimal nilpotent orbits, see Appendix \ref{app-minimal-nilpotent-orbits}. 
This uses the exponential map and representation-theoretic methods. 
For split groups, this relation to minimal nilpotent orbits is well known \cite[\S2.10]{MOV05}. 
For twisted groups, our method extends the method from \cite[\S8]{HR20a} and requires a fine analysis of twisted root systems. 
As a consequence, quasi-minuscule Schubert varieties are never normal if $p \mid \# \pi_1(G_{\text{der}})$. 
From here we use our key observation in Lemma \ref{normality.lemm.intro} along with combinatorial methods to finish the proof of Theorem \ref{main.thm.intro}.
In particular, we reprove in Proposition \ref{bruhat.order.on.coweights} some recent results from \cite[Thm.~4.1]{BH21} for split groups and also extend these to the case of twisted groups. 

Let us mention two other contributions of this paper to the understanding of the geometry of affine flag varieties: reducedness and ind-flatness. 
As we stated earlier, simply connected affine flag varieties are reduced and a similar result holds for all semisimple groups $G$ such that $p \nmid \# \pi_1(G)$ by \cite[Thm.~6.1]{PR08}. On the other hand, affine flag varieties of non-semisimple reductive groups are non-reduced. In \cite[Rem.~6.4]{PR08}, it is indicated that the affine flag variety for $\text{PGL}_2$ in characteristic $2$ is non-reduced.  The result below generalizes this.

\begin{thm}[Prop.~\ref{reduced.prop}, Prop.~\ref{reducedness_twisted}] \label{reduced.main.thm.intro}
The partial affine flag variety $\Fl_{G,\bbf}$ is reduced if and only if $G$ is semisimple and $p \nmid \# \pi_1(G)$.
\end{thm}

We give two different proofs of this result. If $G$ is split, we use the module of distributions, that is, higher differential operators of $\Fl_{G,\bbf}$ supported at the origin $e$, and we prove that the homomorphism $\text{Dist}(\Fl_{G_\scon,\bbf},e)\rightarrow \text{Dist}(\Fl_{G,\bbf},e)$ is not surjective in bad characteristic, implying non-reducedness, by essentially analyzing the effect of the multiplication-by-$p$ map on Grassmannians. If $G$ is tamely ramified, we factor the homomorphism $\Res_{F/F^p} G_\scon \rightarrow \Res_{F/F^p}G$ of pseudo-reductive groups as an epimorphism to a pseudo-reductive group $\overline{G}$ and a closed immersion whose image is strictly smaller than $\Res_{F/F^p}G$ - this works under the hypothesis that $G$ is semisimple and $p$ divides $\# \pi_1(G)$. Then we use the recently developed Bruhat-Tits theory for pseudo-reductive groups from \cite{Lou21} to prove that $\Fl_{\overline{G},\bbf} \rightarrow \Fl_{\Res_{F/F^p}G,\bbf}= \Fl_{G,\bbf}$ is a closed immersion, but not an isomorphism, for Lie-algebraic reasons.

Another natural question concerns the behavior of the integral realizations $\uFl_{G,\bbf}$ of the affine flag varieties over the Witt vectors (or just affine Grassmannians of split groups over integers). 
We are able to show:

\begin{thm}[Prop.~\ref{schubert.var.int.prop}, Prop.~\ref{ind.flat.prop}, Prop.~\ref{ind.flat.prop.twisted}] \label{flatness.thm.intro}
The ind-scheme $\uFl_{G,\bbf}$ is ind-flat over ${\mathbb W}(k)$.
It is reduced if and only if $G$ is semisimple.
In general, the reduced locus $(\uFl_{G,\bbf})_\red$ coincides with the union of the integral Schubert varieties $\uS_w=\uS_{w}(\bba,\bbf)$, $w\in W/W_\bbf$.
Furthermore, for fixed $w\in W/W_\bbf$, the following are equivalent:
\begin{enumerate}
\item The Schubert variety $S_w$ over $k$ is normal; 
\item its integral realization $\uS_w$ over ${\mathbb W}(k)$ is normal;
\item the special fiber $\uS_w\otimes k$ is reduced, hence equals $S_w$.
\end{enumerate}
\end{thm}

The proof of ind-flatness relies on computing the formal completion of the affine flag variety along the identity section, similar to Faltings' work \cite{Fal03}. 
For this, we compare the affine flag variety to its flat closure and it suffices, as both are ind-Noetherian, to show that their functors restricted to strictly Henselian Artinian local rings coincide, see Lemma \ref{what.faltings.knows}. 
This can be achieved by translating with the positive loop group and representative ${\mathbb W}(k)$-sections $\dot{w}$ of the Iwahori-Weyl group, so that those rings are supported at the identity section. 
Here we employ the fake open cell to reduce the ind-flatness to the cases of tori and unipotent groups where it is easy to check. 
The determination of the reduced locus is an immediate consequence because partial affine flag varieties for semisimple groups in characteristic $0$ are reduced. 
Finally, the equivalent conditions characterizing normal Schubert varieties are easily deduced by standard methods, see Proposition \ref{schubert.var.int.prop}.
We also refer the reader to Appendix \ref{appendix_PGL2} for the calculation of an integral Schubert variety whose reduction to characteristic $2$ is not reduced. 

Theorem \ref{flatness.thm.intro} is strongly connected with the theory of local models as follows. 
Let $F$ temporarily denote a discretely valued, complete field of characteristic $0$ with ring of integers $\calO_F$ and algebraically closed residue field $k$ of characteristic $p>0$, 
$G$ a tamely ramified reductive $F$-group and $S$ a maximal $F$-split torus of $G$. 
For each facet $\bbf$ in the appartment $\scrA(G,S,F)$ associated with $S$ of the Bruhat--Tits building $\scrB(G,F)$, we know by \cite{BT84} that there exists a canonical smooth, affine group $\calO_F$-scheme $\calG_\bbf$ with connected fibers whose generic fiber equals $G$ and whose $\calO_F$-valued points fix $\bbf$. Additionally, Pappas-Zhu \cite{PZ13} have constructed a smooth, affine, geometrically connected group scheme $\ucG_\bbf$ over $\calO_F[t]$ lifting $\calG_{\bbf}$ along the specialization $\calO_F[t]\rightarrow \calO_F$, sending $t$ to a preferred choice of uniformizer $\varpi \in \calO_F$, see \cite[Thm.~4.1]{PZ13}.
We note that the construction of this group scheme for split groups is easy and that the essential difficulty lies in its construction for twisted groups, see \cite[Exam.~3.3]{MRR20}. 
This group scheme is then used together with the Beilinson-Drinfeld affine Grassmannian \cite{BD91} to construct so-called Pappas-Zhu local models $\bbM=\bbM(G,\{ \mu\},\calG_\bbf)$ where $\{\mu\}$ is a geometric conjugacy class of cocharacters of $G$.
Recall that the reduced special fiber of a PZ local model is always given by the admissible locus $\calA(G, \{\mu\},\calG_{\bbf})$ (\cite[Thm.~6.12]{HR21}), that is, by a certain explicit union of Schubert varieties in the partial affine flag variety over $k$.

\coro[Cor.~\ref{local.models.corollary}] \label{local.models.corollary.intro}
Assume $p$ divides the order of $\pi_1(G_{\on{der}})$. 
\begin{enumerate}
\item If every Schubert variety in the admissible locus $\calA(G, \{\mu\},\calG_{\bbf})$ is normal, then $\bbM$ is normal and its special fiber is reduced. This is the case when $\bar{\mu}$ is minuscule for the \'echelonnage roots and $\bbf$ contains a special vertex in its closure.
\item If any Schubert variety inside the admissible locus $\calA(G, \{\mu\},\calG_{\bbf})$ is not normal, then $\bbM$ is not normal and its special fiber is not reduced.
\end{enumerate}
\xcoro

For details on part (1) we refer to Proposition \ref{minuscule} below and \cite[Thm.~2.1, Rem.~2.2]{HR22}. 
Here $\bar{\mu}$ is the image of a representative of the conjugacy class $\{\mu\}$ under the projection to inertia coinvariants $X_*(T) \rightarrow X_*(T)_I$.
For (2), suppose one of the Schubert varieties inside $\calA(G, \{\mu\},\calG_{\bbf})$ is not normal.
Then the irreducible component containing this Schubert variety is not normal as well by our key observation in Lemma \ref{normality.lemm.intro}. 
By comparing the Pappas-Zhu local model with its normalization (which is the Pappas-Zhu local model of some $z$-extension of $G$), we see that the special fiber cannot be reduced: compute global sections of line bundles and compare with the generic fiber by flatness. 
Hence, the Pappas-Zhu local model itself is not normal and its special fiber is not reduced, see also \cite[Rem.~2.4]{HR22}. 
In fact, this nuisance appears even if we assume that $\{ \mu\}$ is minuscule but $\bar{\mu}$ is sufficiently large for the échelonnage root system (which is possible if the ramification degree of $G$ is sufficiently large). More concretely, we give examples with restriction of scalars along ramified extensions or for unitary groups along ramified extensions, see Examples \ref{Weil.restriction.example} and \ref{unitary.example}.

Finally, we use Corollary \ref{local.models.corollary.intro} to classify in Proposition \ref{Hodge.classification} all tamely ramified PZ local models of Hodge type which are normal. We refer to Section \ref{Consequences.for.local.models} for the definition of Hodge and of abelian type, and we emphasize that in that section our groups are defined over a discretely valued complete field of characteristic $0$ with algebraically closed residue field of characteristic $p>0$. 

\prop[Prop.~\ref{Hodge.classification}]\label{Hodge.classification.intro}
Let $(G, \{\mu\})$ be of abelian type with a Hodge central lift $(G_1, \{\mu_1\})$, and let $\bbM_1$ be the PZ local model attached to $(G_1, \{\mu_1\},\calG_{\bbf,1})$.
Then the following properties hold:
\begin{enumerate}
\item If $p>2$ or $G_{\on{ad}}$ has no $D$-factors, then $\bbM_1$ is normal.

\item 
\label{boese.boese}
If $p=2$ and $(G_{\on{ad}},\{\mu_{\on{ad}}\})$ is simple of type $D_n^\bbH$, $n\geq 5$, then $\bbM_1$ is non-normal for all sufficiently large $\bar{ \mu}$.
\item If $p=2$ and $(G_{\on{ad}},\{\mu_{\on{ad}}\})$ is simple of type $D_{2m+1}^\bbR$, $m\geq 2$, then $\bbM_1$ is normal.
\item If $p=2$ and $(G_{\on{ad}},\{\mu_{\on{ad}}\})$ is simple of type $D_{2m}^\bbR$, $m\geq 2$, and $\bar{\mu}$ is sufficiently large, then $\bbM_1$ is normal if and only if $G_{1,\on{der}}=G_{1,\on{sc}}$.
\end{enumerate}
\xprop
 
These realizations are viewed via a Hodge embedding as flat closed subschemes of ordinary partial affine flag varieties for $\text{GL}_n$, as was done in \cite{KP18}.
The upshot is that, for $(G, \{\mu\})$ of abelian type, the Hodge embedding can always be arranged to give normal PZ local models except in Case \eqref{boese.boese}.
Also note that a Hodge embedding induces a corresponding closed immersion of the normalized local models only if $\bbM_1$ is normal, but that the corresponding morphism of topological spaces is always a closed immersion.

\subsection{Leitfaden (How to read this paper):} 
Assume ${\rm char}(k) \mid \#\pi_1(G_{\rm der})$.
Readers who are mainly interested in the existence and abundance of non-normal Schubert varieties in $\Fl_{G, \bbf}$ should start with \Cref{normality.prop}, which gives an elementary criterion for normality in terms of the tangent spaces, and which, together with its immediate corollaries, quickly shows that ``most'' Schubert varieties in $\Fl_{G, \bbf}$ are non-normal, once we know there is at least one non-normal Schubert variety (see the proof of \Cref{fin.many.normal.sch.vars}). 
The existence of a non-normal Schubert variety is explained in the introduction, and an alternative argument can be found in \Cref{5.13_rem}. \Cref{fin.many.normal.sch.vars} shows more precisely that when $G$ is absolutely almost simple, then $\Fl_{G, \bbf}$ contains only finitely many normal Schubert varieties. 
\Cref{geom_props} proves the equivalence of geometric properties asserted in \Cref{normality.prop}; the key new ingredient is the proof that all Schubert varieties are regular in codimension 1. 
\Cref{appendix_PGL2} works out the equations for the quasi-minuscule Schubert variety in the already non-trivial case of $G = {\rm PGL}_2$.
This and the above can be read independently of the rest of the paper.

Our second, more effective, criterion for normality of Schubert varieties is \Cref{normality.criterion.cor}, and this is used to develop an upper bound on the finite set of normal Schubert varieties attached to absolutely almost simple tamely ramified $k\rpot{t}$-groups $G$ in \Cref{general.facets.prop}. 
The latter relies on \Cref{quasi.minuscule} and its corollary which shows that the quasi-minuscule Schubert variety in an absolutely special affine Grassmannian is non-normal.
Much of Sections \ref{negative.loop.grps}-\ref{tangent.spaces} feeds into these propositions. 
\Cref{normality.criterion.cor} expresses the criterion in terms of negative loop groups and tangent spaces. 
It is essential for the proof of \Cref{quasi.minuscule} to develop both of these directions over the $p$-typical Witt vectors $\mathbb W = \mathbb W(k)$. More precisely, we construct a smooth, affine group scheme $\underline{\calG}_{\bbf}$ with connected fibers over $\mathbb W\pot{t}$ lifting any parahoric group schemes $\calG_{\bbf}$ over $k\pot{t}$ (see \Cref{PZ_group_scheme}). Their associated Schubert varieties $\underline{S}_w$ and partial affine flag varieties $\underline{\Fl}_{G,\bbf}$ over $\mathbb W$ are constructed in \Cref{negative.loop.grps}; (these lifts and \Cref{tangent.spaces} crucially reduce \Cref{quasi.minuscule} to the characteristic $0$ setting). 
The negative loop group $L^{--}_{\mathbb W} \underline{\calG}_{\bbf}$ is defined in \Cref{Neg_Loop_Grp_Twisted} and its isomorphism with the ``big open cell'' in $\underline{\Fl}_{G,\bbf}$ is proved in \Cref{Open_Cor}. With the goal of making \Cref{normality.criterion.cor} applicable, \Cref{tangent.space.cor} (resp.,\,\Cref{tangent.schubert.lemm}) gives an explicit description of $T_e \underline{\Fl}_{G_{\scon}, \bbf}$ (resp.\,of its subspace $T_e \underline{S}_{\scon, w}$). The latter relies on the comparison with Kac-Moody flag varieties over $\mathbb W$ (\Cref{kac.moody.comp}) and an extension of the tangent space description in that setting (\Cref{tangent.space.formula.and.independence.of.characteristic}) due to Kumar and Polo. \Cref{tangent.space.formula.and.independence.of.characteristic} rests in turn on \Cref{embeddings-given-by-quadrics-and-lines}, which gives the equations cutting out Schubert varieties under a projective embedding and in one another, extending results of Ramanathan and Mathieu; \Cref{Frob_ind-split} provides the technical ingredients of Frobenius splittings for ind-schemes. 
The cases of non-split groups in \Cref{quasi.minuscule} rely on root-theoretic computations for minimal nilpotent orbits done in \Cref{app-minimal-nilpotent-orbits}. 
 Finally, combinatorial results (\Cref{bruhat.order.on.coweights}, \Cref{explicit.bounds.prop}) complete the proof of \Cref{general.facets.prop}.

The remainder of the paper concerns reducedness, ind-flatness, and applications to Pappas-Zhu local models. The criterion for reducedness of $LG$ is proved in \Cref{reduced.prop} (for $G$ split) and in \Cref{reducedness_twisted} (in general), and this is independent of the other results in this paper. Similarly independent, the ind-flatness result is proved in  \Cref{ind.flat.prop} (for ${\rm Gr}_{G,\mathbb Z}$ attached to Chevalley groups $G$ over $\mathbb Z$) and in \Cref{ind.flat.prop.twisted} (for $\underline{\Fl}_{G, \bbf}$ attached to $\underline{\calG}_{\bbf}$ over $\mathbb W$).
Finally, \Cref{Consequences.for.local.models} deals with applications to local models of Shimura varieties. \Cref{local.models.corollary} makes the connection between properties of local models and the Schubert varieties in their special fibers, and \Cref{Hodge.classification} gives a classification of Pappas-Zhu local models of Hodge-type which are normal. 
Here we do assume background knowledge from the literature on local models.

\subsection{Acknowlegements} 
We thank Ulrich G\"{o}rtz, George Pappas, Michael Rapoport, Peter Scholze, Rong Zhou, and Xinwen Zhu for their interest and comments related to this work. We are also grateful to Shrawan Kumar and Peter Littelmann for interesting email exchanges and for providing us with some help in trying to understand an argument of Polo, see Remark~\ref{polo.remark}. 
In addition, we thank the referees for many remarks which led to enormous improvements in our exposition. 
Finally, the first and third authors heartily thank the second author for correcting an earlier version of this paper, which led to the present joint work.
\subsection{Conventions on ind-schemes}\label{conventions.sec}
We recall some basic results on ind-schemes, see \cite[\S1]{Ric20} for details.
An ind-scheme is a functor $X\co \AffSch^\opp\to \text{Sets}$  from the category of affine schemes such that there exists a presentation as functors $X=\on{colim}X_i$ where $\{X_i\}_{i\in I}$ is a filtered system of schemes $X_i$ with transition maps being closed immersions. 
Maps of ind-schemes are natural transformations of functors.
We denote by $\IndSch$ the category of ind-schemes which is locally small (i.e.,\,the {\rm Hom} classes are sets). 
 It contains the category of schemes as a full subcategory, is closed under fiber products and has $\Spec(\bbZ)$ as final object.
We identify $\AffSch^\opp$ with the category of rings whenever convenient. 
Note that every ind-scheme defines an fpqc sheaf on the category of affine schemes. 
Moreover, if $X=\on{colim}_iX_i$, $Y=\on{colim}_jY_j$ where each $X_i$ is quasi-compact, then
\[
\Hom_{\IndSch}(X,Y)\;=\; \on{colim}_i\on{lim}_j\Hom_{\Sch}(X_i,Y_j).
\]
If $S$ is a scheme, then an $S$-ind-scheme $X$ is an ind-scheme $X$ together with a map of functors $X\to S$. 
If $S=\Spec(R)$ is affine, we also use the term $R$-ind-scheme.

\section{A normality criterion} \label{normality.criterion}
Let $k$ be an algebraically closed field, and let $F=k\rpot{t}$ denote the Laurent series field. 
Let $G$ be a (connected) reductive $F$-group which splits over a tamely ramified extension of $F$. 

Let $\bbf\subset \scrB(G,F)$ be a facet in the Bruhat-Tits building, and denote by $\calG_\bbf$ the associated parahoric $\calO_F$-group scheme. 
The loop group $LG$ (resp.~$L^+\calG_\bbf$) is the functor on the category of $k$-algebras $R$ defined by $LG(R)=G(R\rpot{t})$ (resp.~$L^+\calG_\bbf(R)=\calG_\bbf(R\pot{t})$). Then $L^+\calG_\bbf\subset LG$ is a subgroup functor, and the {\it twisted affine flag variety} is the \'etale quotient
\[
\Fl_{G, \bbf}\defined LG/L^+\calG_\bbf,
\]
which is representable by an ind-projective $k$-ind-scheme. 
When $G$ is understood, we will often abbreviate by writing $\Fl_\bbf$ in place of $\Fl_{G, \bbf}$.

Let $S\subset G$ be a maximal $F$-split torus whose apartment $\scrA=\scrA(G,S,F)$ contains $\bbf$. 
We fix an alcove $\bba\subset \scrA$ which contains $\bbf$ in its closure. Fixing also a special vertex ${\bf 0}$ in the closure of ${\bf a}$, we may identify $\scrA$ with the vector space $X_*(S)_{\mathbb R}$, and, following Bruhat-Tits, we obtain an action of the Iwahori-Weyl group $W = W(G, S, F)$ on $\scrA$ and thus an isomorphism $W \overset{\sim}{\rightarrow} W_{\rm aff} \rtimes \Omega_{\bf a}$ where $W_{\rm aff}$ denotes the affine Weyl group and where $\Omega_{\bf a}$ is the subgroup of $W$ preserving ${\bf a}$. 
These basic notions related to Iwahori-Weyl groups can be found, for example, in \cite{PRS13, Ric16}. 
The left $L^+\calG_\bba$-orbits inside $\Fl_\bbf$ are enumerated by the quotient $W/W_\bbf$, where $W_\bbf\subset W_\aff$ is the subgroup of the affine Weyl group generated by the reflections fixing $\bbf$. 
For each class $w\in W/W_\bbf$, we define the Schubert variety  
\[
S_w=S_w(\bba,\bbf) \subset \Fl_{G,\bbf}
\]
as the reduced $L^+\calG_\bba$-orbit closure of $\dot{w}\cdot e$ where $e\in\Fl_{G,\bbf}(k)$ is the base point and $\dot{w}\in LG(k)$ is any representative of the class $w$. 
Then $S_w$ is a projective $k$-variety. The choice of $\bba$ equips the quotient $W/W_\bbf$ with a length function $l=l(\bba, \bbf)$ and a Bruhat partial order $\leq$ satisfying $\dim(S_w)=l(w)$, and $S_v\subset S_w$ if and only if $v\leq w$ for $v,w \in W/W_\bbf$, see \cite[Prop.~2.8]{Ric13}.

Let $\phi\co G_\scon\to G_\der\subset G$ be the simply connected cover. 
Let $T$ be the centralizer of $S$ in $G$ (a maximal torus by Steinberg's theorem) and let $T_{\rm der} := T \cap G_{\rm der}$.
Then $S_\scon:=\phi^{-1}(S)^o\subset \phi^{-1}(T_{\rm der})^o=:T_\scon$ is a maximal $F$-split torus contained in a maximal torus of the group $G_\scon$.
This induces a map on apartments $\scrA(G_\scon,S_\scon,F)\to \scrA(G,S,F)$ under which the facets correspond bijectively to each other. 
We denote the preimage of $\bbf$ by the same letter. The map $G_\scon\to G$ extends to a map on parahoric group schemes $\calG_{\scon,\bbf}\to \calG_{\bbf}$, and hence to a map on twisted partial affine flag varieties $\Fl_{G_\scon,\bbf}\to \Fl_{G,\bbf}^o$ onto the neutral component. 
We are interested in comparing their Schubert varieties.  

The natural map on Iwahori-Weyl groups 
\[
W_\scon=W(G_\scon,S_\scon,F)\longto W(G,S,F)=W,
\]
is injective and its image identifies with the affine Weyl group $W_\aff$ compatibly with the subgroup $W_\bbf$. 
Thus, for each class $w\in W_\aff/W_\bbf$ we get a map of projective $k$-varieties
\begin{equation}\label{schubert.map}
S_{\scon,w}=S_{\scon, w}(\bba,\bbf)\longto S_{w}(\bba,\bbf)=S_w.
\end{equation}

\prop \label{normality.prop}
For each class $w\in W_\aff/W_\bbf$, the following statements are equivalent:
\begin{enumerate}
\item The Schubert variety $S_w$ is normal \textup{(}resp.~weakly normal, resp.~$(S2)$, resp.~Cohen-Macaulay, resp.~Frobenius split if $\on{char}(k)>0$\textup{)}. 
\item The map \eqref{schubert.map} is an isomorphism.
\item The map \eqref{schubert.map} induces an injective map on tangent spaces at the base points. 
\end{enumerate}
\xprop

\pf
We will establish the equivalence of the geometric properties listed in (1) in Appendix \ref{geom_props}. It remains to show the following implications: \smallskip\\
(1)$\Rightarrow$(2): The map \eqref{schubert.map} is a finite birational universal homeomorphism by \cite[Prop.~3.5]{HR22}, and thus is an isomorphism whenever $S_w$ is normal.\smallskip\\
(2)$\Rightarrow$(1): Since $G$ splits over a tamely ramified extension of $F$, the Schubert variety $S_{\scon,w}\subset \Fl_{G_\scon,\bbf}$ is normal by \cite[Thm.~0.2]{PR08}, and so is $S_w$ whenever \eqref{schubert.map} is an isomorphism. \smallskip\\
(2)$\Rightarrow$(3): This is trivial. \smallskip\\
(3)$\Rightarrow$(2): The locus in $S_{\scon,w}$, where \eqref{schubert.map} is an isomorphism, is non-empty, open and $L^+\calG_{\scon,\bba}$-invariant. 
Thus, it suffices to show that the map of local rings at the base points
\[
\calO:=\calO_{S_{w}, e} \longto \calO_{S_{\scon, w},e}=:\calO_\scon
\] 
is an isomorphism. 
Here $e$ denotes the base point of both $\Fl_{G_\scon,\bbf}$ and $\Fl_{G,\bbf}$.
As \eqref{schubert.map} is a finite birational map between integral schemes, the map $\calO\hookto\calO_\scon$ is a finite ring extension which induces an isomorphism on fraction fields. 
Since we are assuming that the map \eqref{schubert.map} induces an injection on tangent spaces at the base points, we know that it is unramified at the base points by \cite[0B2G]{StaProj} so that $\frakm\calO_\scon=\frakm_\scon$ for the maximal ideals. 
An application of Nakayama's Lemma \cite[00DV (6)]{StaProj} to the finite map $\calO\hookto\calO_\scon$ (both viewed as $\calO$-modules) shows that this map is surjective as well. 
This finishes the proof of the proposition.
\xpf

\coro \label{normality.cor}
Let $w\in W_\aff/W_\bbf$. If $S_w$ is normal, then $S_v$ is normal for all $v\leq w$.
\xcoro
\pf
This is immediate from Proposition \ref{normality.prop} (3).
\xpf

Let us also record the following two useful results. 
We have an isomorphism $W \overset{\sim}{\rightarrow} W_{\rm aff} \rtimes \Omega_{\bf a}$, where $\Omega_{\bf a}$ is the subgroup of $W$ preserving ${\bf a}$, see Section \ref{absolutely.special}.  The properties of Schubert varieties indexed by $W$ reduce to those indexed by $W_{\rm aff}$, as follows: for any $w \in W_{\rm aff}$ and $\tau \in \Omega_{\bf a}$, there is an isomorphism
$$
S_{w\tau} ({\bf a}, {\bf f}) \cong S_w({\bf a}, \,\tau{\bf f}),
$$
where $\tau \bf f$ is the image of $\bf f$ under the action of $\tau$ on facets in the boundary of $\bf a$.

\prop \label{permanence}
Let $w \in W_{\rm aff}$ and $\tau \in \Omega_{\bf a}$. The following are equivalent:
\begin{enumerate}
\item $S_{w\tau\eta}(\bba, \bba)$ is normal for all $\eta \in W_{\bbf}$;
\item $S_{w\tau \eta_0}(\bba, \bba)$ is normal for $\eta_0 \in W_{\bbf}$ such that $w\tau \eta_0 $ is right $\bbf$-maximal;
\item $S_{w\tau}(\bba, \bbf)$ is normal.
\end{enumerate}
\xprop
\pf
By the above discussion, we immediately reduce to the case $\tau = e$. The implication (1)$\Rightarrow$(2) is obvious, and the opposite implication follows from Corollary \ref{normality.cor}. For (2)$\Leftrightarrow$(3) we use the fact that the inverse image of $S_w(\bba, \bbf)$ under the smooth surjective morphism $\Fl_{\bba} \to \Fl_{\bbf}$ is the Schubert variety $S_{w\eta_0}(\bba, \bba)$. We conclude by observing that normality is local for the smooth topology, see \cite[Tag 034F]{StaProj}.
\xpf

\lemm \label{normality.stabilizer.lem}
Let $\tau\in \Om_{\bba}$. 
For each class $w\in W_\aff/W_\bbf$, the $(\bba,\bbf)$-Schubert variety $S_w\subset \Fl_{G,\bbf}$ is normal \textup{(}resp.~smooth\textup{)} if and only if the $(\bba,\tau\bbf)$-Schubert variety $S_{\tau w\tau^{-1}}\subset \Fl_{G,\tau\bbf}$ is normal \textup{(}resp.~smooth\textup{)}.
\xlemm
\pf First note that the class of $\dot{\tau}\cdot \dot{w}\cdot \dot{\tau}^{-1}$ inside $W_\aff/W_{\tau\bbf}$ is well-defined where $\dot{\tau}, \dot{w}\in LG(k)$ are any representatives. 
Thus, the $(\bba,\tau\bbf)$-Schubert variety $S_{\tau w\tau^{-1}}$ is well-defined. 
Further, the isomorphism $LG\to LG, g\mapsto \dot{\tau}g\dot{\tau}^{-1}$ descends to an isomorphism $\Fl_{G,\bbf}\to \Fl_{G,\tau\bbf}$ mapping the $(\bba,\bbf)$-Schubert variety $S_w$ isomorphically onto the $(\bba,\tau\bbf)$-Schubert variety $S_{\tau w\tau^{-1}}$. 
This proves the lemma.
\xpf

Let us state one of the main results of the paper. In this paper, a group will be termed {\em absolutely almost simple} if its absolute Dynkin diagram is connected, and it will be called {\em absolutely simple} if it is absolutely almost simple and adjoint.

\theo \label{fin.many.normal.sch.vars}
Suppose $G$ is an absolutely almost simple, semisimple group such that its simply connected cover is a non-\'etale isogeny. 
Then $\Fl_{G, \bbf}$ contains only finitely many Schubert varieties which are normal.
\xtheo
\begin{proof}
 We already stated above that if $p$ divides $\#\pi_1(G_{\rm der})$, then there is at least one non-normal Schubert variety $S
_w$ in $\Fl_{G,\bbf}$. 
Corollary \ref{normality.cor} then implies there are infinitely many elements $w \in W/W_\bbf$ such that $S_w$ is not normal. In fact, as pointed out by a referee, we can go further and deduce already that there are only finitely many $w \in W/W_{\bbf }$ such that $S_w$ is normal.  Indeed, to show this, we can reduce to $\bbf = \bba$ and $w \in W_{\rm aff}$, and show that given an element $w \in W_{\rm aff}$, all but finitely many elements $v \in W_{\rm aff}$ satisfy $w \leq v$. For each $s \in S_{\rm aff}$, the parabolic subgroup of $W_{\rm aff}$ generated by
$S_{\rm aff} \backslash \{s\}$ is finite, so there is an integer $L_s$ such that all reduced words
in the alphabet $S_{\rm aff} \backslash \{s\}$ have length $< L_s$. Therefore any reduced word
of length $\geq L_s$ in the alphabet $S_{\rm aff}$ contains the letter $s$. Write $w$ as a
reduced word $s_1 \cdots  s_q$. Then any reduced word $v$ in the alphabet $S_{\rm aff}$ of
length $\geq L_{s_1} + \cdots + L_{s_q}$ contains the word $w$, in the sense that $v\geq w$
in the Bruhat order.
\end{proof}

The argument above does not give any explicit indication of how large $w$ must be in order that $S_w$ is necessarily non-normal. 
In Section \ref{classification}, we give explicit examples of non-normal $S_w$ for each absolutely almost simple group $G$ such that $p$ divides $\#\pi_1(G_{\rm der})$.  
Additionally, we use the intervening material as well as the appendices to give an explicit upper bound (in the Bruhat order) on the elements $w \in W/W_{\bbf}$ such that $S_w$ is normal (see Proposition \ref{general.facets.prop}).

Notice that we can easily find semisimple groups which are not absolutely almost simple whose affine flag varieties contain infinitely many normal Schubert varieties. 
However, it is still true that the great majority of them are not normal: indeed, as soon as all of their projections to the partial affine flag variety of the adjoint factors of $G$ have sufficiently large dimension, then the Schubert varieties are not normal.

\section{Tame liftings and negative loop groups} \label{negative.loop.grps}
In this section, we explain how to lift Schubert varieties from characteristic $p>0$ to characteristic $0$, and set the stage for the calculation of tangent spaces.

\subsection{Tame liftings of groups} \label{tame_lift_subsec}
Let $k$ be an algebraically closed field, and let $F=k\rpot{t}$ denote the Laurent series field. 
Let $G$ be an absolutely almost simple, tamely ramified reductive $F$-group, and assume that $G$ has the same splitting field as its simply connected group (equivalently, as its adjoint group). 
We follow the presentation in \cite[\S7]{PR08}, but we no longer assume that $G$ is simply connected. Let $F'/F$ be the tamely ramified splitting field of $G$. The extension $F'/F$ is a cyclic Galois extension of degree $e=1$, $e=2$ or $e=3$, cf.~\cite[\S7]{PR08}. We fix a uniformizer $u\in F'$ such that $u^e=t$, and a generator $\tau$ for the group $\lan\tau\ran=\Gal(F'/F)$. We have $\tau u=\zeta\cdot u$ where $\zeta$ is a primitive $e$-th root of unity.

Fix a Chevalley group $H$ over $\bbZ$ together with an isomorphism $G\otimes_FF'=H\otimes_\bbZ F'$ compatible with pinnings on both sides. The pinning for $G$ over $F$ is denoted $(G,B,T,X)$, and for $H$ over $\bbZ$ it is denoted $(H,B_H,T_H,X_H)$.  Here $T\subset G$ is the centralizer of the maximal $F$-split torus $S$ as above, and $X$ (similarly $X_H$) denotes the sum of a choice of root vectors in the Lie algebra of $G$ corresponding to the simple roots $\Delta = \Delta(G,B,T)$ for $(G,B,T)$. Recall that we fixed an alcove $\bba\subset \scrA(G,S,F)$ containing a facet $\bbf$ in its closure.
 
The automorphism ${\rm id} \otimes \tau$ of $G \otimes_F F'$ induces an automorphism $\sig$ on $H\otimes_\bbZ F'$ which can be written in the form $\sig=\sig_0\otimes \tau$ where $\sig_0 \in {\rm Aut}(X^*(T_H), \Delta_H, X_*(T_H), \Delta^\vee_H)$ is viewed as an order $e$ automorphism of $H$. 
Here $\Delta_H$ (resp.~$\Delta^\vee_H$) denotes the simple roots (resp.~coroots) for $(H, B_H, T_H)$. 
Explicitly, we have $G = \Res_{F'/F}(H_{F'})^\sigma$. There are identifications of buildings $\scrB(G,F) = \scrB(\Res_{F'/F}(H_{F'}), F)^\sigma = \scrB(H_{F'}, F')^\sigma$ (see \cite{PY02, HR20b}). 
We fixed a facet $\bbf\subset \scrB(G,F)$, which corresponds to $\sigma$-stable facet in $\scrB(H_{F'},F')$, also denoted $\bbf$.
The parahoric group scheme can now be written in the form
\begin{equation}\label{Res_Scalar_Grp}
\calG_\bbf\,=\, (\Res_{k\pot{u}/k\pot{t}}\calH_\bbf)^{\sig,o}, 
\end{equation}
where $\calH_\bbf$ is the parahoric group scheme associated with the $\sigma$-stable facet $\bbf \subset \scrB(H_{F'},F')$. 
Here $(\str)^o$ denotes the fiberwise neutral component which only plays a role if $G$ is not simply connected. 

 This leads to the identifications of loop groups
\begin{equation}\label{Relation_Loop_Grps}
LG\,=\, (LH_{k\rpot{u}})^\sig \;\;\; \text{and}\;\;\; L^+\calG_\bbf\,=\, (L^+\calH_\bbf)^{\sig,o},
\end{equation}
where we refer to the discussion below \eqref{Relation_Loop_Grps2} for the second equality. 

We now lift \eqref{Res_Scalar_Grp} and \eqref{Relation_Loop_Grps} to the Witt vectors. 
For this, assume that $k$ is of characteristic $p>0$, and denote by ${\mathbb W}={\mathbb W}(k)$ the ring of $p$-typical Witt vectors together with the natural map ${\mathbb W}\to k$. Let ${\mathbb K}=\on{Frac}({\mathbb W})$ be the field of fractions. Following the arguments in \cite[\S7]{PR08} (for Iwahori group schemes), or \cite[\S4.2.2 (a)]{PZ13}, we have the `parahoric group scheme' $\ucH_\bbf$ over the ring ${\mathbb W}\pot{u}$ such that $\ucH_\bbf\otimes k\pot{u}=\calH_\bbf$. The group $\ucH_\bbf$ is by construction `horizontal along the ${\mathbb W}$-direction', so that $\ucH_\bbf\otimes {\mathbb K}\pot{u}$ is an Iwahori group scheme of the same type as $\calH_\bbf$ (but now the residue field ${\mathbb K}$ is of characteristic zero). Note that $\ucH_\bbf\otimes {\mathbb W}\rpot{u}=H\otimes_\bbZ {\mathbb W}\rpot{u}$. Further, the automorphism $\tau$ lifts to the automorphism $\tau\co {\mathbb W}\pot{u}\to {\mathbb W}\pot{u}$, $u\mapsto [\zeta]\cdot u$ where $[\str]$ denotes the Teichm\"uller lift. Again we have the automorphism $\sig=\sig_0\otimes \tau$ on $\Res_{{\mathbb W}\pot{u}/{\mathbb W}\pot{t}}(\ucH_\bbf)$ so that we can define the fiberwise neutral component (cf.~[SGA 3, $\on{VI}_B$, \S3] for general base schemes)
\begin{equation}\label{Res_Scalar_Grp_Witt}
\ucG_\bbf\defined (\Res_{{\mathbb W}\pot{u}/{\mathbb W}\pot{t}}\ucH_\bbf)^{\sig,o}.
\end{equation}
By [SGA3, $\on{VI}_B$, Thm.~3.10], this is a smooth ${\mathbb W}\pot{t}$-group scheme with connected fibers such that $\ucG_\bbf\otimes_{\mathbb{W}\pot{t}} k\pot{t}=\calG_\bbf$, and Lemma \ref{PZ_group_scheme} below shows that it is affine as well. We define $\uG:=\ucG_\bbf\otimes_{\mathbb{W}\pot{t}} {\mathbb W}\rpot{t}$, and $\uH:=\ucH_\bbf\otimes_{\mathbb{W}\pot{u}} {\mathbb W}\rpot{u}=H\otimes_\bbZ {\mathbb W}\rpot{u}$. We have by base change 
\begin{equation}\label{descent_reductive}
\uG\;=\; (\Res_{{\mathbb W}\rpot{u}/{\mathbb W}\rpot{t}}\uH)^\sig,
\end{equation}
and since ${\mathbb W}\rpot{u}/{\mathbb W}\rpot{t}$ is \'etale the latter is a reductive group scheme over ${\mathbb W}\rpot{t}$ (with connected fibers).

\begin{lem} \label{PZ_group_scheme}
The ${\mathbb W}\pot{t}$-group scheme $\ucG_\bbf$ is a Bruhat-Tits group scheme for $\uG$ in the sense\footnote{As in \cite[Thm.~4.1]{PZ13}, the reductive group scheme $\uG$ is defined over ${\mathbb W}[t,t^{-1}]$. It is more convenient for us to consider the base change along ${\mathbb W}[t,t^{-1}]\to {\mathbb W}\rpot{t}$.} of \cite[Thm.~4.1]{PZ13}. In particular, it is smooth and affine with connected fibers. 
\end{lem}
\begin{proof} If $\uG=\uH$ is split, then by construction $\ucG_\bbf=\ucH_\bbf$ is such a Bruhat-Tits group scheme, cf.~\cite[\S4.2.2 (a)]{PZ13}. This is the first step in showing that \eqref{Res_Scalar_Grp_Witt} agrees with the construction in \cite[p.~180, middle]{PZ13} in general: starting from $\uG$, we may follow \cite[\S4.2]{PZ13} and construct the group scheme analogous to the one Pappas-Zhu denotes as $\calG_\Om=(\calG_\Om')^o$, where $\calG'_\Om$ is defined on the bottom of p.\,187. We observe the following: if we start from $\calG^\#_\bbf:=\ucG_\bbf\otimes_{{\mathbb W}\pot{t}, t\mapsto p} {\mathbb W}$, which is a parahoric group scheme in mixed characteristic, and apply the construction of Pappas-Zhu \cite[Thm.~4.1]{PZ13} to it, then we recover the group scheme $\ucG_\bbf$.  We use along the way that there is a canonical identification
\[
\calG^\#_\bbf\otimes_{\mathbb W}{\mathbb K}[p^{1\over e}] \;=\; H\otimes_\bbZ {\mathbb K}[p^{1\over e}]
\]
coming from \eqref{Res_Scalar_Grp_Witt} or \eqref{descent_reductive}, that is, under $\Gal(F'/F)=\Gal({\mathbb K}[p^{1\over e}]/{\mathbb K})$ the Galois actions on the group $\Res_{{\mathbb W}\pot{u}/{\mathbb W}\pot{t}}(\ucH_\bbf)$ induced from $\calG_\bbf$ resp.~$\calG^\#_\bbf$ agree. 
\end{proof}

We define the loop groups as the functor on the category of ${\mathbb W}$-algebras $R$ given by $L_{\mathbb W}\uG(R)=\uG(R\rpot{t})$ (resp.~$L_{\mathbb W}^+\ucG_\bbf(R) =\ucG_\bbf(R\pot{t})$), and $L_{\mathbb W}\uH(R)=\uH(R\rpot{u})$ (resp.~$L_{\mathbb W}^+\ucH_\bbf(R)=\ucH_\bbf(R\pot{u})$). This leads to the identifications 
\begin{equation}\label{Relation_Loop_Grps2}
L_{\mathbb W}\uG\,=\, (L_{\mathbb W}\uH)^{\sig} \;\;\; \text{and}\;\;\; L_{\mathbb W}^+\ucG_\bbf\,=\, (L_{\mathbb W}^+\ucH_\bbf)^{\sig,o}.
\end{equation}
For the second equality, we note that $L_{\mathbb W}^+(\Res_{{\mathbb W}\pot{u}/{\mathbb W}\pot{t}}\ucH_\bbf)^\sig\,=\, (L_{\mathbb W}^+\ucH_\bbf)^{\sig}$ which is a countably infinite successive extension of the $\sig$-fixed points of
\[
\Res_{{\mathbb W}\pot{u}/{\mathbb W}\pot{t}}\ucH_\bbf\otimes_{{\mathbb W}\pot{t},t\mapsto 0} {\mathbb W}
\]
by vector groups which only depend on a neighborhood of the unit section and so are the same for $L_{\mathbb W}^+\ucG_\bbf$, cf.~\cite[Prop.~A.4.9, (A.4.11)]{RS20}. Since vector groups are fiberwise connected, we obtain the desired equality using that taking fiberwise connected components commutes with base change, cf.~[SGA 3, $\on{VI}_B$, Prop.~3.3].  

\coro\label{special.fiber.groups.cor}
As group ind-schemes $L_{\mathbb W}\uG\otimes_{\mathbb W}k=LG$ compatible with the subgroup schemes $L_{\mathbb W}^+\ucG_\bbf\otimes_{\mathbb W}k=L^+\calG_\bbf$.
\xcoro
\pf This is immediate from \eqref{Relation_Loop_Grps2} and \eqref{Relation_Loop_Grps}.
\xpf

\subsection{Tame liftings of Schubert varieties} Being ind-schemes, the loop groups \eqref{Relation_Loop_Grps2} define fpqc (in particular \'etale) sheaves on the category $\AffSch/{\mathbb W}$ of affine ${\mathbb W}$-schemes.

\begin{lem}\label{etale_loc_triv}
The \'etale quotient $\uFl_{G,\bbf}:=L_{\mathbb W}\uG/L_{\mathbb W}^+\ucG_\bbf$ is an fpqc sheaf on $\AffSch/{\mathbb W}$ which is represented by an ind-projective ${\mathbb W}$-ind-scheme. There is an identification
\[
\uFl_{G,\bbf}\otimes_{\mathbb W} k=\Fl_{G,\bbf}.
\]
\end{lem}
\begin{proof} The proof is the same as in \cite[Lem.~5.3.2 (i)]{RS20}: Let $T'\to T$ be a faithfully flat map in $\AffSch/{\mathbb W}$. Let $T' \gets P' \to L_{\mathbb W}\uG$ be an object in $\uFl_{G,\bbf}(T')$ together with a descent datum along $T'\to T$. By effectivity of descent for affine schemes \cite[0244]{StaProj}, the torsor $T' \gets P'$ descends to a fpqc-locally trivial torsor $T\gets P$ represented by affine schemes. The map $P'\to L_{\mathbb W}\uG$ descends as well because every ind-scheme is an fpqc-sheaf. By \cite[Prop.~A.4.9, Exam.~A.4.12 iii.(a)]{RS20} every fpqc-locally trivial $L_{\mathbb W}^+\ucG_\bbf$-torsor is \'etale-locally trivial. Thus, $T \gets P \to L_{\mathbb W}\uG$ is an object of $\uFl_{G,\bbf}(T)$. Now the representability of $\uFl_{G,\bbf}$ is a special case of \cite[Prop.~6.5]{PZ13} in view of Lemma \ref{PZ_group_scheme}. Finally, the displayed formula is immediate from Corollary \ref{special.fiber.groups.cor} noting that sheafification commutes with fiber products. 
\end{proof}

We can also provide something close to a lift of Schubert varieties. 
First, it is well known that the Iwahori-Weyl group not only does not depend on $p$ but it admits an integral realization (see \cite[Prop. 3.4.1]{Lou19}). 
Indeed, we have a canonical isomorphism:
$$\uN({\mathbb W}(k)\rpot{t})/\uT({\mathbb W}(k)\pot{t}) \xrightarrow{\sim} N(k\rpot{t})/T(k\pot{t})=W $$
where $N$ is the normalizer of $S$ in $G$ and its underlined counterpart is its canonical lift to a closed subgroup of $\uG$. 
For any representative $\dot{w} \in \uN({\mathbb W}(k)\rpot{t})$ of $w \in W/W_\bbf$, we denote by $\uS_{w}=\uS_{w}(\bba,\bbf)$ the scheme-theoretic image of the morphism $L^+\ucG_{\bba} \rightarrow \uFl_{G,\bbf}, g \mapsto g\cdot \dot{w} \cdot e$.

\prop\label{schubert.var.int.prop}
For any $w\in W/W_\bbf$, the ${\mathbb W}(k)$-scheme $\uS_{w}=\uS_w(\bba,\bbf)$ is projective, integral, geometrically unibranch and flat over the base. 
The following are equivalent: the map $S_w \to \uS_w \otimes k$ is an isomorphism; $\uS_w \otimes k$ is reduced; $S_w$ is normal.
\xprop

The latter properties hold whenever $p \nmid \# \pi_1(G_\text{der})$ and only very rarely otherwise.

\pf
We may and do assume that $w \in W_{\text{aff}}$. Projectivity follows from the existence of Demazure resolutions (see \cite[Eq.~(9.18)]{PR08}), whereas being integral and flat over the base follows from the similar properties for the smooth finite type quotients of the positive loop group $L^+\ucG_{\bba}$. 
For the remaining claims, we consider the morphism
\begin{equation}\label{integral.sc.schubert.eq}
\uS_{\scon,w}= \uS_{\scon,w}(\bba,\bbf)\;\longto\; \uS_{w}(\bba,\bbf)=\uS_{w}.
\end{equation}
As for $S_{\scon, w} \to S_{w}$, this can be shown to be finite, birational, and a universal homeomorphism (see \cite[pf.~of Prop.~3.5]{HR22}). 
Moreover, we know that $\uS_{\scon,w}$ is geometrically normal over ${\mathbb W}(k)$ and that its special fiber is nothing more than $S_{\scon,w}$, by \cite[\S9]{PR08}. 
Hence $\uS_{w}$ is at least geometrically unibranch. 

Suppose now that $S_{w}$ is normal. 
The scheme $\uS_{w}$ is reduced, hence has a normalization, which can be identified with the canonical morphism $c\co \uS_{\scon, w} \to \uS_{w}$ from \eqref{integral.sc.schubert.eq}. 
The canonical closed immersion $S_{w} \to \uS_{w} \otimes k$ fits into a commutative diagram
$$
\xymatrix{
S_{\scon, w} \ar[r]^a \ar[d]_{\wr} & S_{w} \ar@{^(->}[d] \\
\uS_{\scon, w} \otimes k \ar[r]^b & \uS_{w} \otimes k.}
$$
The map $a$ is an isomorphism since $S_{w}$ is normal, hence $b$ is a closed immersion. 
The map $c$ is therefore fiberwise a closed immersion and a homeomorphism, hence is by Nakayama's lemma a closed immersion of reduced schemes. It follows that $c$ is an isomorphism, and then the diagram shows that $S_{w} =  \uS_{w} \otimes k$.

 On the other hand, if $\uS_{w}\otimes k =S_{w},$ equivalently, the special fiber $\uS_{w} \otimes k$ is reduced, then we have an equality
$$\dim_k H^0(S_{w}, \calL^N)= \dim_k H^0(S_{\scon, w}, \calL^N)$$
for any ample line bundle $\calL$ on $\uS_{w}$ and $N>0$ sufficiently large, by transporting the claim to the generic fiber using flatness. 
This implies that the map $S_{\scon, w} \to S_{w}$ is an isomorphism and thus $S_{w}$ is normal, so that we are done.
\xpf

\rema  \label{S_w_formation}  
Assume $S_w$ is normal, so that, as above,  $\uS_w$ is normal.  
In this case, one can show more generally that the formation of $\uS_w$ is compatible with arbitrary base change, in the following sense.   
Let $Z$ be any ${\mathbb W}(k)$-scheme. 
Then $\uS_{w,Z}:=\uS_w\times_{\Spec({\mathbb W}(k))}Z$ is equal to the scheme theoretic image of the map
\begin{equation}\label{Act_Map}
L_{\mathbb W}^+\ucG_{\bbf,Z} \,\longto\, \uFl_{G,\bbf, Z},\;\; b\mapsto b\cdot \dot{w}\cdot e_Z,
\end{equation}
where $e_Z$ denotes the base point in $\uFl_{G, {\bf f}, Z} : = \uFl_{G, \bf f} \times_{\Spec({\mathbb W}(k))}Z$. 
The main ingredient in the essential case of ${\bf f} = {\bf a}$ is the fact that, if $\underline{\pi}_w: \uD_{\tilde{w}} \to \uS_w$ is the Demazure resolution,  the formation of $\underline{\pi}_{w,*}\calO_{\underline{D}_{\tilde{w}}}$ commutes with arbitrary base change, cf.~\cite{Fal03}, \cite[Lem.~3.13, Prop.~3.15 ff.]{Gor03}.
\xrema

\subsection{Negative loop groups}\label{Neg_Loop_Grp_Sec} 
We continue with the same notation and assumptions. The base point ${\bf 0}\in \scrA(H, T_H ,F')$ defined by $H\otimes_\bbZ\calO_{F'}$ is invariant under the Galois group, 
and defines a special point also denoted ${\bf 0}\in \scrA(G,S,F)$ (because by construction $\sigma_0$ preserves the pinning $(G,B, T,X)$). After conjugation by an element in $W_\aff$, we may assume that the alcove $\bba$ contains ${\bf 0}$ in its closure, and lies in the chamber defined by the Borel $B_H$. 

We adapt the notion of the negative loop group from \cite[\S3.6]{dHL18} to our set-up as follows: the Iwahori $\calH_\bba$ corresponds now to the Borel subgroup $B_H\subset H$, more precisely, $\calH_\bba$ is the N\'eron blow up (resp.~dilatation) of $H\otimes_\bbZ\calO_{F'}$ in $B\otimes_\bbZ k$, cf.~\cite[Exam.~3.3]{MRR20}. 
We let $B^{\on{op}}_H=T_H\ltimes U^{\on{op}}_H$ denote the opposite Borel subgroup. The negative loop group is the functor on the category of ${\mathbb W}$-algebras $R$ given by $L^-_{\mathbb W}H(R):=H(R[u^{-1}])$. We define $L^{--}_{\mathbb W}H:=\ker(L^-_{\mathbb W}H\to H_{\mathbb W})$, $u^{-1}\mapsto 0$, and further we define {\em strictly negative loop group}
\[
L^{--}_{\mathbb W}\ucH_\bba := L^{--}_{\mathbb W}H\rtimes U^{\on{op}}_{H,{\mathbb W}},
\]
which is a subgroup ind-scheme of the ind-affine ind-scheme $L_{\mathbb W}\uH$ over ${\mathbb W}$. Finally, for the facet $\bbf$ contained in the closure of $\bba$, we define the {\em strictly negative loop group}
\begin{equation}\label{strictly.negative.split}
L^{--}_{\mathbb W}\ucH_\bbf \defined \bigcap_{w\in W_{H,\bbf}} {^w}\big(L^{--}_{\mathbb W}\ucH_\bba\big),
\end{equation}
where $W_{H,\bbf}$ denotes the subgroup of the affine Weyl group $W_{H,\aff}$ corresponding to the unique facet containing $\bbf\subset \scrA(H,T_H,F')$; see the beginning of Section \ref{normality.criterion}. As $H$ is split, each element $w\in W_{H,\aff}$ has a representative in $\dot{w}\in H({\mathbb W}\rpot{u})$. We set
\[
{^w}\big(L^{--}_{\mathbb W}\ucH_\bba\big):=\dot{w}\cdot \big(L^{--}_{\mathbb W}\ucH_\bba\big)\cdot \dot{w}^{-1}\;\subset \; L_{\mathbb W}\uH,
\]
and the intersection \eqref{strictly.negative.split} is taken inside $L_{\mathbb W}\uH$. The strictly negative loop group has the following key property.

\begin{lem}\label{Open_Lem_Split}
The map $L^{--}_{\mathbb W}\ucH_\bbf\to L_{\mathbb W}\uH/L_{\mathbb W}^+\ucH_\bbf$, $h^-\mapsto h^-\cdot e$ is representable by a quasi-compact open immersion.
\end{lem}
\begin{proof}
Equivalently, we have to show that the multiplication map 
\begin{eqnarray}
L^{--}_{\mathbb W}\ucH_\bbf\times L_{\mathbb W}^+\ucH_\bbf \to L_{\mathbb W}\uH
\end{eqnarray} 
is a quasi-compact open immersion (to check the equivalence we use that $L_{\mathbb W}^+\ucH_\bbf \to\Spec({\mathbb W})$ is fpqc, and that quasi-compact immersions are of effective fpqc descent \cite[02JR, 0246]{StaProj}).

This in turn is equivalent to being representable in schemes by a finitely presented \'etale monomorphism (see \cite[Thm.~17.9.1]{EGAIV4}). 
This was already known in the case $\bbf=0$ (see \cite[pf.~of Lem.~3.1]{HR21}) or working over a field (see \cite[pf.~of Prop.~3.7.1]{dHL18}).

The representability follows by writing the morphism as a composition of a closed immersion
$$L^{--}_{\mathbb W}\ucH_\bbf\times L_{\mathbb W}^+\ucH_\bbf \to L_{\mathbb W}\uH \times L_{\mathbb W}^+\ucH_\bbf$$
followed by the group multiplication $$L_{\mathbb W}\uH \times L_{\mathbb W}^+\ucH_\bbf \to L_{\mathbb W} \uH$$ which is representable, because the functor $L_{\mathbb W}^+\ucH_\bbf$ also is. For finite presentedness, we simply observe that both $L^{--}_{\mathbb W}\ucH_\bbf$ and $L_{\mathbb W}\uH/L^+_{\mathbb W}\ucH_\bbf$ are of ind-finite type.

Next we show that the map is a monomorphism, that is, that the finite type ${\mathbb W}$-group subscheme
$$L_{\mathbb W}^{--}\ucH_{\bbf} \cap L_{\mathbb W}^{+}\ucH_{\bbf} $$
of $L_{\mathbb W}H$ equals its unit section. We can do this in two different ways: either check it on both fibers, see \cite[Prop.~3.7.1]{dHL18}, which implies that the defining ideal is $p$-divisible and $p$-power torsion, hence trivial; or we check that its field valued points are trivial, again by \cite[Prop.~3.7.1]{dHL18}, and that its Lie algebra with coefficients in any ${\mathbb W}$-algebra $R$ vanishes.

Actually, we are going to show more generally that we have a triangular decomposition $$\text{Lie}\,L_{\mathbb W}^{--}\ucH_\bbf \oplus \text{Lie}\,L_{\mathbb W}^+\ucH_\bbf= \text{Lie}\,L_{\mathbb W}\uH.$$
This can be easily calculated using our choice of pinning (especially if $\bbf=\bba$ or $0$); comp.~\cite[Prop.~3.6.4]{dHL18}. 
Alternatively, we can observe that we have equalities
\begin{eqnarray}
\text{Lie}\,L_{\mathbb W}^+\ucH_\bbf= \sum_{w\in W_{H,\bbf}} {^w}\big(\text{Lie}\,L^{+}_{\mathbb W}\ucH_\bba\big) \\
\text{Lie}\,L_{\mathbb W}^{--}\ucH_\bbf= \bigcap_{w\in W_{H,\bbf}} {^w}\big(\text{Lie}\,L^{--}_{\mathbb W}\ucH_\bba\big)
\end{eqnarray}
almost by definition. 
This reduces the decomposition to the case $\bbf = \bba$, where it is clear. 
Indeed, for the purposes of reducing general facets to alcoves, we can further decompose into $T_H$-weights and the intersections of $\text{Lie}\,L_{\mathbb{W}}U_a$ for any root subgroup $U_a$ with a $W_{H,\bbf}$-conjugate of $\text{Lie}\,L^{+}_{\mathbb W}\ucH_\bba$ (respectively, $\text{Lie}\,L^{--}_{\mathbb W}\ucH_\bba$) are totally ordered by inclusion, so their sum (respectively, intersection) equals one of them (respectively, its complement), comp.\,\cite[Prop.\,3.6.4]{dHL18}.
Note that formation of the Lie algebra commutes with (completed) base change for all functors under consideration. 
Finally, it is enough to remark that this decomposition implies \'etaleness of the original map, via translating back to the origin (here we must use that all functors are formally smooth).
\end{proof}

\rema
We could have also argued via a bundle interpretation as in \cite[Lem.~3.1]{HR21} by constructing an appropriate opposite parahoric group scheme over ${\mathbb W}[t^{-1}]$. This is done in \cite[Def. 4.2.8, Cor. 4.2.11]{Lou19}.
\xrema

We now want to descend the result.

\begin{lem} 
The subgroup ind-scheme $L_{\mathbb W}^{--}\ucH_\bbf\subset L_{\mathbb W}\ucH_\bbf$ is $\sig$-invariant.
\end{lem}
\begin{proof} As the base point $0$ is $\sig$-invariant, one finds that the subgroup $L^{--}_{\mathbb W}H$ is $\sig$-invariant. The $\sig$-invariance of the alcove $\bba\subset \scrA(H,T_H,F')$ implies that the opposite unipotent radical $U_{H,{\mathbb W}}^{\on{op}}$ is $\sig$-invariant. Note that $\sig$ acts on $U_{H,{\mathbb W}}^{\on{op}}$ through the automorphism ${\sig_0}$. The lemma follows from the definition \eqref{strictly.negative.split} using the $\sigma$-invariance of $\bbf$.
\end{proof}

We define the {\em twisted strictly negative loop group} as 
\begin{equation}\label{Neg_Loop_Grp_Twisted}
L^{--}_{\mathbb W}\ucG_\bbf\defined (L_{\mathbb W}^{--}\ucH_\bbf)^{\sig,o}
\end{equation}

\begin{cor} \label{Open_Cor}
The map $L^{--}_{\mathbb W}\ucG_\bbf\to L_{\mathbb W}\uG/L_{\mathbb W}^+\ucG_\bbf=\uFl_{G,\bbf}$, $g^-\mapsto g^-\cdot e$ is representable by a quasi-compact open immersion.
\end{cor}
\begin{proof} 
As in the proof of Lemma \ref{Open_Lem_Split}, it is enough to show that the multiplication map $L_{\mathbb W}^{--}\ucG_\bbf\times_{\mathbb W} L_{\mathbb W}^+\ucG_\bbf\to L_{\mathbb W}\uG$ is a quasi-compact open immersion.
There is a Cartesian diagram
\[
\begin{tikzpicture}[baseline=(current  bounding  box.center)]
\matrix(a)[matrix of math nodes, 
row sep=1.5em, column sep=2em, 
text height=1.5ex, text depth=0.45ex] 
{L_{\mathbb W}^{--}\ucH_\bbf\times_{\mathbb W} L_{\mathbb W}^+\ucH_\bbf& L_{\mathbb W}\uH  \\ 
(L_{\mathbb W}^{--}\ucH_\bbf)^\sig\times_{\mathbb W} (L_{\mathbb W}^+\ucH_\bbf)^\sig& (L_{\mathbb W}\uH)^\sig, \\}; 
\path[->](a-1-1) edge node[above] {}  (a-1-2);
\path[->](a-2-1) edge node[left] {}  (a-1-1);
\path[->](a-2-1) edge node[below] {}  (a-2-2);
\path[->](a-2-2) edge node[right] {}  (a-1-2);
\end{tikzpicture}
\]
where the horizontal maps are given by multiplication, and the vertical maps are the canonical closed immersions. As the top arrow is an open immersion by Lemma \ref{Open_Lem_Split}, the bottom arrow is an open immersion as well. The corollary now follows from \eqref{Relation_Loop_Grps2}, \eqref{Neg_Loop_Grp_Twisted} by passing to neutral components. 
\end{proof}


\section{Kac-Moody flag varieties and projective embeddings}\label{Kac-Moody-Interlude}

In this subsection, we aim at generalizing Ramanathan's methods \cite{Ram87} via Frobenius splitting to gather information on the homogeneous ideals that define Kac-Moody Schubert varieties inside projective spaces or their Schubert overvarieties. 
We will also follow the treatment of Mathieu and use some ideas of \cite{Mat89}. 
An important result for us is Corollary \ref{tangent.space.formula.and.independence.of.characteristic} which gives a formula for the tangent spaces of Schubert varieties in arbitrary characteristic.
All Kac-Moody algebras below are assumed to be symmetrizable.

Let us start by recalling the definition of a (symmetrizable) Kac-Moody algebra. These are (mostly infinite-dimensional) Lie algebras $\frakg$ over $\bbC$ associated to symmetrizable generalized Cartan matrices, i.e. finite integer-valued square matrices $A=(a_{ij})$ with $a_{ii}=2$ and $a_{ij}\leq 0$, $i\neq j$, which become symmetric after multiplication by an invertible diagonal matrix, see \cite[\S1.1]{Kac90}. To that end, one starts with the notion of a realization $(\frakh, \Pi, \Pi^{\vee})$ of the given generalized Cartan matrix, consisting of a finite dimensional $\bbC$-vector space $\frakh$, a linearly independent set of roots $\alpha_i \in \frakh^{\vee}$, $i=1, \dots,n$ and coroots $h_i:=\alpha_i^{\vee}\in \frakh$ such that $\langle \alpha_i^\vee, \alpha_j \rangle =a_{ij}$ and $\text{dim}\,\frakh=n+\text{corank}A$, see \cite[\S1.1]{Kac90}. Then, we extend the abelian Lie algebra $\frakh$ to a Kac-Moody algebra $\frakg$ by freely adding generators $e_i$, $f_i$ for each positive simple root $\alpha_i$ and then by imposing the relations: $[h,e_i]=\alpha_i(h)e_i$, $[h,f_i]=-\alpha_i(h)f_i$, $[e_i,f_j]=\delta_{ij}h_i$, $\text{ad}^{-a_{ij}+1}(e_i)(e_j)=0$ and similarly for $f_i$, $f_j$, cf.~\cite[p.~16--17]{Mat88}. 

We have root and coroot lattices $Q=\sum_{i=1}^n\bbZ \alpha_i \subseteq \frakh^\vee$, $Q^{\vee}= \sum_{i=1}^n\bbZ \alpha_i^\vee\subseteq \frakh$. 
It turns out that $\frakg$ factors into a sum of weight spaces 
$$\frakg=\bigoplus_{\alpha \in \Phi \cup 0 } \frakg_{\alpha} $$
for the adjoint action of $\frakh$. Here $\Phi \subset Q$ denotes the root system of $\frakg$, for which the $\alpha_i$ form a basis (see \cite[Thm.~1.2]{Kac90} for these assertions). It admits a natural partition into real roots, that is, those that are conjugate to a positive simple root under the Weyl group, and imaginary roots, that is, the rest of them. 
Note that there is a triangular decomposition $\frakg=\frakn_+ \oplus \frakh \oplus \frakn_-$, where $\frakn_+$, resp.~$\frakn_-$ denotes the sum of all positive resp.~negative weight spaces; we denote by $\frakb_+$ the positive Borel subalgebra. 
Finally, we fix once and for all a weight lattice $P \subseteq \frakh^\vee$ and a coweight lattice $P^\vee \subseteq \frakh$ given by taking $\bbZ$-duals, such that there are inclusions of abelian groups $Q \subset P$ and $Q^\vee \subset P^\vee$ such that the latter is saturated (that is, with flat cokernel), compare also with \cite[p.~16]{Mat88}. 

Consider the category $\calO$ of finitely generated $\frakg$-modules which decompose into $\frakh$-weight spaces and whose finitely generated $\frakn_+$-submodules have finite dimension as vectors spaces. We have a highest weight module $V(\lambda)$ with maximal dominant weight $\lambda \in P$ for the Bruhat order. 
This arises as the unique irreducible quotient of the universal Verma module $\frakU(\frakg)\otimes_{\frakU(\frakb_+)}\bbC_\lambda$. The extremal weights of $V(\lambda)$ are the conjugates $w\lambda$ of the highest weight under the Weyl group action and they have multiplicity $1$. 
Demazure modules are the cyclic $\frakb_+$-modules generated by $V(\lambda)_{w\lambda}$.

In order to study arithmetic related to Kac-Moody algebras and groups, Tits introduced a $\bbZ$-form $\frakU_\bbZ(\frakg)$ of the universal enveloping algebra and a fortiori a $\bbZ$-form $\frakg_\bbZ$ of the Lie algebra.
In \cite{Mat88} and \cite{Mat89}, Mathieu uses this to define a certain ind-affine $\bbZ$-group ind-scheme $\frakG$ whose Hopf algebra of distributions supported at the origin (also known as the hyperalgebra) is given by the completion $\widehat{\frakU}_\bbZ(\frakg)$ of $\frakU_\bbZ(\frakg)$ for the obvious descending filtration (compare with \cite[Lem.~2, Lem.~3]{Mat89}). 
It comes equipped with a canonical maximal split torus $\frakT$ corresponding to $\frakh_\bbZ$, as well as a positive Borel subgroup $\frakB^+=\frakT\ltimes \frakU^+$ containing it. 
Let us mention that $\frakB^+$ is an affine, non-finitely presented, flat, closed subgroup scheme of $\frakG$ with underlying hyperalgebra given by $\widehat{\frakU}_\bbZ(\frakb^+) $. 
The fppf quotient $\frakF:=\frakG/\frakB^+$ is representable by a reduced ind-projective $\bbZ$-ind-scheme. 
It is known as Mathieu's flag variety associated to the Kac-Moody algebra $\frakg$ (together with the rest of the chosen data). Given an admissible set $J \subseteq I$ of positive simple roots (i.e.~such that the subgroup $W_J$ generated by the corresponding reflections is finite), we  can associate to it a standard parabolic subgroup $\frakP^+_J=\frakL_J \ltimes \frakU_J^+$ containing $ \frakB^+$ as a closed subgroup scheme. 
We have a partial flag variety $\frakF_J:=\frakG/\frakP^+_J$, which is still representable by an ind-projective $\bbZ$-ind-scheme. 
For each $w \in W/W_J$, we may consider the Schubert variety $ \frakS_{w,J} \subseteq \frakF_J$ obtained as the scheme-theoretic image of the orbit map $\frakB^+ \rightarrow \frakF_J, b\mapsto b\cdot \dot{w}\cdot e$, where $e\in \frakF_{J}(\bbZ)$ is the base point and $\dot{w} \in \frakG(\bbZ)$ some representative of the class $w$.
It is a fundamental theorem of Mathieu \cite{Mat88} and Littelmann \cite{Lit98} that the $ \frakS_{w,J}$ are geometrically normal over $\bbZ$, i.e., the structural map $S_w\to \Spec(\bbZ)$ is a normal morphism in the sense of \cite[038Z]{StaProj}.

Next we are going to introduce the negative (resp.,\,strictly negative) parabolic subgroups $\frakP_J^-$ (resp.\,$\frakU_J^-$), seemingly a novelty in the literature: in \cite[p.~45]{Mat89}, Mathieu mentions that $\frakP^-_J$ is not defined. We will define $\frakP_J^-$ and $\frakU_J^-$ using $\mathbb G_{\on{m}}$-actions on $\frak G$. We refer the reader to \cite{Ric19a} or \cite[\S2.1]{CGP15} for basic facts on $\bbG_{\on{m}}$-actions. We will freely use the notation $X^+$ (resp.~$X^-$, resp.~$X^0$) as in \cite{Ric19a} (and \cite[Section 2]{HR20b} for ind-schemes) to denote the attractor (resp.~repeller, resp.~fixed point) ind-scheme of an ind-scheme $X$ over $\bbZ$ equipped with $\bbG_{\on{m}}$-action. We use similar notation applied to functors which are not known to be ind-schemes. (We emphasize that the superscripts $+,-$ already appear in connection with subgroups generated by positive or negative affine root groups. We hope that this will not cause confusion: we will eventually show that the different superscript meanings are compatible with each other.)

 Let us define $\frakP_J^-$ and $\frakU_J^-$. Fix any dominant $\mu\co \bbG_{\on{m},\bbZ} \to \frakT$, which is $J$-regular, i.e.,\,whose composition with the positive simple roots $\alpha$ in $J$ (resp. $I\setminus J$) is zero (resp. strictly positive).
As a functor, $\frakP_J^-$ is given by the repeller locus of $\mu$ acting via conjugation on $\frak G$, whereas its unipotent radical $\frakU_J^-$ is the strict repeller, i.e., the fiber over the identity section of the map $\frakP_J^-\to \frakG^0$ given by evaluation at $t=\infty$, where $\frakG^0$ denotes the sub-functor of fixed points for the $\bbG_{\on{m}}$-action. Note that  $\frak G^0 \subset \frakP_J^-$ and we have a semi-direct product decomposition $\frakP_J^- = \frakG^0 \ltimes \frakU_J^-$ of group functors.

\begin{lem}\label{negative.unipotent.lem}
The negative unipotent subgroup $\frakU_J^-$ is representable by an ind-affine closed subgroup ind-scheme of $\frakG$ of ind-finite presentation, which does not depend on the choice of the dominant, $J$-regular cocharacter $\mu$. 
The multiplication morphism $\frakU_J^- \times \frakP_J^+ \rightarrow \frakG$ is a quasi-compact open immersion.
\end{lem} 

\pf

We consider the $\bbG_{\on{m}}$-action on $\frakF_J$ induced by $\mu$. 
The action is Zariski locally linearizable in the sense of \cite[Thm.~1.8]{Ric19a} and \cite[Thm.~2.1]{HR21} for ind-schemes (implying representability of attractors, repellers and fixed points). We consider the ind-closed embedding $\frakF_J \to \bbP(V(\lambda)_\bbZ)$, where $\lambda$ is an integral dominant $J$-regular weight; comp.\,Theorem \ref{embeddings-given-by-quadrics-and-lines} below.
 The embedding is $\bbG_{\on{m}}$-equivariant when equipping $\bbP(V(\lambda)_\bbZ)$ with the linear action induced by the adjoint action of $\bbG_{\on{m}}$ via $\mu$ on the Weyl module $V(\lambda)_\bbZ$ given as the image of $ \frakU_\bbZ(\frakg)\otimes_{\frakU_\bbZ(\frakb_+)}\bbZ_\lambda$ in $V(\lambda)$ (for suitable integral structures).
 So any ind-affine $\bbG_{\on{m}}$-stable open cover of $\bbP(V(\lambda)_\bbZ)$ yields one of $\frakF_J$.
 Moreover, we observe that the identity section $e\in \frakF_J(\bbZ)$ is scheme-theoretically an isolated fixed point because the same is true in $\bbP(V(\lambda)_\bbZ)$, as is seen by considering the action of the dominant, $J$-regular cocharacter $\mu$ on the line spanned by the highest weight space in $V(\lambda)_\bbZ$. 
 So the fiber $\frakV$ of $\frakF_J^- \to \frakF_J^0$ over the identity section $e$ is equal to the open non-vanishing locus of the section $v_\la^\vee\in \Gamma(\frakF_J,\frakL)$ killing all the non-highest weight spaces, where $\frakL$ is the pullback of the line bundle $\calO(1)$ on $\bbP(V(\lambda)_\bbZ)$;  comp.\,also Corollary \ref{tangent.space.formula.and.independence.of.characteristic} below.

	We return to the ind-affine group scheme $\frakG$. 
	It has a $\bbG_{\on{m}}$-stable presentation by affine schemes $\frakG_{w,J}$, which are $\frakP_J^+$-torsors over $\frakS_{w,J}$, so \cite[Lem.~1.9]{Ric19a} applies to show that $\frakU_J^-\subset \frakP_J^- \subset \frakG$ are representable by closed immersions. 
	We have $\frakG^0\supset \frakL_J$ and we claim that $\frakG^0=\frakL_J$. 
	We verify the equality first at the level of $k$-valued points, where $k$ is any field: by the Bruhat decomposition, any $g \in \frakG^0(k) $ lies in some double coset $\frakP^+_J(k)\dot{w}\frakP_J^+(k)$ with the lift $\dot{w}$ normalizing $\frakT$.  As $\frakP_J^+ \subset \frakG^+$ we deduce that $\dot{w} \in \frakG^+(k)$. 
	The resulting morphism 
	\begin{equation}\label{conjugation_weyl_elts}
	\mu \dot{w}\mu^{-1} =(\mu\,^w\mu^{-1})\dot{w} \colon \bbG_{\on{m},k} \to \frakG
	\end{equation}
	takes values in the closed subscheme $\frakT \dot{w}$ and extends to $\bbA_k^1 = \bbG_{\on{m},k} \cup \{0\}$ because $\dot{w} \in \frakG^+$.  On the other hand, $\mu\,^w\mu^{-1}$ extends to $\bbA_k^1$ exactly when it is trivial, which happens when $w\in W_J$ by regularity of $\mu$. This shows that $g \in \frakP_J^+(k) = \frakL_J(k) \ltimes \frakU^+_J(k)$, which in turn implies $g \in \frakL_J(k)$, since the $\frakU^+_J(k)$-factor contracts to the identity.
	This shows $\frakG^0(k)=\frakL_J(k)$ as desired.
	In particular, the natural map $\frakG^0\to \frakF_J$ factors through $\frakV$. 
	Since the scheme-theoretic fixed points of the $\bbG_{\on{m}}$-stable open neighborhood $ \frakV$ of $e \in \frakF_J$ coincide with the origin section $e$, it follows that $\frakG^0 \subset \frakP^+_J$, the fiber of $\frakG\to \frakF_J$ above $e$.
	But $\frakP^+_J=\frakL_J \ltimes \frakU_J^+$ and the limit when $t=0$ of the attractive $\frakU_J^+$ is the identity, so this implies $\frakG^0=\frakL_J$, and thus $\frakP_J^-=\frakL_J\ltimes \frakU_J^-$.
	
	Before proceeding, we calculate the $k$-valued points of $\frakU_J^-$ for any field $k$. 
	Set $G=\frakG(k)$, $P_J^\pm=\frakP_J^\pm(k)$, $L_J=\frakL_J(k)$, $U_J^\pm=\frakU_J^{\pm}(k)$, $U_\alpha=\frakU_\alpha(k)$ for real roots $\alpha$. 
	Then the Birkhoff decomposition \cite[Prop.~3.16]{Rou16} reads
	\begin{equation} \label{Birkhoff_J_new}
	G=\bigsqcup_{w \in W_J\setminus W/W_J} Q_J^- \dot{w} P_J^+,
	\end{equation}
	where $Q_J^-\subset P_J^-$ is the subgroup generated by $L_J$ and the $U_\alpha$ with $\alpha$ being $J$-negative. 
	(Here a $J$-positive real root is a positive real root in which at least one simple root $\alpha_i$ in $I\backslash J$ appears with positive multiplicity; a $J$-negative real root is the negative of a $J$-positive real root.)
	We claim $P_J^-=Q_J^-$.
	Indeed, if $P_J^-$ were different from $Q_J^-$, then it would intersect $\dot{w}P_J^+$ with $\dot{w} \notin L_J$. 
	But this is impossible as $\frakP_J^-$ maps to the contracting open $\frakV$, whose unique fixed point is the identity. 
	In particular, $U_J^-$ is generated by the $U_\alpha$ for all $J$-negative real roots $\alpha$. 
	Also, using the existence of Zariski local sections of $\frakG \to \frakF_J$, see \cite[Lemme 7]{Mat89}, we get that $\frakF_J(k)=G/P_J^+$, and hence \eqref{Birkhoff_J_new} implies that $U_J^-$ maps bijectively to $\frakV(k)$. Indeed, an element of the form $y \dot{w} e \in G/P^+_J$ with $y = ul \in P_J^- = U^-_J L_J$ can only contract to $e$ under the $\bbG_{\on{m}}$-action if $w \in W_J$, and then $y\dot{w}e = ue$.

	Finally, in order to show that $\frakU_J^-\times \frakP_J^+\to \frakG$ is a quasi-compact open immersion, it suffices to see that $\frakU_J^-\to \frakV$ is an isomorphism. Since we already know from the paragraph above that it is a universally bijective  monomorphism, it is enough to show that $\frakU_J^-\to \frakV$ is representable by a smooth morphism.
	Let $\frakP_J^+(n) \subset \frakP_J^+$ be the subgroup generated by the root groups associated with $J$-positive roots of height at least $n \in \bbN$. 
	The $\bbG_{\on{m}}$-equivariant map of ind-schemes $\frakG/\frakP_J^+(n)\to \frakF_J$ is affine (in fact a torsor for a finite-type group scheme over $\bbZ$, since $J$ is admissible), so the source has a Zariski locally linearizable $\bbG_{\on{m}}$-action. 
	We see using \eqref{Birkhoff_J_new} that the injection $P_J^-\to (\frakG/\frakU_J^+)^-(k)$ is bijective using the following steps. First, note that the pre-image of $\frakV(k) \cong U_J^-$ along $(\frakG/\frakU_J^+)(k)=G/U_J^+ \to \frakF_J(k)=G/P_J^+$ coincides exactly with $P_J^- = U_J^- L_J$. Second, we observe that $P_J^- \subseteq (\frakG/\frakU^+_J)^-(k)$ and the latter is stable under left multiplication by $P^-_J$.  Third, we notice $(\frakG/\frakU_J^+)^-(k)$ maps to $\frakV(k) \simeq U_J^-$, because otherwise \eqref{Birkhoff_J_new} implies that $(\frakG/\frakU_J^+)^-$ would have non-empty intersection with the $\bbG_{\on{m}}$-stable closed subscheme $\dot{w}\frakL_J$ of $\frakG/\frakU_J^+$ for $w \notin W_J$, whereas $\dot{w}\frakL_J$ has empty fixed points, because it does not meet $\frakL_J=\frakG^0$ when regarded as a subscheme of $\frakG$.

	Using the above, we see $(\frakG/\frakP_J^+(n))^-\to (\frakG/\frakU_J^+)^-$ is universally bijective because the fiber of $ \frakG/\frakP_J^+(n)\to \frakG/\frakU_J^+$ over the identity equals $\frakU_J^+/\frakP_J^+(n)$ and intersects trivially with the repeller; in other words, because $\frak P_J^- \to \frakP^-_J \frakU_J^+ /\frakP^+_J(n) \, \cap \, (\frakG/\frakP^+_J(n))^-$ is an isomorphism. By \cite[Lem.~2.2]{HR21}, the morphism $(\frakG/\frakP_J^+(n))^-\to (\frakG/\frakU_J^+)^-$ is also smooth, thus necessarily étale (as the geometric fibers are singletons), so it must be an isomorphism by \cite[Tag 025G]{StaProj}. 
	We claim that this yields an isomorphism $\frakP_J^- \cong (\frakG/\frakU_J^+)^-$ in the limit.
	For this, we pass to  $\bbG_{\on{m}}$-stable scheme presentations of the ind-schemes. 
	So we consider the closed subscheme $\frakG_w \subset \frakG$ obtained as pullback of $\frakS_w \subset \frakF_{J}$ along the projection $\frakG\to \frakF_{J}$. 
	We have isomorphisms
	\[ \frakG_w^- \cong \lim_{n\geq 1} (\frakG_w/\frakP_J^+(n))^- \cong (\frakG_w/\frakU_J^+)^-,  \]
	where the second follows from the previous identification for fixed $n$ (so all transition maps are isomorphisms), and the first from the fact that $\frakG_w\cong \lim_{n\geq 1}\frakG_w/\frakP^+_J(n)$ \cite[Lem.~3]{Mat89} and that attractors commute with cofiltered limits of schemes along affine maps. 
	Indeed, passing to an affine cover, the claimed commutation follows from the induced $\bbZ$-gradings commuting with filtered colimits of rings. 
	Taking now colimits for varying $w$, we get the claimed isomorphism $\frakP_J^- \cong (\frakG/\frakU_J^+)^-$, so the projection $\frakP_J^-\to \frakV$ is representable by a smooth surjection.
	This finishes the proof of the lemma.
	\xpf

For each $J$-regular dominant weight $\lambda$, we may consider the line bundle $\frakL(\lambda):=\frakG \times_{\frakP_J^+} \bbZ_{-\lambda }$ obtained from the natural $\frakP_J^+$-bundle $\frakG \rightarrow \frakF_J$ via the character $-\lambda$ of $\frakP_J^+$ (compare with \cite{Mat88} and \cite{Mat89}, which use the opposite sign convention). 
This is a very ample line bundle on $\frakF_J$ and we have a natural identification between $\Gamma(\frakF_J, \frakL(\lambda))^{\vee}:=\on{colim}_{w \in W/W_J}\Gamma(\frakS_{w,J}, \frakL(\lambda))^{\vee}$ and the canonical $\bbZ$-form $V(\lambda)_{\bbZ}$ of the highest weight module: indeed, at each finite step, the submodule $\Gamma(\frakS_{w,J}, \frakL(\lambda))^{\vee}$ is identified with the integral Demazure module $V_w(\lambda)_{\bbZ}$, see \cite[Thm.~5]{Mat88}. 
Finally, we recall that these constructions exhaust after linearization the Picard group of $\frakF_J$, by virtue of the isomorphism $P_J/P_I\simeq \text{Pic}(\frakF_J)$, where $P_J=\{ \lambda \in P\,|\, \lambda(\alpha_j^\vee)=0 \, \forall j \in J\}$ (see \cite[Prop.~28]{Mat88}).

The following result goes back originally to \cite{Ram87} in the case of classical flag varieties, by using the relatively new at the time method of Frobenius splitting. 
Credit is also due to Mathieu for showing the existence of an ind-splitting of the diagonal of the flag variety (see \cite[Prop.~1]{Mat89}). 
We are pretty convinced that Littelmann's path model also yields the same type of results, see~\cite{Lit98}.

\begin{thm}[Ramanathan, Mathieu]\label{embeddings-given-by-quadrics-and-lines}
Given an ample line bundle $\frakL$ of $\frakS_w$, the corresponding morphism $\frakS_w \rightarrow \bbP(\Gamma(\frakS_w, \frakL^{\vee}))$ is a closed immersion defined by quadrics. 
Moreover, the closed immersion $\frakS_u\rightarrow \frakS_w$, $u \leq w$ is linearly defined with respect to $\frakL$.
\end{thm}
 Here we allow $w=\infty$ to get the entire partial affine flag variety $\frakS_{\infty}=\frakF_J$.
\begin{proof}
The statement refers to the behavior of the graded algebra $\Gamma(\frakS_w, \frakL^{\bullet})$ and its graded module $\Gamma(\frakS_u, \frakL^{\bullet})$ as in Definition \ref{quadratic-graded-algebras-and-modules}. 
By upper semicontinuity, it suffices to base change to any positive characteristic $p$ field $k$. 
Since every ample line bundle on $\frakS_w$ extends to $\frakS_{\infty}$, we are reduced to showing, by Proposition \ref{ind-split-diagonal} and Proposition \ref{ind-split-triple-diagonal}, that the compatibly ind-split $\frakS_u \subseteq \frakS_w \subseteq \frakS_{\infty}$ satisfy: the diagonal $\Delta_{\frakS_{\infty}}$ is compatibly ind-split with $\frakS_{\infty}^2$; the partial mixed diagonals $\Delta_{\frakS_{\infty}}\times \frakS_{\infty}$, $\frakS_{\infty}\times\Delta_{\frakS_{\infty}}$, $\frakS_{w}\times\Delta_{\frakS_{\infty}}$ and $\frakS_{u}\times\Delta_{\frakS_{\infty}}$ are simultaneously compatibly ind-split with $\frakS_{\infty}^3$. 
We may and do assume $J=\emptyset$ by pushing forward the splitting along the obvious projection.

 For this, we need the convoluted flag varieties $\frakS_{\infty}^{\tilde{\times}n}:=\frakS_{\infty} \tilde{\times}\dots \tilde{\times} \frakS_{\infty}=\frakG\times^{\frakB_+}\dots \times^{\frakB_+} \frakF$. 
 Note that herein we have convoluted Schubert varieties $\frakS_{w_1,\dots, w_n}:=\frakS_{w_1} \tilde{\times} \dots \tilde{\times} \frakS_{w_n} $ which are all compatibly Frobenius split, as one can observe by using appropriate Demazure resolutions (see \cite[Lem.~9, Lem.~10]{Mat89}). In particular, under the natural isomorphism $\frakS_{\infty}^{\tilde{\times}n} \cong \frakS^n_{\infty}$ given by $(m_1, \dots, m_n)$, where $m_i$ denotes the product of the first $i$ coordinates, the diagonal $\Delta_{\frakS_{\infty}}$ is identified with $\frakS_{\infty,1}$ (compare with \cite[Prop.~1]{Mat89}), and the partial mixed diagonals $\Delta_{\frakS_{\infty}}\times \frakS_{\infty}$, resp.~$\frakS_{\infty}\times\Delta_{\frakS_{\infty}}$, resp.~$\frakS_{w}\times\Delta_{\frakS_{\infty}}$, resp.~$\frakS_{u}\times\Delta_{\frakS_{\infty}}$ are identified with $\frakS_{\infty,1,\infty} $, resp. $\frakS_{\infty,\infty,1}$, resp.~$\frakS_{w,\infty,1}$, resp.~$\frakS_{u,\infty,1}$.
\end{proof}

The following corollary gives an explicit formula describing the tangent space, which goes back to work of Kumar \cite{Kum96} in characteristic 0, and Polo \cite{Pol94} for classical flag varieties.

\begin{cor}[Kumar, Polo]\label{tangent.space.formula.and.independence.of.characteristic}
Let $k$ be a field of arbitrary characteristic. 
Let $\lambda$ be any fixed $J$-regular dominant weight, with associated ample line bundle $\frakL := \frakL(\lambda)$.
The $k$-vector space $T_e\frakS_{w,J}\otimes k$ consists of all $X \in T_e\frakF_{J}\otimes k$ such that $X v_{\lambda} \in V_w(\lambda)_k$.
\end{cor}

\pf
We start by noticing that $\frakU_J^-$ becomes naturally identified with the distinguished open subset $D_+(v_{\lambda}^\vee)\subseteq \frakF_J$ associated with the dual section $v_{\lambda}^{\vee} \in \Gamma(\Fl_J, \frakL)$ killing all weight spaces different from $V(\lambda)_\lambda=kv_\lambda$. 
Indeed, this is a consequence of \cite[Lem.~8.3]{Kum96}.

By Theorem \ref{embeddings-given-by-quadrics-and-lines}, the closed immersion $\frakU_{w,J}^-:=\frakS_{w,J}\cap \frakU_J^- \rightarrow \frakU_J^-$ is defined by the coefficients $\varphi_{\xi}$, where $\xi \in \Gamma(\frakF_J, \frakL)\otimes k$ runs over all vectors perpendicular to $V_w(\lambda)_k$ and  $\varphi_{\xi}(u)=\xi(uv_{\lambda})$ (see also \cite[Prop.~3.1]{Pol94}). 
Given a tangent vector $X \in T_e\frakF_J\otimes k=T_e\frakU_J^-\otimes k$, one can check that it lies in $T_e\frakS_{w,J} \otimes k$ if and only if the associated distribution of $k[ \frakU_J^-]:=\varprojlim k[\frakU_{w',J}^-]$ kills all the $\varphi_{\xi}$ designated above. 
Representing $X$ by a $k[\varepsilon]$-valued point $u$ of $\frakU_J^{-}$, we get 
$$
\varphi_{\xi}(u)=\xi(uv_{\lambda})=\varepsilon\,\xi(Xv_{\lambda}).
$$ 
The right side is obviously zero for all $\xi$ if $Xv_{\lambda}\in V_w(\lambda)_k$ and the converse follows from the isomorphism $V_w(\lambda)_k^\vee=\Gamma(\frakF_J, \frakL)\otimes k/V_w(\lambda)_k^\perp$. 
This proves our description of the tangent space (compare our argument with Polo's \cite[Thm.~3.2]{Pol94}).
\xpf

\rema\label{polo.remark}
Polo claims in \cite[Cor.~4.1]{Pol94} that the dimension of $T_e\frakS_{w,J}\otimes k$ does not depend on $p=\on{char}k$. 
He invokes the fact that $T_e\frakU_J^-$ has a natural $\bbZ$-model given by $\frakn_{J,\bbZ}^-$, whereas the integral model $V_w(\lambda)_{\bbZ}$ of $V_w(\lambda)_{k}$ is a direct summand of the model $V(\lambda)_{\bbZ}$ for $V(\lambda)_{k}$.
Now the independence of $p$ of the tangent space dimension is equivalent to the flatness of the cokernel of $\frakn_{J,\bbZ}^-\to V(\lambda)_{\bbZ}/V_w(\lambda)_{\bbZ}$, $X\mapsto Xv_\lambda$, or equivalently the saturatedness of its image. 
An argument for this is missing in the proof of \cite[Cor.~4.1]{Pol94}.
It would be interesting to clarify this point. 
This result has at least been invoked in \cite[Rem.~8.10]{Kum96} by Kumar in order to generalize his smoothness criterion to positive characteristic. 
\xrema

Assume $\frakg$ is an affine Kac-Moody algebra, that is, the corank of the corresponding generalized Cartan matrix is equal to $1$. 
These are classified by affine Dynkin diagrams and admit very explicit realizations as some mildly modified loop algebras or their fixpoints under order $2$ or $3$ automorphisms, see \cite[\S7--8]{Kac90}. 
More explicitly, the quotient of $[\frakg_{ \bbC}, \frakg_{\bbC}]$ by its non-trivial center (this is a phenomenon particular to infinite-dimensional Kac--Moody algebras) can be identified with the graded Lie algebra of the group scheme $L_{\mathbb W}\underline{G}\otimes \bbC$ constructed in \eqref{descent_reductive} for a given embedding ${\mathbb K}=\on{Frac}{\mathbb W}\hookrightarrow \bbC$, where $G$ is the only simply connected absolutely almost simple semisimple group over $k\rpot{t}$ having the same affine Dynkin diagram as $\frakg$. Under this correspondence, the Borel subalgebra $\frakb$ is mapped to the Lie algebra of $\ucG_\bba $ and every standard parabolic $\frakp_J$ to the Lie algebra of $\ucG_{\bbf_J} $ for some facet $\bbf_J$ in the boundary of $\bba$. 
We also pick the usual weight lattice $P$ in the Cartan subalgebra $\frakh$. 
We have the following important comparison result, which can be found in \cite[\S9.h]{PR08} in a weaker form (see also \cite[\S2.5]{Zhu17} for an exposition).

\begin{propo}\label{kac.moody.comp}
Let $k$ be an algebraically closed field of positive characteristic with ring of Witt vectors ${\mathbb W}={\mathbb W}(k)$.
There are natural isomorphisms 
\[
\frakP^+_J\otimes {\mathbb W} \cong \widehat{L^+_{\mathbb W}\ucG_{\bbf_J}}\rtimes \bbG_{\on{m}}^{\emph{rot}} \;\;\;\;\text{and}\;\;\;\;\frakG \otimes {\mathbb W} \cong \widehat{L_{\mathbb W}\uG}\rtimes \bbG_{\on{m}}^{\emph{rot}},
\]
inducing an equivariant isomorphism $\frakF_J \otimes {\mathbb W} \cong \uFl_{G,\bbf_J}$  compatible with $\frakU_J^-\otimes {\mathbb W}\cdot e\cong L^{--}_{\mathbb W}\ucG_{\bbf_J}\cdot e$.
\end{propo}

Here the hat loop groups are the central extensions of the respective loop groups by $\bbG_{\on{m}}$ given by parametrizing pairs $(g,\al)$ of group elements $g$ and isomorphisms $\al\co g^*\calO(1)\cong \calO(1)$ where $\calO(1)$ is the generator of $\text{Pic}(\uFl_{G,\bf 0})$, see \cite[Eqn. (4.3.29)]{Lou19}.
In fact, in place of a generator, we may use any line bundle of a partial affine flag variety with central charge $1$.
The rotation $\bbG_{\on{m}}$ is induced by automorphisms of the formal disk $R\pot{u}$ (as opposed to $R\pot{t}$).

\pf[Idea of proof]
First of all, we construct the isomorphism $\frakP^+_J\otimes {\mathbb W} \cong \widehat{L^+_{\mathbb W}\ucG_{\bbf_J}}\rtimes \bbG_{\on{m}}^{\text{rot}}$. 
This essentially amounts to verifing that their algebras of distributions match (see either \cite[\S8.d]{PR08} or \cite[Thm. A.2.4]{Lou19}).

Next we identify the Kac-Moody setting Demazure varieties $\frakD_{\tilde{w}}$ with those denoted $\uD_{\tilde{w}}$ above, in the natural way, induced by the twisted product decomposition in terms of parabolics for the left side and jet groups for the right side. 
By means of a topological argument, we can now get an equivariant identification $\frakS_w \cong \uS_{w}(\bba,\bba)$ (see \cite[Lem.~32, Lem.~33]{Mat88} and compare them to the references in the previous paragraph). 
This yields the desired equivariant identification of full flag varieties $\frakF \otimes {\mathbb W} \cong \uFl_{G,\bba}$, as both are geometrically reduced over ${\mathbb W}$, and then of the partial counterparts by taking quotients.

Now we deduce $\frakG \otimes {\mathbb W} \cong \widehat{L_{\mathbb W}\uG}\rtimes \bbG_{\on{m}}^{\rm rot}$. Consider the $\frakB^{+}$-bundle $\frakG_w$ over $\frakS_w$ obtained as the affine hull with respect to $\frakS_w$ of the canonical $\frakB^{+}$-bundle over the Demazure variety $\frakD_{\tilde{w}}$, see \cite[Ch.~XI]{Mat88}. Here, the affine hull of a morphism $f\colon Y \to X$ is the factorization $Y\to \Spec f_*\calO_Y \to X$, where the middle object is the relative spectrum of the quasi-coherent $\calO_X$-algebra  $f_*\calO_Y$.
By the universal property of the relative affine hull, we get a $\frakB^{+}\otimes_\bbZ \mathbb{W} \cong \widehat{L^+_{\mathbb W}\ucG_{\bba}}\rtimes \bbG_{\on{m}{\mathbb{W}}}^{\text{rot}}$-equivariant morphism towards the preimage of $\uS_{w}(\bba,\bba)$ in $\widehat{L_{\mathbb W}\uG}\rtimes \bbG_{\on{m}{\mathbb{W}}}^{\text{rot}}$, which then must be an isomorphism. 
Taking direct limits now recovers the isomorphism $\frakG \otimes {\mathbb W} \cong \widehat{L_{\mathbb W}\uG}\rtimes \bbG_{\on{m}{\mathbb{W}}}^{\text{rot}}$.

Finally let us prove $\frakU_J^-\otimes {\mathbb W}\cong L^{--}_{\mathbb W}\ucG_{\bbf_J}$ inside $L_{\mathbb W}\uG$. First we note that $\frakU_J^-\otimes {\mathbb W}$ lies in $\widehat{L_{\mathbb W}\uG}$ by naturality of strict repellers, because any $J$-regular dominant coweight $\mu$ has positive image in $\bbG_{\on{m}}^{\on{rot}}$. Recall from the proof of Lemma \ref{negative.unipotent.lem} that $\frakU_J^-(\kappa)$ is generated by the $J$-negative real root subgroups, and thus $\frakU_J^-(\kappa)$ maps into $L^{--}_{\mathbb W}\ucG_{\bbf_J}(\kappa)$ for any algebraically closed field $\kappa$ which is a ${\mathbb W}$-algebra.
Now the composition $\frakU_J^-\otimes {\mathbb W}\subset \widehat{L_{\mathbb W}\uG}\to L_{\mathbb W}\uG$ is a (representable) quasi-compact monomorphism of reduced ind-schemes.
Since $L^{--}_{\mathbb W}\ucG_{\bbf_J}$ is a closed sub-ind-scheme of $L_{\mathbb W}\uG$, we obtain a map $p\co \frakU_J^-\otimes {\mathbb W}\to L^{--}_{\mathbb W}\ucG_{\bbf_J}$ which is a quasi-compact monomorphism. 
Now consider the commutative diagram of ind-schemes
\[
\begin{tikzpicture}[baseline=(current  bounding  box.center)]
\matrix(a)[matrix of math nodes, 
row sep=1.7em, column sep=3em, 
text height=1.5ex, text depth=0.45ex] 
{\widehat{L_{\mathbb W}\uG}& \frakU_J^-\otimes {\mathbb W} & \frakF_J\otimes {\mathbb W}  \\ 
L_{\mathbb W}\uG & L^{--}_{\mathbb W}\ucG_{\bbf_J} & \uFl_{G,J} \\}; 
\path[->](a-1-2) edge node[above] {{\rm\tiny closed}}  (a-1-1);
\path[->](a-1-2) edge node[above] {{\rm\tiny open}}  (a-1-3);
\path[->](a-2-2) edge node[above] {{\rm\tiny closed}}  (a-2-1);
\path[->](a-2-2) edge node[above] {{\rm\tiny open}}  (a-2-3);
\path[->](a-1-1) edge node[left] {{\rm\tiny smooth}} node[right] {$q$}  (a-2-1);
\path[->](a-1-2) edge node[right] {$p$}  (a-2-2);
\path[->](a-1-3) edge node[right] {$\cong$}  (a-2-3);
\end{tikzpicture}
\]
The right square implies that $p$ is an open immersion. 
Since $q^{-1}(\frakU_J^-\otimes {\mathbb W})=\frakU_J^-\otimes {\mathbb W} \times \bbG_{\on{m}{\mathbb{W}}}^{\on{cent}}$ is closed in $\widehat{L_{\mathbb W}\uG}$, fppf (or smooth) descent for closed immersions implies that $p$ is a closed immersion as well. 
Since $\frakU_J^-\otimes {\mathbb W}$ is non-empty and $L^{--}_{\mathbb W}\ucG_{\bbf_J}$ is connected, the map $p$ is an isomorphism.
\xpf

\rema
The above picture extends to the integers $\bbZ$ (in other words, to the wildly ramified cases) by work of the second-named author, resting on complicated group-theoretic constructions described in \cite[\S3, App.]{Lou19}.
\xrema

\rema
For any algebraically closed field $k$ which is a ${\mathbb W}$-algebra, we can identify $\frakU_J^-(k)\cong L^{--}_{\mathbb W}\ucG_{\bbf_J}(k)$ inside $L_{\mathbb W}\uG(k)$, by means of a combinatorial argument. 
For this, recall that $\frakU_J^-(k)$ is generated by the $J$-negative real root subgroups it contains. The fact that $L^{--}_{\mathbb W}\ucG_{\bbf_J}(k)$ shares the same generation property is more or less implicit, though seemingly never explicitly proved, in Kac-Moody theory (compare with \cite[App.~2]{Tit85} and \cite[\S1.3, \S4]{Tit89}); we will give a proof, for general $J$, for completeness. 
First of all, assume that $J=\emptyset$, so that $\bbf_J=\bba$ is an alcove.  Let $S$ be a maximal $k\rpot{t}$-split torus of $G$ whose corresponding apartment contains ${\bf a}$, and
note that $L^{-}_{\mathbb W}\ucG_{\bba}(k)=S(k) \ltimes L^{--}_{\mathbb W}\ucG_{\bba}(k)$ fits into an adequate Birkhoff decomposition by \cite[Prop.~1.1]{HNY13} which we can compare with the induced Birkhoff decomposition on $L_{\mathbb W}\uG$ coming from \eqref{Birkhoff_J_new} for the Kac-Moody group in the case $J = \emptyset$ (this is legitimate because the two groups differ by at most $\bbG_{\on{m}}$-factors in the maximal torus).
Write $g\in L^{--}_{\mathbb W}\ucG_{\bba}(k)$ as $u \dot{w}b$ with $u \in U^-_\emptyset$ and $b \in L^{+}_{\mathbb W}\ucG_{\bba}(k)$;
we see that $w=1$ and $b$ lies in $L^{--}_{\mathbb W}\ucG_{\bba}(k) \cap L^{+}_{\mathbb W}\ucG_{\bba}(k)=1$ (recall Corollary \ref{Open_Cor}). 
This proves $L^{--}_{\mathbb W}\ucG_{\bba}(k) = U^-_\emptyset$, hence the generation result for the left hand side.  
In general, the subgroup $U^-_\emptyset$ is the semi-direct product of $U_J^-$ and the group generated by all negative real root groups contained in $L_J \subseteq L^+_{\mathbb W}\ucG_{\bbf_J}(k)$, so the inclusion $L^{--}_{\mathbb W}\ucG_{\bbf_J}(k) \subseteq U^-_\emptyset$ and the triviality of the intersection $L^+_{\mathbb W}\ucG_{\bbf_J} \cap L^{--}_{\mathbb W}\ucG_{\bbf_J}$ forces the desired equality $L^{--}_{\mathbb W}\ucG_{\bbf_J}(k)= U_J^-$.
\xrema

\rema 
If we transport the $\bbG_{\on{m}}$-action on $\frakG$ given by a $J$-regular dominant weight to $LG$, we can almost immediately deduce that $L^+\ucG_{\bbf_J}$ is the attractor locus and $L^{--}_{\mathbb W}\ucG_{\bbf_J}$ is the strict repeller locus.
Suffice it to say that our proof of Lemma \ref{negative.unipotent.lem} was heavily inspired by the dynamical method of \cite[\S2.1]{CGP15} and should be regarded as an infinite-dimensional generalization thereof.
\xrema

\section{Tangent spaces at base points}\label{tangent.spaces}

Here we combine the results from Sections \ref{negative.loop.grps}--\ref{Kac-Moody-Interlude} to give an effective criterion for the normality of Schubert varieties in Section \ref{normality.crit.sec}. 
With a view toward a future classification of all non-normal Schubert varieties, we state our results over the ring of $p$-typical Witt vectors, which provides a possibly useful link between the characteristic $p$ and characteristic $0$ settings.

\subsection{Preliminaries on tangent spaces}
We start with some general properties of tangent spaces of (ind)-schemes over a general base equipped with a section.

\defi
Let $S$ be a scheme and $X$ be a sheaf of sets on the category of $S$-schemes equipped with $x \in X(S)$. The tangent space $T_{x}X$ of $X$ at $x$ is the sheaf which associates an $S$-scheme $T$ to the pre-image $p_{\varepsilon}^{-1}(x_T)$ of $x_T \in X(T)$ induced by $x$ along the map $p_{\varepsilon}:X(T[\varepsilon])\rightarrow X(T)$. Here by definition $T[\varepsilon] = T \times {\Spec} \, \mathbb Z[\varepsilon]$ where $\varepsilon^2 = 0$.
\xdefi

If $X$ is representable by a scheme, our tangent space coincides with the (implicit) definition of Demazure-Gabriel (see \cite[II, \S4, Cor.~3.3]{DG70}):

\prop\label{tangent_space_prop}
Let $X\to S$ be a scheme endowed with a section $x\co S \to X$. For any $S$-scheme $T$, there is a natural bijection of sets $T_{x}X(T)=\Hom_{\calO_T}(x_T^*\Omega_{X_T/T},\calO_T)$. In particular, $T_{x}X(T)$ has a natural structure of an $\Ga(T,\calO_T)$-module. 
\xprop
\pf Since $S$ is an arbitrary scheme, we reduce to the case $S=T$. To give $f\in T_{x}X(S)$, i.e., a morphism $f\co S[\varepsilon]\to X$ compatible with $x$, is the same as to give an $S$-derivation $d\co \calO_X\to x_*\calO_S$: since $|S[\varepsilon]|=|S|$ on topological spaces, such an $f$ is the same as a morphism of sheaves of rings $f^{\#}\co \calO_X\to x_*\calO_{S[\varepsilon]}=x_*\calO_S\oplus\varepsilon x_*\calO_S$. The compatibility with $x$ implies $f^{\#}=x^{\#}+\varepsilon d_f$ and it is easily verified that $d_f\in \Der_S(\calO_X,x_*\calO_S)$ is an $S$-derivation. We thus get natural bijections 
\[
T_xX(S)=\Der_S(\calO_X,x_*\calO_S)=\Hom_{\calO_X}(\Omega_{X/S},x_*\calO_S)=\Hom_{\calO_S}(x^*\Omega_{X/S},\calO_S),
\]
where the second identification is \cite[01UR]{StaProj}, and the last identification is the adjunction. 
\xpf

Hence, if $T=\Spec(R)$ is an affine scheme, then $T_{x}X(R)$ is an $R$-module by the preceding proposition.

\coro\label{injective_tangent_cor}
Let $R$ be a ring and $i:(X,x) \rightarrow (Y, y)$ be a monomorphism of pointed $R$-schemes. Then the induced homomorphism $i_*\co T_{x}X(R) \rightarrow T_{y}Y(R)$ of $R$-modules is injective. If $i$ is an open immersion, then this homomorphism is bijective.
\xcoro
\pf
This is immediate from the definition, and the fact that monomorphisms are formally unramified, and open immersions are formally \'{e}tale (see \cite[Prop.~17.1.3.(i)]{EGAIV4}), combined with the exact sequence $i^*\Omega_{Y/R} \to \Omega_{X/R} \to \Omega_{X/Y} \to 0$.
\xpf

\coro \label{ind.scheme.tangent}
If $X=\on{colim} X_i$ is a strict pointed ind-scheme over $R$, then $T_xX=\on{colim} T_xX_i$ \textup{(}here $T_xX_i:=0$ if $x\not\in X_i(R)$\textup{)} is an $R$-module independent of the chosen presentation as a strict ind-scheme.
\xcoro
\pf This is immediate from Corollary \ref{injective_tangent_cor}.
\xpf

It is worth noting that the tangent space does not commute with base change in general, whereby we mean the equality $T_{x}X(R)\otimes_R R' \rightarrow T_{x}X(R')$ for all $R$-algebras $R'$, but we still have the following:

\lemm\label{base_change_tangent_lem}
Maintain the notation of Proposition \ref{tangent_space_prop}, and suppose moreover that $X\to S$ is of finite type and that $T=\Spec \, R$ is a Dedekind scheme. Then, for all $R$-algebras $R'$, the canonical homomorphism $T_{x}X(R)\otimes_R R' \rightarrow T_{x}X(R')$ is injective. Moreover, it is bijective for all $R$-algebras $R'$ if and only if $x^*\Omega_{X/R}$ is torsion-free.
\xlemm
\pf 
We may assume $S = T$ and hence that $X \to S$ is of finite presentation. After localizing, we can write $x^*\Omega_{X/R}=R^n \oplus M$ where $M$ is a (finitely generated) torsion module, cf.~\cite[01V3]{StaProj}. Let $X':=X_{R'}$ with induced section denoted $x'$. Since $x'^*\Omega_{X'/R'}=R'^n \oplus (M\otimes_RR')$ by \cite[01UV]{StaProj}, we get 
\[
T_{x}X(R')=\text{Hom}_{R'}(x'^*\Omega^1_{X'/R'}, R')=R'^n \oplus \text{Hom}_{R'}(M\otimes_RR', R'). 
\]
Using $\Hom_R(M,R)=0$ (because $R$ is torsion-free), the lemma follows. This also shows that bijectivity is equivalent to $\Hom_R(M,R')=0$ for all $R$-algebras $R'$, which in turn amounts to asking $M=0$ by Nakayama's lemma.
\xpf

\lemm \label{flat.coker.tangent}
Suppose $i\co (X,x) \rightarrow (Y,y)$ is a closed immersion of pointed ind-schemes of ind-finite type over a Dedekind ring $R$. Then the cokernel of $i_*\co T_{x}X(R) \rightarrow T_{y}Y(R)$ is a flat $R$-module.
\xlemm
\pf
By Corollary \ref{ind.scheme.tangent}, we may and do assume that $X$ and $Y$ are finite-type schemes. After localization, we may also assume that $R$ is a discrete valuation ring with uniformizer $\pi$. Assume there is a $v \in T_{y}Y(R) \setminus T_{x}X(R)$ such that $\pi v \in T_{x}X(R)$. By Corollary \ref{injective_tangent_cor} and Lemma \ref{base_change_tangent_lem}, we have have injections
\[
T_{x}X(R)/\pi T_{x}X(R) \;\subset\; T_{x}X(R/\pi R) \;\subset\; T_{y}Y(R/\pi R).
\]
Since $\pi v=0$ in $T_{y}Y(R/\pi R)$, we get $\pi v\in \pi T_{x}X(R)$, i.e., the existence of some $w \in T_{x}X(R)$ such that $\pi w=\pi v$. As $T_{y}Y(R) = \text{Hom}_R(y^*\Omega^1_{Y/R}, R)$ is free of finite rank (so in particular $R$-torsion free) by the proof of Lemma \ref{base_change_tangent_lem}, we reach a contradiction. This proves the lemma.
\xpf

\subsection{Tangent spaces of affine flag varieties}\label{tangent.space.flag.sec}
Let us give a description of the tangent space of the affine flag variety.
We proceed with the assumptions and notations of Section \ref{negative.loop.grps}, except that we will henceforth often abbreviate $L^{--}_{\mathbb W}\uG$, etc, by omitting the subscript ${\mathbb W}$ and writing $L^{--}\uG$, etc. In what follows, ${\frakh}$ (respectively $L^{--}\frakh$, respectively $\fraku^{\rm op}_H$) denotes the Lie algebra of $H$ (respectively $L^{--}H$, respectively $U^{\rm op}_H$) taken over ${\mathbb W}(k)$.
 
\lemm\label{tangent.space.cor} For the tangent space at the base point $e\in \uFl_{G,\bbf}(R)$ with values in a ${\mathbb W}(k)$-algebra $R$, one has as $R$-modules
\[
T_e\uFl_{G,\bbf}(R)\;=\;T_eL^{--}\ucG_\bbf(R)\;=\; \Bigg(\bigcap_{w\in W_{H,\bbf}} {^w\big(}(L^{--}\frakh) \oplus \fraku_H^\opp\big)\otimes R\Bigg)^\sig
\]
This is a free $R$-module and its formation commutes with arbitrary base change.
\xlemm
\pf
By Corollary \ref{Open_Cor}, the map $L^{--}\ucG_\bbf\to \uFl_{G,\bbf}, g\mapsto g\cdot e$ is representable by a quasi-compact open immersion. This immediately implies 
 \[
 T_e\uFl_{G,\bbf}(R)=T_eL^{--}\ucG_\bbf(R)=T_e(L^{--}\ucH_\bbf)^\sig(R).
 \]
 Using that $T_e(\str)$ commutes with taking fixed points $(\str)^\sig$ and intersections, the corollary follows from Definition \eqref{strictly.negative.split}.
 
 Next we observe that the $R$-module $T_e\uFl_{G,\bbf}(R)$ is projective, as the $\sigma$-averaging map furnishes a retraction to its inclusion in the free module
 $$\bigcap_{w\in W_{H,\bbf}} {^w}\big((L^{--}\frakh)\oplus \fraku_H^\opp\big)\otimes R.$$
(Note that the order of $\sigma$ (which is $e$) is a unit in ${\mathbb W}(k)$, hence in $R$ as well.) A similar argument shows that the tangent space is compatible with base change, i.e. the natural map $T_e\uFl_{G,\bbf}({\mathbb W}(k))\otimes R \rightarrow T_e\uFl_{G,\bbf}(R)$ is an isomorphism for all ${\mathbb W}(k)$-algebras $R$ (use the $\sigma$-averaging retraction applied to the obvious equality in the split case). Hence it suffices to observe that $T_e\uFl_{G,\bbf}({\mathbb W}(k))$ is free, which follows from Kaplansky's theorem, because ${\mathbb W}(k)$ is local.
\xpf

\exam \label{absolutely.special.tangent} Assume that $\bbf_J =: \bbf=\bf 0$ is the base point of $\scrA(G,S,F)$ which is an absolutely special vertex. Recall that absolutely special vertices exist for all quasi-split groups by \cite[Lem.~5.2]{HR20a}. In this case, $L^{--}\calG_\bbf=(L^{--}H)^{\sig,o}$ so that we obtain
\[
T_e\uFl_{G,\bbf}(R)=(L^{--}\frakh\otimes R)^{\sig}=\bigoplus_{i\geq 1} \big(\frakh\otimes R[u^{-i}] \big )^\sig.
\]
(Here $R[u^{-i}]$ is just notation meaning the $R$-span of the monomial $u^{-i}$.)
\xexam

\subsection{Tangent spaces of Schubert varieties}\label{schubert.tangent.sec}
Within this subsection, we additionally assume $G$ to be simply connected. 
By Proposition \ref{kac.moody.comp}, there is a canonical isomorphism of ${\mathbb W}$-ind-schemes $\underline\Fl_{G, \bbf}\cong \frakF_J\otimes {\mathbb W}$ inducing isomorphisms of integral Schubert varieties $\underline S_w\cong \frakS_w\otimes {\mathbb W}$ for all $w\in W/W_\bbf$. 
Given any ample line bundle $\calL$ on $\uFl_{G,\bbf}$, it admits a unique equivariant action of $\widehat{L\uG}$, which in particular naturally acts on $\Gamma(\uFl_{G,\bbf},\calL)^\vee:=\on{colim}_w\Gamma(\uS_w,\calL)^\vee$.  
Restricting this action to $L^{--}\ucG_{\bbf}$ and taking the tangent spaces at base points, we obtain the action of $T_e\uFl_{G,\bbf}$ on $\Gamma(\uFl_{G,\bbf},\calL)^\vee$.
Under the isomorphisms in Proposition \ref{kac.moody.comp}, this is nothing but the Kac-Moody action used in Corollary \ref{tangent.space.formula.and.independence.of.characteristic}.

\lemm\label{tangent.schubert.lemm}
The $R$-valued tangent space $T_e\underline S_w(R)$ identifies with the submodule of $T_e\uFl_{G,\bbf}(R)$ consisting of those $X$ such that $X \Theta^\vee_\calL$ lies in $\Gamma(\uS_w,\calL)^\vee$, where $\Theta_\calL \in \Gamma(\uFl_{G,\bbf},\calL)$ is the usual theta divisor attached to $\calL$ with support given by the complement of $L^{--}\ucG_{\bbf}\cdot e$ and $\Theta^\vee_\calL$ denotes the unique element in the dual weight space sending the theta divisor to $1$ and all other weight spaces to $0$.
\xlemm
\pf
 Now that we have defined a general notion of tangent spaces for any ring $R$, we can repeat the proof of Corollary \ref{tangent.space.formula.and.independence.of.characteristic} for arbitrary $R$, using the isomorphism $\frakU_J^-\otimes {\mathbb W}\cdot e\cong L^{--}\ucG_{\bbf}\cdot e$ from Proposition \ref{kac.moody.comp}.
\xpf

\subsection{Application to the normality criterion}\label{normality.crit.sec}
Let us now turn to our effective criterion for normality, namely Corollary \ref{normality.criterion.cor} below. 
Let $G$ be a tamely ramified, absolutely almost simple, semisimple $F$-group which has the same splitting field as its simply connected cover $G_\scon\to G$. 
The set-up of Section \ref{negative.loop.grps} applies to both groups $G_\scon$, $G$ and we use it to determine the kernel of the map $\Fl_{G_\scon,\bbf}\to \Fl_{G,\bbf}$ on tangent spaces at the base points, cf.~Corollary \ref{normality.criterion.cor}. 

We proceed with the notation of Section \ref{negative.loop.grps}. The map $H_\scon\to H$ on Chevalley groups extends to a map on parahoric $k\pot{u}$-group schemes $\calH_{\scon,\bbf}\to \calH_\bbf$. This induces a map on strictly negative loop groups $L^{--}\calH_{\scon,\bbf}\to L^{--}\calH_\bbf$ over $k$, and hence a map on twisted strictly negative loop groups
\begin{equation}\label{map.strictly.negative.twisted}
L^{--}\calG_{\scon,\bbf}\to L^{--}\calG_\bbf.
\end{equation}
(In this subsection, we abbreviate the functor $L^{--}_{k}$, which has the meaning  analogous to $L^{--}_{\mathbb W}$,  by $L^{--}$.)
We want to determine the kernel of \eqref{map.strictly.negative.twisted}. There is a central extension of flat affine $\bbZ$-group schemes
\[
1\to Z_H\to H_\scon\to H\to 1,
\]
where $Z_H$ is a suitable $\sigma_0$-invariant subgroup of the center of $H_\scon$. Then $Z_H$ is a finite flat $\bbZ$-group scheme of multiplicative type which is \'etale over $\bbZ_{(p)}$ if and only if $p\nmid \# Z_H = \#\pi_1(G)$. 

\defi Let $Z$ be the kernel of $G_\scon\to G$. The strictly negative loop group for $Z$ over $k$ is the subgroup functor of $LG_\scon$ defined as 
\[
L^{--}Z\defined (L^{--}Z_{H})^{\sig,o}\subset (LH_{\scon})^\sig=LG_\scon.
\]
\xdefi
Note that $L^{--}Z$ is representable by a closed subgroup ind-scheme of $LG_\scon$.

\lemm \label{short.exact.lem}
There is a short exact sequence of group functors
\[
1\to L^{--}Z \to L^{--}\calG_{\scon,\bbf}\to L^{--}\calG_\bbf.
\]
\xlemm
\pf 
Clearly, there is a short exact sequence $1\to L^{--}Z_{H}\to L^{--}H_{\scon}\to L^{--}H$. 
Using that $U_{H_\scon}^\opp=U_H^\opp$ for the opposite unipotent radicals and that $W_{H_\scon,\bbf}=W_{H,\bbf}$ in \eqref{strictly.negative.split}, we obtain a short exact sequence
\[
1\to L^{--}Z_{H}\to L^{--}\calH_{\scon,\bbf}\to L^{--}\calH_{\bbf}.
\]
The lemma now follows from Definition (\ref{strictly.negative.split}) by passing to $\sig$-invariants (which is left exact) and by taking neutral components. 
\xpf

By Corollary \ref{tangent.space.cor}, we obtain a $k$-vector subspace
\[
T_{e}L^{--}Z\;\subset\; T_{e}L^{--}\calG_{\scon,\bbf}\;=\; T_{e}\Fl_{G_\scon,\bbf},
\]
where $e\in\Fl_{G_\scon,\bbf}(k)$ denotes the base point. 
Recall from \eqref{schubert.map} that there is a map of Schubert varieties $S_{\scon,w}=S_{\scon,w}(\bba,\bbf)\to S_w(\bba,\bbf)=S_w$ for each $w\in W_\aff/W_\bbf$.

\coro \label{normality.criterion.cor}
For each class $w\in W_\aff/W_\bbf$, the following are equivalent:
\begin{enumerate} 
\item The Schubert variety $S_w\subset \Fl_{G,\bbf}$ is normal.
\item One has
\[
(T_eL^{--}Z)\cap (T_eS_{\scon, w})\;=\; 0
\]
as $k$-vector subspaces of $T_e\Fl_{G_\scon,\bbf}$.
\end{enumerate}
\xcoro
\pf By Proposition \ref{normality.prop}, part (1) is equivalent to $\ker(T_{e}S_{\scon,w}\to T_eS_w)=0$ where $e$ denotes the base point of both $\Fl_{G_\scon,\bbf}$ and $\Fl_{G,\bbf}$. Lemma \ref{short.exact.lem} implies that there is an exact sequence of $k$-vector spaces
\[
0\to T_eL^{--}Z\to T_eL^{--}\calG_{\scon,\bbf}\to T_eL^{--}\calG_\bbf,
\]
so that $\ker(T_{e}S_{\scon,w}\to T_eS_w)=(T_eL^{--}Z)\cap (T_eS_{\scon, w})$. This proves the corollary.
\xpf

\rema \label{5.13_rem}
By \cite[Thm.~0.2]{PR08}, the ind-scheme $\Fl_{G_\scon,\bbf}$ is reduced so that $\Fl_{G_\scon,\bbf}=\on{colim}_w S_{\scon, w}$ is a presentation. 
Thus, Corollary \ref{normality.criterion.cor} shows that, if $L^{--}Z$ is non-trivial, there are infinitely many $(\bba,\bbf)$-Schubert varieties inside $\Fl_{G,\bbf}$ which are not normal, hence not weakly normal, not Frobenius split and not Cohen-Macaulay.
\xrema

\section{Towards a classification of normal Schubert varieties}\label{classification}
Let $k$ be an algebraically closed field of characteristic $p>0$, and let $G$ be a tamely ramified, absolutely simple group over $F=k\rpot{t}$. 
Examining the tables in \cite[Ch.~VI, Planche IX]{Bou68} and \cite[\S4]{Tit79}, here is the list of all such pairs $(G,p)$ such that $p\mid \#\pi_1(G)$. 
Split groups:
\begin{itemize}
\item type $A_{n}$, $n\geq 1$ and $p\,|\,n+1$; 
\item type $B_n$, $n\geq 2$ and $p=2$;
\item type $C_n$, $n\geq 2$ and $p=2$;
\item type $D_n$, $n\geq 3$ and $p=2$;
\item type $E_6$ and $p=3$;
\item type $E_7$ and $p=2$.
\end{itemize}
The split groups $E_8$, $F_4$ and $G_2$ have connection index $1$, and hence are excluded from the list. 
Twisted groups: 
\begin{itemize}
\item type $B\str C_{n}$, $n\geq 3$ and $p\,|\,2n$, $p\not =2$ (even unitary); 
\item type $C\str BC_n$, $n\geq 1$ and $p\,|\,2n+1$ (odd unitary);
\item type $F_4^I$ and $p=3$ (ramified $E_6$);
\item type $G_2^I$ and $p=2$ (ramified triality);
\end{itemize}
The twisted orthogonal groups $C\str B_n$, $n\geq 2$ are excluded by our tamely ramified hypothesis. 

The methods developed in the preceding paragraphs allow us to give a quantitative criterion for the normality of Schubert varieties in general partial affine flag varieties, see Propositions \ref{explicit.bounds.prop} and \ref{general.facets.prop}.
The key input is the computation of the tangent spaces of quasi-minuscule Schubert varieties in twisted affine Grassmannians for absolutely special vertices in Section \ref{absolutely.special}.
In Section \ref{example.pgl2} we discuss the example of $\PGL_2$ in characteristic $2$ which is much easier.
In general the classification of all finitely many normal Schubert varieties in the flag variety for each pair $(G,p)$ in the above list seems to be a challenging problem, see Section \ref{classification.rmk} for some further discussion.

\subsection{Absolutely special vertices} \label{absolutely.special}
We proceed with the assumptions and notations of Section \ref{normality.crit.sec}. 
In particular, $G$ is a tamely ramified, absolutely almost simple, semisimple $F$-group which has the same splitting field as its simply connected cover $G_\scon\to G$.

We further assume $\bbf=\bf 0$ is the fixed absolutely special vertex in $\scrA(G,S,F)$. 
Our aim is to give an effective criterion for the normality of $(\bba,\bf 0)$-Schubert varieties inside the neutral component of the twisted affine Grassmannian $\Gr:=\Fl_{G,\bf 0}$. 
For this, we study the tangent spaces of $(\bba, \bf 0)$-Schubert varieties inside $\Gr_\scon:=\Fl_{G_\scon,\bf 0}$.
The $L^+\calG_{\scon,\bba}$-orbits inside $\Gr_\scon$ are enumerated by the set $W_\aff/W_{\bf 0}=X_*(T_\scon)_I$, the coinvariants under the Galois group $I:=\Gal(F'/F)$ where $F'/F$ is the splitting field. 
For each $\bar{\mu}\in X_*(T_\scon)_I$, we denote by $S_{\scon, \bar\mu}\subset \Gr_\scon$ the corresponding $(\bba,\bf 0)$-Schubert variety. 
In view of Corollary \ref{normality.criterion.cor} we have to determine exactly those $\bar\mu\in X_*(T_\scon)_I$ such that $(T_eL^{--}Z)\cap (T_eS_{\scon, \bar\mu})=0$ inside 
\begin{equation}\label{tangent.aff.grass}
T_e\Gr_\scon\;=\; \bigoplus_{i\geq 1} \big( \frakh_\scon[u^{-i}]\big )^\sig, 
\end{equation}
cf.~Example \ref{absolutely.special.tangent} (in particular note that $\frakh_\scon[u^{-i}]$ is just our notation for $u^{-i}\frakh_\scon$). 
Our normality criterion rests on the following key calculation.

\prop\label{quasi.minuscule}
 Let $\bar\mu\in X_*(T_\scon)_I$ be the unique $B$-dominant, quasi-minuscule element. Then
\[
T_eS_{\scon, \bar\mu}\;\supset\; \big(\frakh_\scon[u^{-1}]\big)^\sig
\]
as $k$-vector subspaces of \eqref{tangent.aff.grass}, and equality holds if $\on{char}(k)=0$. 
\xprop

\pf For the proof of this inclusion, we may and do assume $\on{char}(k)=0$ by Lemmas \ref{base_change_tangent_lem} and \ref{flat.coker.tangent} combined with Proposition \ref{schubert.var.int.prop}, all applied to the normal Schubert variety $\uS_{\scon, \bar{\mu}}$.
The equality follows then from our work with minimal nilpotent orbits in Appendix \ref{app-minimal-nilpotent-orbits}, which extends previous results of \cite[\S2.9]{MOV05} and \cite[\S8]{HR20a}. 
(Note that $\bar{\mu} = \theta^\vee$, where $\theta$ is the highest root in the \'{e}chelonnage root system for $G$; see \cite{HR08} and Section \ref{classification.rmk}.)
For convenience of the reader, let us just note that the inclusion
$T_eS_{\scon, \bar\mu}\;\supset \; \big(\frakh_\scon[u^{-1}]\big)^\sig$
is much simpler - and this is all we will need to prove the important Corollary \ref{qmin_in_Gr} below. 
Indeed, the intersection of both sides is certainly non-trivial, as we see by looking at the $(\bba,0)$-Schubert variety $S_{w}$ for the affine simple reflection $w=s_0$. 
Moreover, both tangent spaces carry an action by the split group $H^{\sig_0}_\scon$. 
Now we use that the right side is an irreducible $H^{\sig_0}_\scon$-module: this is obvious in the split case, because we get the adjoint representation; in the twisted case, it is proved in Proposition \ref{repn.nilp.orbit.reduced}, Proposition \ref{repn.nilp.orbit.non.reduced} and \cite[Lem.~8.4]{HR20a}.
\xpf

\coro  \label{qmin_in_Gr} 
If $p\mid \#\pi_1(G)$, then the quasi-minuscule Schubert variety inside $\Gr_G$ is not normal. 
\xcoro
\pf 
Let $\frakz_{H}$ denote the Lie algebra over $k$ of the kernel $Z_H$ of $H_\scon\to H$; cf.~Section \ref{normality.crit.sec}. Note $\frakz_{H}$ is nonzero since $Z_H$ is not \'etale over $k$ by assumption.
Combining Proposition \ref{quasi.minuscule} with Corollary \ref{normality.criterion.cor} (2), it is enough to show that the subspace $(\frakz_H[u^{-1}])^{\sigma}\subset (\frakh_\scon[u^{-1}])^\sigma$ is non-trivial.
If $G$ is split, so that $\sig$ acts trivially, then $\frakz_H[u^{-1}]$ is clearly non-trivial. 
If $G$ is non-split, we go through the possible types for $H$ listed in the beginning of Section \ref{classification}. 
First for simplicity assume $G$, {and hence $H$}, is adjoint, so that $Z_H$ is the center of $H_\scon$ and $\frakz_H$ is the center of $\frakh_\scon$, that is, the kernel of the adjoint representation, see \cite[Prop.~3.3.8 ff.]{Con14}.
If $H$ is of type $A_n$, then $\frakz_H$ is spanned by the element
$$ 
\sum_{i=1}^n i\alpha_i^\vee.
$$
Then we notice the congruence $n+1-i\equiv-i $ modulo $p$ and use that $\al_i^\vee\mapsto \al_{n+1-i}^\vee$ and $u^{-1}\mapsto -u^{-1}$ under $\sig$. 
If $H$ is of type $D_4$, then $\sigma_0$ permutes the roots as follows: $\alpha_1 \mapsto \alpha_3 \mapsto \alpha_4 \mapsto \alpha_1$, $\alpha_2 \mapsto \alpha_2$.  
It follows that $\frakz_H$ contains in characteristic $p=2$ the element
$$ 
\alpha_1^\vee+\zeta^{-1}\alpha_3^\vee+\zeta^{-2}\alpha_4^\vee,
$$
which becomes $\sigma$-invariant after multiplying by $u^{-1}$. Here $\zeta$ is a primitive $3$rd root of unity in the notation of Section \ref{tame_lift_subsec}. To check the containment, multiply the Cartan matrix by the column vector $(1,0,\zeta^{-1},\zeta^{-2})$ and show that the sum of the entries in each row vanishes modulo $2$.
Finally in the $E_6$ type, a similar argument in characteristic $p=3$ shows that $\frakz_H$ contains
$$\alpha_1^\vee+2\alpha_3^\vee-2\alpha_5^\vee-\alpha_6^\vee. $$
To show that this element becomes $\sigma$-invariant after multiplying by $u^{-1}$ use that $\alpha_1^\vee\leftrightarrow \alpha_6^\vee$, $\alpha_3^\vee\leftrightarrow \alpha_5^\vee$ and $u^{-1}\mapsto -u^{-1}$ under $\sig$.
This proves the corollary in the adjoint case.

To handle the general non-split cases where $H$ is of type $A_n$, note that $\frakz_H\not = 0$ must be the entire $1$-dimensional center of $\frakh_\scon$.
In type $E_6$, the group $Z_H$ is a non-trivial subgroup of the center $Z_{H_\scon}=\Bmu_3$ of $H_\scon$, hence equal, so that $G$ must be adjoint in this case.
If $H$ is of type $D_4$, then the center of $H_\scon$ is $Z_{H_\scon}={\Bmu}_2 \x {\Bmu}_2$, with Lie algebra $\frakz_{H_\scon} = k \oplus k$.  
The order $3$ automorphism $\sigma_0$ preserves $Z_{H_\scon}$ and acts on $\frakz_{H_\scon}$ (up to choice of basis) by $e_1 \mapsto e_1 + e_2$, $e_2 \mapsto e_1$. Now if $1 \subsetneq Z_H \subsetneq Z_{H_\scon}$ were $\sigma_0$-invariant and non-\'etale, then $Z_H \cong {\Bmu}_2$ (use Cartier duality), and hence $\sigma_0$ acts trivially on $Z_H$ (since ${\rm Aut}({\Bmu}_2) \cong {1}$). In this case $\sigma_0$ would fix a vector in $\frakz_{H_\scon} = k \oplus k$, a contradiction. 
It follows that only $Z_H = Z_{H_\scon}$ occurs so that $G$ must be adjoint in this case as well.
This proves the corollary.
\xpf

Using the absolutely\footnote{For this discussion, any special vertex will do.} special vertex $\bf 0 \in \bar{\bf a}$, we identify $\scrA = \scrA(G,S,F)$ with $X_*(T)_{I,\mathbb R}$, where $I = {\rm Gal}(\bar{F}/F)$. Recall that the Iwahori-Weyl group $W$ acts by affine linear transformations on $\scrA$. We use the Bruhat-Tits convention: $t \in T(F)$ acts by translation by $-\kappa_T(t)$\footnote{More precisely, it acts by the image of this element in $X_*(T)_{I, \mathbb R}$; recall $X_*(T)_I$ might have torsion.}, where $\kappa_T\co T(F) \twoheadrightarrow X_*(T)_I$ is the Kottwitz homomorphism constructed in \cite[\S7]{Kot97}. Following \cite[Prop.~13, Lem.~14]{HR08}, we get isomorphisms 
$$W \overset{\sim}{\rightarrow} W_{\rm aff} \rtimes \Omega_{\bf a} \cong X_*(T)_I \rtimes W_{\bf 0},$$ 
where $\Omega_{\bf a}$ is the subgroup of $W$ preserving ${\bf a}$, where the map ${\rm Norm}_GT(F) \rightarrow X_*(T)_I \rtimes W_{\bf 0}$ extends $\kappa_T\co T(F) \twoheadrightarrow X_*(T)_I$.  
In particular, we have the group embedding $X_*(T)_I \hookrightarrow W$ denoted $\nu\mapsto t^\nu$ where $t^\nu$ is characterized by the property $\kappa_T(t^\nu) = \nu$ (if $T$ is split, $t^\nu = \nu(t)\,\,\mbox{mod ${\rm ker}(\kappa_T)$}$). According to the Bruhat-Tits convention, the element $t^\nu$, and hence $\nu$, acts on $\scrA$ by translation by the image of $-\nu$ in $X_*(T)_{I, \mathbb R}$. 
We may view $W_{\rm aff}$ as the Coxeter group generated by the reflections through the walls of ${\bf a}$.
Using the isomorphism, we transport the Bruhat order on $W_{\rm aff} \rtimes \Omega_{\bf a}$ to one on $W$; this induces the Bruhat order on $W/W_{\bf 0}$. Our choice of embedding $X_*(T)_I \hookrightarrow W$ induces a well-defined bijection of sets $X_*(T)_I \overset{\sim}{\rightarrow} W/W_{\bf 0}$, and we consider the transported Bruhat order on $X_*(T)_I$. We are going to need the following combinatorial description of the Bruhat order on $X_*(T)_I$, which can be found in \cite[Thm.~2.5]{BH21} for split groups. 

Recall (cf.~\cite{HR08}) that $W_{\rm aff} = W_{\rm aff}(\Sigma)$ for the \'echelonnage roots $\Sigma = \Sigma(G, S, F)$; these have the property that the hyperplanes annihilated by the affine roots $\Phi_{\rm af}(G,S,F)$ of Tits \cite[\S1.6]{Tit79} are in bijection with those annihilated by the affine functionals of the form $\beta + n$ for $\beta \in \Sigma$, $n \in \mathbb Z$. Let $Q^\vee = \mathbb Z[\Sigma^\vee]$ be the \'echelonnage coroot lattice; it may be identified with $X_*(T_{\rm sc})_I$. In what follows, all finite and affine roots mentioned will be \'echelonnage (affine) roots. Let $C^+$ be the Weyl chamber in $\scrA$ which contains ${\bf a}$ and has apex ${\bf 0}$. We say a finite root $\beta$ (resp.,\,affine root $\beta +n$) is positive (and write $\beta > 0$, resp.,\,$\beta + n > 0$) if it takes positive values on $C^+$ (resp.,\,${\bf a}$). Recall that $W_{\rm aff}$ is the Coxeter group generated by the reflections $s_{\beta +n}$ in the simple positive affine roots $\beta + n$.

\begin{propo}[Besson-Hong]\label{bruhat.order.on.coweights}
Given two coweights $\lambda, \mu \in X_*(T)_I$, the inequality $\lambda \leq \mu$ holds if and only if $\lambda- \mu \in Q^\vee$ and there is a sequence of coweights $\mu_i \in X_*(T)_I$ such that $\mu_0=\mu$, $\mu_r=\lambda$ and satisfying the following: there is a positive root $\alpha_i$ such that either $\mu_{i+1}=\mu_i-k\alpha_i^{\vee}$ with $0 \leq k \leq  \langle \alpha_i, \mu_i\rangle$ or $\mu_{i+1}=\mu_i+k\alpha_i^{\vee}$ with $0 \leq k < - \langle \alpha_i,  \mu_i \rangle$. 

\end{propo}

It was already well-known that, if $\lambda$ and $\mu$ lie in a common Weyl chamber, then the Bruhat order described above coincides with the usual dominance partial order with respect to the given Weyl chamber (cf.~\cite[Lem.~3.8, Prop.~3.5]{Rap05}, \cite[Thm.~4.1]{BH21}).

\begin{proof}

By definition $\lambda \leq \mu$ if and only if $t^\lambda \leq t^\mu$ in the Bruhat order on $W/W_{\bf 0}$. Let $w_\nu \in W$ be the minimal length element in $t^\nu W_{\bf 0}$. The Bruhat order on $W/ W_{\bf 0}$ is generated by the following relation between $w_{\nu'}, w_{\nu}$ for pairs of distinct elements $\nu, \nu' \in X_*(T)_I$: there is an affine reflection $s_{\beta + n}$ with $\beta + n$ positive such that 
$$
w_{\nu'} > s_{\beta + n}w_{\nu'}
$$
in the Bruhat order on $W$, and $s_{\beta + n}w_{\nu'} W_{\bf 0} = w_\nu W_{\bf 0}$; (when this happens we write $w_{\nu'}W_{\bf 0} > s_{\beta + n} w_{\nu'}W_{\bf 0}$.) 
This is the same as saying that $s_{\beta + n}(-\nu') = -\nu$, and the point $-\nu'+ {\bf 0}$ and the alcove ${\bf a}$ are on opposite sides of the affine hyperplane $H_{\beta + n}$, that is, $-\langle \beta, \nu' \rangle + n < 0$.

Therefore, $t^\lambda < t^\mu$ if and only if $\lambda - \mu \in Q^\vee$ and there exists a sequence of reflections $s_i = s_{\beta_i + n_i}$, ($0 \leq i \leq r-1, \, \beta_i + n_i >0$), such that as elements in $X_*(T)_{I,\mathbb R}$ we have $-\mu_0 = -\mu$, $-\mu_r = -\lambda = s_{r-1}\cdots s_0 (-\mu)$, and where, for each $i \geq 0$, if $-\mu_i := s_{i-1}\cdots s_0 (-\mu_0)$, then $-\langle \beta_i, \mu_i \rangle + n_i < 0$.  Of course, we may assume $\mu_0, \dots, \mu_r$ has no repetitions.

By definition $-\mu_{i+1} = s_i(-\mu_i)$, that is,
$$
-\mu_{i+1} = -\mu_i - \big(\langle \beta_i, -\mu_i \rangle + n_i\big) \beta_i^\vee.
$$
Because $\beta_i + n_i$ is a positive affine root, we have $n_i \geq 0$ and $n_i = 0 \Rightarrow \beta_i > 0$.

\begin{enumerate}
\item[(1)]If $\beta_i >0$ then $n_i \geq 0$ and  $\mu_{i+1} = \mu_i - k\beta_i^\vee$ where $k = \langle \beta_i,\mu_i \rangle - n_i$. Note that $0 < k \leq \langle \beta_i, \mu_i \rangle$. Set $\alpha_i = \beta_i$.
\item[(2)] If $\beta_i < 0$ then $n_i \geq 1$, and $\mu_{i+1} = \mu_i + k(-\beta_i)^\vee$, where $k = \langle \beta_i, \mu_i \rangle - n_i$. Note that $0 <k < -\langle -\beta_i, \mu_i \rangle$. Set $\alpha_i = -\beta_i$.
\end{enumerate}
Conversely, given the positive root $\alpha_i$ and integer $k$ satisfying the given restrictions, we may define the positive affine root $\beta_i + n_i$ using (1) or (2), for which we have $-\langle \beta_i, \mu_i \rangle + n_i < 0$.
\end{proof}

In the following we apply this to uniformly bound the subset of normal Schubert varieties for absolutely almost simple semisimple groups such that $p \mid \# \pi_1(G) $, that is, for those semisimple groups $G$ such that $G_{\rm sc} \rightarrow G$ is a non-\'{e}tale isogeny.

\prop\label{explicit.bounds.prop}
Let $G$ be an absolutely almost simple semisimple group whose simply connected cover is a non-\'etale isogeny. Then the set of $\lambda \in Q^{\vee}$ such that $S_\lambda$ is normal is finite. More precisely, it is contained in the finite complement of all $\lambda\in Q^{\vee}$ such that $\lambda \geq -2\theta^{\vee}$, where $\theta$ denotes the highest root for the \'echelonnage root system $\Sigma(G, S, F)$.
\xprop

\pf
We start by observing that $-2\theta^{\vee}$ is bigger than $\theta^{\vee}$. 
Indeed, $-\langle \theta, -2\theta^{\vee} \rangle=4$ and thus $\theta^{\vee}=-2\theta^{\vee}+3\theta^{\vee}$ is less than $-2\theta^{\vee}$ for the partial Bruhat order, see Proposition \ref{bruhat.order.on.coweights}.
By Corollary \ref{normality.cor} combined with Corollary \ref{qmin_in_Gr}, this gives the proposition as soon as we know that the complement of $\{\lambda \in Q^{\vee}\;|\; \lambda \geq -2\theta^{\vee} \}$ in $Q^{\vee}$ is finite. 

Suppose $C_1$ and $C_2$ are two adjacent closed Weyl chambers such that $C_1$ lies in a minimal gallery connecting the dominant Weyl chamber to $C_2$. Then there is a positive root $\alpha$ such that the wall of the reflection $s_{\alpha}$ bounds $C_1$ and $C_2$, in such a way that $C_1$ lies on the nonnegative side with respect to $\alpha$. In particular, if $\lambda \in C_1$, then $s_{\alpha}\lambda \in C_2$ and the inequality $s_{\alpha}\lambda \leq \lambda$ holds, again by Proposition \ref{bruhat.order.on.coweights}, as $\langle \alpha, \lambda \rangle \geq 0$.

Let $Q^\vee_+$ denote the dominant elements in the coroot lattice $Q^\vee$.
The above argument reduces us to considering only antidominant $\lambda$, that is, to showing that the set $\{ \lambda \in Q^{\vee}_+ \;|\; -\lambda \ngeq -\lambda_0 \}$ is finite for any fixed $\lambda_0 \in Q^\vee_+$. 
We will show the equivalent statement that $\{ \lambda \in Q^\vee_+ \; | \; \lambda \ngeq \lambda_0\}$ is finite. 
Dominance ensures we may write $\lambda=\sum n_i \alpha_i^{\vee}$ and $\lambda_0=\sum n_{0,i} \alpha_i^{\vee}$, where $n_i, n_{0,i} \geq 0$ for all $i$.
Writing $\lambda=\sum n_i \alpha_i^{\vee}$ and $\lambda_0=\sum n_{i,0} \alpha_i^{\vee}$, it is enough to prove that for all $j$, $n_j \leq 2^r \,{\rm max}_i \{ n_{0,i}\}$ whenever $\lambda \not\geq \lambda_0$, where $r$ is the number of nodes of the Dynkin diagram for $\Sigma(G,S,F)$.
In this case, by Proposition \ref{bruhat.order.on.coweights} there is some $i$ such that $n_i < n_{0,i}$. 
For $j \neq i$, set $r_{ij} = -\langle \alpha_i, \alpha^\vee_j \rangle \in \mathbb Z$. 
Assuming $\alpha_j$ is adjacent to $\alpha_i$ in the Dynkin diagram, $r_{ij} \geq 1$. 
By the dominance of $\lambda$, we see that $2n_i - r_{ij}n_j \geq \langle \alpha_i, \lambda \rangle \geq 0$, which implies $n_j \leq 2n_i < 2\,{\rm max}_i \{ n_{0,i}\}$. Repeating this argument shows $n_j \leq 2^r \, {\rm max}_i \{n_{0,i} \}$ for all $j$ because the Dynkin diagram for $\Sigma(G,S,F)$ is connected.
\xpf

\subsection{General facets}\label{classification.rmk}
Here we keep virtually all notation introduced in the previous section, in particular we require $G$ is absolutely almost simple, but we no longer assume that $\bbf=\bf0$. We rather assume that $\bbf$ and $\bf0$ are subfacets of the dominant base alcove $\bba$. For any $\lambda \in X_*(T)_I$, let $w^\lambda$ (resp.\,$w_\lambda$) denote the maximal (resp.\,minimal) length element in $t^\lambda W_{0}$.  Let $\theta$ be the highest \'{e}chelonnage root of $G$. Then $\bar{\mu} = \theta^\vee$ is the unique quasi-minuscule coweight for the \'{e}chelonnage root system $\Sigma(G, S, F)$. Fix any regular antidominant element $\delta \in X_*(T)_I$ such that $\delta \geq \theta^\vee$ in the Bruhat order on $X_*(T)_I$.

\prop\label{general.facets.prop}
Let $\tau \in \Omega_{\bf a}$. All but finitely many elements of the form $x\tau \in W_{\rm aff}\tau/W_{\bf f}$ satisfy $x\tau  \geq w_{\delta-\theta^\vee} \tau$ in the Bruhat order on $W/W_{\bf f}$, and for any such element $S_{x\tau}({\bf a}, {\bf f})$ is not normal if $G_{\rm sc} \to G$ is a non-\'{e}tale isogeny.
\xprop

Note that this proposition proves Theorem \ref{fin.many.normal.sch.vars}.

\pf
We can immediately reduce to the case $\tau  = 1$. Since $\delta$ is regular antidominant, we see easily that $w_{\delta-\theta^\vee} > w^{\delta} \geq w^{\theta^\vee}$ in the Bruhat order on $W_{\rm aff}$.  Indeed, since $\delta$ is regular antidominant we have $t^{\delta} = w^{\delta}$ (it is known that $l(t^\delta w) = l(t^\delta) - l(w),\,\, \forall w \in W_{\bf 0}$, by e.g.\,\cite[Prop.~1.23]{IM65}).  The element $t^{-\theta^\vee} s_\theta \in W$ acts on $\scrA(G, S, F)$ by the simple affine reflection $s_0$, and so $t^{\delta} s_0({\bf a}) = t^{\delta-\theta^\vee} s_\theta({\bf a})$ is separated by a wall of type $s_0$ from $t^{\delta}({\bf a})$, with $t^{\delta}({\bf a})$ closer to the base alcove (remember that $t^\delta$ acts by translation by the regular dominant vector $-\delta$). So $l(t^{\delta} s_0) = l(t^{\delta}) + 1$.  It follows that $t^{\delta} s_0 = w_{\delta-\theta^\vee}$, and from this that $w^{\delta} < w^{\delta} s_0 = w_{\delta - \theta^\vee }$.  Finally, observe that $\delta \geq \theta^\vee$ is equivalent to $w_{\delta} \geq w_{\theta^\vee}$, which is equivalent to $w^\delta \geq w^{\theta^\vee}$.

Since $S_{\theta^\vee}({\bf a}, 0)$ is not normal (when $G_{\rm sc} \to G$ is non-\'{e}tale), we deduce that $S_{w^{\theta^\vee}}({\bf a}, {\bf a})$ is not normal, hence also $S_x({\bf a}, {\bf a})$ is not normal whenever $x \geq w_{\delta - \theta^\vee}$ in the Bruhat order on $W_{\rm aff}$.  

Finally we prove that all but finitely many $x \in W_{\rm aff}$ satisfy $x \geq w_{\delta-\theta^\vee}$.  By the proof of Proposition \ref{explicit.bounds.prop}, all but at most finitely many $\lambda \in Q^\vee$ satisfy $w_\lambda \geq w_{\delta-\theta^\vee}$.  For any $w \in W_{0}$ and any such $\lambda$, we have $t^\lambda w \geq w_{\delta-\theta^\vee}$. We are done.
\xpf

\subsection{The example of $\PGL_2$} \label{example.pgl2}

Our methods allow us to give a complete classification of normal Schubert varieties for $\PGL_2$ in characteristic $2$. 
In this subsection, let $k$ be a field of characteristic $2$.

\lemm \label{quasi.min.pgl2.lemm} 
The quasi-minuscule Schubert variety inside the affine Grassmannian for ${\PGL_2}$ is not normal.
More precisely, an open affine neighborhood of the base point is isomorphic to the spectrum of the $k$-algebra
\[
k[x,y,v,w]/(vw+x^2y^2, v^2+x^3y, w^2+xy^3, xw+yv).
\] 
\xlemm

\pf 
Since $2$ divides $\#\pi_1(\PGL_2)=2$, the non-normality is a special case of Corollary \ref{qmin_in_Gr}.
It remains to prove the displayed formula for the coordinate ring. 
By putting $v=xz$, $w=yz$, this $k$-algebra identifies with the subalgebra of $k[x,y,z]/(z^2+xy)$ generated by $x,y, xz, yz$.
Now the lemma follows from the calculations in Appendix \ref{appendix_PGL2}, see Corollary \ref{quasi_min_PGL2_special}.
\xpf

Let $\Fl:=\Fl_{\PGL_2,\bba}$ be the affine flag variety. 
For each $w$ in the Iwahori-Weyl group $W$, we denote by $S_w\subset \Fl$ the associated $(\bba,\bba)$-Schubert variety.

\coro\label{classification.pgl2}
For $w\in W$, the Schubert variety $S_w$ is normal if and only if $\dim(S_w)\leq 2$ in which case it is smooth.
\xcoro
\pf 
After possibly multiplying $w\in W$ with an element in the stabilizer of $\bba$, we may and do assume that $w\in W_\aff$, i.e., $S_w$ lies in the neutral component of $\Fl$.   
The affine Weyl group $W_\aff$ is the free group with generators $s_0, s_1$ and relations $s_0^2=s_1^2=1$. 
Here $s_0$ is the simple affine reflection, and $s_1$ the simple finite reflection. 
Consider the projection $\pi\co \Fl\to \Gr:=\Gr_{\PGL_2}$, a smooth proper morphism of relative dimension $1$. 
Let $S_\mu\subset \Gr$ be the quasi-minuscule Schubert variety, which is not normal by Lemma \ref{quasi.min.pgl2.lemm}. 
Hence, the Schubert variety $\pi^{-1}(S_\mu)=S_w$, $w=s_1s_0s_1$ is not normal. 
By Corollary \ref{normality.cor}, all other Schubert varieties $S_v$ with $v\geq w$ are not normal as well. 
In particular, all Schubert varieties with $\dim(S_w)\geq 4$ are not normal. 
If $\dim(S_w)\leq 1$, i.e., either $w=1$, or $w=s_0$, or $w=s_1$, then $S_w$ is clearly smooth, hence normal. 
In order to treat the remaining cases where $\dim(S_w)=2$ or $\dim(S_w)=3$, we observe that the $(\bba,\bba)$-Schubert variety $S_w$ is normal (resp.~smooth) if and only if the $(\bba,\bba)$-Schubert variety $S_{\tau w\tau^{-1}}$ is normal (resp.~smooth) where $\tau\in W$ is the non-trivial element in the stabilizer of $\bba$, see Lemma \ref{normality.stabilizer.lem}.
We have $\tau w\tau^{-1}=s_1s_0$ for $w=s_0s_1$ and $\tau w\tau^{-1}=s_0s_1s_0$ for $w=s_1s_0s_1$.
Hence, both $3$-dimensional Schubert varieties are not normal as argued above. 
One of the $2$-dimensional Schubert varieties is the preimage in $\Fl$ of the translated to the neutral component minuscule Schubert variety in $\Gr$.
Hence, both $2$-dimensional Schubert varieties are smooth. 
This proves the corollary.
\xpf

\coro\label{classification.pgl2.grassmannian}
A Schubert variety in the affine Grassmannian for $\PGL_2$ in characteristic $2$ is normal if and only if it is at most $1$-dimensional, in which case it is already smooth. 
\xcoro
\pf
This is immediate from Corollary \ref{classification.pgl2} by considering the smooth projection of relative dimension $1$ from the affine flag variety. 
\xpf

\subsection{Some remarks on the classification} 

Our methods from Section \ref{absolutely.special} do not apply to the case of special, but not absolutely special vertices. 
This is only an issue in the case of odd unitary groups of type $C\str BC_n$, $n\geq 1$ with $p\, | \,2n+1$. 
In this case, there are up to $G_\ad(F)$-conjugation two types of special vertices, where exactly one of them is absolutely special, see \cite[\S5]{HR20a}.
Here separate methods seem to be required to calculate the tangent space of the quasi-minuscule Schubert variety in the corresponding twisted affine Grassmannian.
Furthermore, we note that the normality criterion obtained in Proposition \ref{explicit.bounds.prop} is not effective. 
Indeed, this can be seen already in the case of $\PGL_2$ in characteristic $2$ by comparing with the classification in Corollary \ref{classification.pgl2}.
In principle, Corollary \ref{normality.criterion.cor} (2) together with the tangent space formula of Kumar and Polo (Corollary \ref{tangent.space.formula.and.independence.of.characteristic}) gives an effective way of classifying all normal Schubert varieties. 
Here the main difficulty is the determination of the affine Demazure modules.
The case of, say, $\PGL_3$ in characteristic $3$ already seems quite involved.

\section{Reducedness}

In \cite[Thm.~6.1]{PR08}, the authors show that loop groups (equivalently, their partial affine flag varieties) attached to semisimple groups $G$ over a field $k$ are reduced under the hypothesis $\on{char}(k)\nmid \# \pi_1(G) $. 
We show in Proposition \ref{reduced.prop} (split case) and Proposition \ref{reducedness_twisted} (twisted case) that this hypothesis is necessary.

\subsection{The split case}
Throughout this subsection, let $k$ be an arbitrary field and $G$ be a connected split reductive group over $k$. 
We are going to use the notion of distributions, which should be regarded as higher order differential operators. 
For the theory of distributions for (group) schemes we refer to \cite[II, \S4]{DG70} and \cite[\S7]{Jan03}.

\begin{dfn}
Let $(X,x)$, $x\in X(k)$ be a pointed $k$-ind-scheme. 
The {\it space of distributions} $\text{Dist}(X,x)$ is the $k$-vector space obtained as the filtered colimit of the $k$-vector space duals of all Artinian closed subschemes of $X$ supported at $x$. 
\end{dfn}
 
We record some basic properties.

\lemm 
\label{basic.dist.lemm}
Let $(X,x)$, $(Y,y)$ be pointed $k$-ind-schemes, and let $f\co (Y,y)\to (X,x)$ be a map of pointed $k$-ind-schemes. 
\begin{enumerate}
\item If $(X,x)=\on{colim} (X_i,x)$ is any presentation, then $\Dist(X,x)=\on{colim} \Dist(X_i,x)$ with injective transition maps. Further, each $\Dist(X_i,x)$ only depends on the formal spectrum $\Spf(\calO_{X_i,x})$ viewed as an ind-scheme.
\item The map $f$ induces a map $(df)_y\co \Dist(Y,y)\to \Dist(X,x)$. 
\item There is a natural map $\Dist(X,x)\otimes_k\Dist(Y,y)\to \Dist(X\x_kY,(x,y))$ which is an isomorphism if both $X$, $Y$ are ind-\textup{(}locally Noetherian\textup{)} over $k$.\end{enumerate}
\xlemm
\pf 
Part (1) is immediate because the transition maps $X_i\to X_j$ are closed immersions.
Part (2) and (3) follow from (1) and the case of schemes in \cite[I, \S7.2 \& \S7.4]{Jan03}.
Note that {\it loc.~cit.} is over more general base rings, and that the assumptions are satisfied for locally Noetherian schemes over fields.
\xpf

In particular, for any pointed $k$-ind-scheme $(X,x)$ which is ind-(locally Noetherian), e.g., of ind-(finite type), the space of distributions $\Dist(X,x)$ is a cocommutative counital $k$-coalgebra whose coalgebra structure is induced from the diagonal $X\to X\x_k X$ and Lemma \ref{basic.dist.lemm} (3), cf.~\cite[I, \S7.4 (3)]{Jan03} for details.
If $X$ is a $k$-group ind-scheme --possibly of ind-(infinite type)-- then we define
\begin{equation}
\Dist(X)\defined \Dist(X,1),
\end{equation}
where $1\in X(k)$ denotes the neutral section. 
In this case, the action morphism $X\x_kX\to X$ (combined with Lemma \ref{basic.dist.lemm}) induces on $\Dist(X)$ the structure of an associative $k$-algebra under the convolution of distributions, cf.~\cite[I, \S7.7]{Jan03} for details.

For the next lemma recall that a quasi-compact morphism of schemes is called {\it scheme-theoretically dominant} if its scheme theoretic image \cite[01R5]{StaProj} is equal to its target. 

\lemm
\label{surjective.dist.lemm}
Let $f\co (Y,y)\rightarrow (X,x)$ be a quasi-compact, scheme-theoretically dominant morphism of locally Noetherian pointed $k$-schemes. 
Then the induced homomorphism $(df)_y\co\emph{Dist}(Y,y) \rightarrow \emph{Dist}(X,x)$ is surjective.
\xlemm
\pf
Since $f$ is quasi-compact and scheme-theoretically dominant, the induced map $\calO_{X,x}\rightarrow \calO_{Y,y}$ on local rings is injective, cf.~\cite[01R8 (1), (2)]{StaProj}. 
Also note that the map $(df)_y$ only depends on the induced map on completed local rings $\hat\calO_{X,x}\to \hat\calO_{Y,y}$, which is injective as well.
By Krull's intersection theorem, the decreasing sequence of ideals $\{\hat\frakm_y^n\cap \hat\calO_{X,x}\}_{n\geq 1}$ has zero intersection, and hence by Chevalley's lemma \cite[Lem.~7]{Che43} defines a cofinal family of Artinian closed subschemes of $\Spec(\hat\calO_{X,x})$ supported at $x$.
This implies the lemma.
\xpf

\rema
\label{Jantzen.flat.rem}
Another interesting example (cf.~also \cite[I, \S7.6]{Jan03}) to which Lemma \ref{surjective.dist.lemm} applies is the case of a map $f\co (Y,y)\rightarrow (X,x)$ of locally Noetherian pointed $k$-schemes which is flat at $y$. 
Indeed, then the induced map $\calO_{X,x}\to \calO_{Y,y}$ is faithfully flat, and hence injective, that is, the map on spectra is scheme-theoretically dense.
Also we find it instructive to check Lemma \ref{surjective.dist.lemm} ``by hand'' in the special cases of the normalization of the cusp, and the (relative) Frobenius morphism in strictly positive characteristic, say, of the affine line.
\xrema

The previous lemma will be used to show that $\Gr_G$ for adjoint non-(simply connected) groups is non-reduced in bad characteristics by noticing that the $k$-vector space of the distributions of its reduction is strictly smaller. 
The following lemma shows that this space can be easily computed at ``infinite level''.
For later use we formulate this lemma in more generality.

\lemm
\label{fake.open.lemm}
Let $G$ be a Chevalley group scheme over $\bbZ$.
Let $T\subset G$ be a split, maximal torus over $\bbZ$, and let $B^\pm=T\ltimes U^\pm$ be Borel subgroups in $G$ over $\bbZ$ such that $B^+\cap B^-=T$. 
Then the multiplication map on strictly negative loop groups
\begin{equation}\label{fake.open.eq}
L^{--}U^-\x_\bbZ L^{--}T\x_\bbZ L^{--}U^+\to \Gr_G, \;\;\; (u^-,t,u^+)\mapsto u^-\cdot t\cdot u^+\cdot e
\end{equation}
is formally \'etale (when viewed as a map of ind-schemes).
The source of this map is called the fake open cell.
This construction is compatible with arbitrary base change $S\to \Spec(\bbZ)$, e.g., for $S=\Spec(k)$ a field.
\xlemm
\pf
The morphism $U^-\x T\x U^+\to G$ given by multiplication is an open immersion \cite[Thm.~5.1.13]{Con14}, and in particular formally \'etale \cite[04FF]{StaProj}. 
Passing to negative loop groups (and using that the $L^-$-construction commutes with products), this immediately implies that the top horizontal map 
\[
\begin{tikzpicture}[baseline=(current  bounding  box.center)]
\matrix(a)[matrix of math nodes, 
row sep=1.5em, column sep=2em, 
text height=1.5ex, text depth=0.45ex] 
{L^{-}U^-\x L^{-}T\x L^{-}U^+& L^-G  \\ 
L^{--}U^-\x L^{--}T\x L^{--}U^+ & L^{--}G, \\}; 
\path[->](a-1-1) edge node[above] {}  (a-1-2);
\path[->](a-2-1) edge node[left] {}  (a-1-1);
\path[->](a-2-1) edge node[below] {}  (a-2-2);
\path[->](a-2-2) edge node[right] {}  (a-1-2);
\end{tikzpicture}
\]
is formally \'etale. 
Here the vertical maps are the natural inclusions, and one checks that the diagram is Cartesian.
Hence, the lower horizontal arrow is formally \'etale as well.
\xpf

\rema
\label{loop.group.rem}
By the same reasoning, the induced map on loop groups $LU^-\x LT\x LU^+\to LG$ is formally \'etale as well.
\xrema

Lemma \ref{fake.open.lemm} implies that every Artinian local ring supported at the base point in $\Gr_G$ uniquely factors through the fake open cell. 
We obtain the following proposition which improves on \cite[Thm.~6.1]{PR08} in the case of split groups.

\prop
\label{reduced.prop}
Let $G$ be a split reductive group over a field $k$. Then the following are equivalent:
\begin{enumerate} 
\item The ind-scheme $LG$ is reduced \textup{(}and then even geometrically reduced\textup{)}.
\item The ind-scheme $\Gr_G$ is reduced \textup{(}and then even geometrically reduced\textup{)}.
\item The group $G$ is semisimple, and $\on{char}(k) \nmid \#\pi_1(G)$. 
\end{enumerate}
\xprop

\begin{proof}
We first show the equivalence of (1) and (2). 
Recall that the quotient map $LG\to \Gr_G$ is a (right) $L^+G$-torsor in the \'etale topology.
Thus, the ind-scheme $LG$ is \'etale locally isomorphic to $\Gr_G\x_k L^+G$.
If $LG$ is reduced, then $\Gr_G$ is reduced because $L^+G\to \Spec(k)$ is flat \cite[06QM]{StaProj}.
Conversely, if $\Gr_G$ is reduced, then $LG$ is reduced because $L^+G$ is geometrically reduced \cite[035Z]{StaProj}. 
This finishes the equivalence of (1) and (2).
Concerning geometrically reducedness, we note that if $\Gr_G$ is reduced, then it admits a presentation by Schubert varieties.
As Schubert varieties are geometrically reduced, because scheme-theoretic closure commutes with flat base change, it follows that $\Gr_G$ is reduced if and only if $\Gr_G$ is geometrically reduced. 	
Since the equivalence of (1) and (2) is valid for any field, this also implies that $LG$ is reduced if and only if $LG$ is geometrically reduced.

It remains to show the equivalence of (2) and (3) for which we may (and do) assume that $k$ is algebraically closed. 
If (3) holds, then (2) holds by \cite[Thm.~6.1]{PR08}.
Conversely, if (2) holds, i.e., if $\Gr_G$ is reduced, then $G$ is semisimple by \cite[Prop.~6.5]{PR08}.
It remains to show that $p:=\on{char}(k)$ does not divide $\pi_1(G)$. 
We may (and do) assume that $p>0$ is strictly positive. 
Let $G_\scon\to G$ be the simply connected covering.
Fix $T\subset G$, and denote by $T_\scon$ its preimage in $G_\scon$. Let $\Gr_G^0$ denote the neutral connected component of $\Gr_G$.
Then the induced map on Schubert varieties $\Gr_{G_\scon,\leq \mu}\to \Gr_{G,\leq\mu}$, $\mu\in X_*(T_\scon)$ is dominant, and hence scheme-theoretically dominant (because the target is reduced by definition).
As both ind-schemes $\Gr_{G_\scon}$, $ \Gr_G^0$ are reduced, they admit presentations by Schubert varieties indexed by dominant $\mu \in X_*(T_\scon)$. 
Thus, Lemma \ref{surjective.dist.lemm} (combined with Lemma \ref{basic.dist.lemm} (1) for the passage to ind-schemes) implies that the map
\begin{equation}\label{surj.dist.eq1}
\Dist(\Gr_{G_\scon},e)\longto \Dist(\Gr_G,e)
\end{equation}
is surjective where $e$ denotes the base point.
This map is calculated using Lemma \ref{fake.open.lemm} as follows. 
Let $B^\pm=T\ltimes U^\pm$ Borel subgroups in $G$ such that $B^+\cap B^-=T$.
Then $B^\pm_\scon=T_\scon\ltimes U^\pm$ are Borel subgroups in $G_\scon$.
By Lemma \ref{fake.open.lemm} (combined with Lemma \ref{basic.dist.lemm} (3) for the compatibility with products), the surjectivity of \eqref{surj.dist.eq1} implies the surjectivity of
\begin{equation}\label{surj.dist.eq2}
\Dist(L^{--}T_\scon)\longto \Dist(L^{--}T).
\end{equation} 
Here we use the principle that a tensor product of linear operators on possibly infinite dimensional vector spaces is surjective if and only if each linear operator is surjective.

To make the connection with $n:=\#\pi_1(G)$, recall that the kernel $Z$ of $G_\scon \to G$ is a finite $k$-group scheme of order $n$.
Clearly, the subgroup $Z$ is contained in $T_\scon$ (in fact in any maximal torus) which shows $Z=\ker(T_\scon\to T)$.
We claim that the surjectivity of \eqref{surj.dist.eq2} implies that $p\nmid n$.
We need to analyze the map $T_\scon\to T$ more carefully. 
Let $r:=\dim(T_\scon)=\dim(T)$ denote the rank of the $k$-tori.
Since $k$ is algebraically closed, passing to cocharacter lattices induces an equivalence between $k$-tori of rank $r$, and finite free $\bbZ$-modules of rank $r$.
Hence, the elementary divisor theorem implies that there exist isomorphisms $\bbG_{\on{m},k}^r\simeq T_\scon$ and $T\simeq \bbG_{\on{m},k}^r$ such that the composite
\begin{equation}\label{surj.dist.eq3}
\bbG_{\on{m},k}^r\simeq T_\scon\longto T\simeq \bbG_{\on{m},k}^r
\end{equation}
is given by $(\la_1,\ldots, \la_r)\mapsto (\la_1^{n_1},\ldots,\la_r^{n_r})$ for positive integers $n_1\geq \ldots \geq n_r\geq 1$.
We necessarily have $n=n_1\cdot \ldots\cdot n_r$.
Hence, the claim $p\nmid n$ is equivalent to the claim $p\nmid n_i$, $i=1,\ldots,r$.
Since \eqref{surj.dist.eq3} splits as a product of maps, we can apply Lemma \ref{basic.dist.lemm} (3) to see that the surjectivity of \eqref{surj.dist.eq2} implies the surjectivity of each map
\begin{equation}\label{surj.dist.eq4}
\Dist(L^{--}\bbG_{\on{m},k})\longto \Dist(L^{--}\bbG_{\on{m},k}),
\end{equation}
which is induced from $\bbG_{\on{m},k}\to \bbG_{\on{m},k}$, $\la\mapsto \la^{n_i}$ for $i=1,\ldots,r$.
Finally, Lemma \ref{last.surj.lemm} below implies $p\nmid n_i$ which finishes the proof of the proposition. 
\end{proof}

\lemm
\label{last.surj.lemm}
Let $k$ be a field of characteristic $p>0$. 
Let $n\geq 1$ be an integer, and consider the morphism of $k$-group schemes $\bbG_{\on{m},k}\to\bbG_{\on{m},k}$, $\la\mapsto \la^n$ given by taking the $n$-th power. 
If the induced morphism $\Dist(L^{--}\bbG_{\on{m},k})\to\Dist(L^{--}\bbG_{\on{m},k})$ is surjective, then $p\nmid n$.  
\xlemm
\pf
We immediately reduce to the case that $n$ is a prime number. 
For a $k$-algebra $R$, the $n$-th power map on $L^{--}\bbG_{\on{m},k}(R)$ is given by
\begin{equation}\label{power.map.eq}
1+\textstyle{\sum_{i\geq 1}}a_iu^{-i} \;\mapsto\; \big(1+\textstyle{\sum_{i\geq 1}}a_iu^{-i} \big)^n,
\end{equation}
where all $a_i\in R$ are nilpotent, and almost all $a_i$ are zero. 
The nilpotency of the $a_i$ shows that there is a presentation $L^{--}\bbG_{\on{m},k}=\on{colim}_{i\geq 1}\Spec\big(k[a_1,\ldots,a_i]/(a_1^i,\ldots,a_i^i)\big)$ where $a_i$ are viewed as formal variables.
In these coordinates, we have a canonical identification $\Dist(L^{--}\bbG_{\on{m},k})=\Dist(\bbA^{\bbN}_k,\{0\})$ where $\bbA^{\bbN}_k=\Spec(k[\{a_i\}_{i\in \bbN}])$ is the infinite-dimensional affine space.
Hence, a distribution is a $k$-linear map $\delta\co k[\{a_i\}_{i\in \bbN}]\to k$ supported at only finitely many monomials. 
We see that the space of distributions has a basis given by $\delta_{\underline r}$ such that $\underline{r}=(r_i)_{i\in \bbN}$ is a sequence of positive integers where almost all $r_i$ are zero.  
Here $\delta_{\underline r}$ takes the value $1$ on the monomial $\sqcap_{i\in \bbN}a_i^{r_i}$ and the value $0$ on all other monomials (by convention $\delta_{(0,0,\ldots)}=0$). 
We need to write down the map \eqref{power.map.eq} in the basis $\Dist(L^{--}\bbG_{\on{m},k})=\on{span}_k\{\delta_{\underline r}\;|\;\underline r\in (\bbZ_{\geq 0})^\bbN\}$. 
Suppose $n=p$ in which case we have to show that the induced map on the spaces of distributions is not surjective. 
Since $k$ has characteristic $p>0$, the formula in \eqref{power.map.eq} becomes
\[
\big(1+\textstyle{\sum_{i\geq 1}}a_iu^{-i} \big)^p=1+\textstyle{\sum_{i\geq 1}}a_i^pu^{-ip}.
\]
This means that the map on spaces of distributions is induced from $\delta_{\underline r}\mapsto \delta_{p\star\underline r}$ where $p\star\underline r\in (\bbZ_{\geq 0})^\bbN$ is the zero vector (hence $\delta_{p\star\underline r}=0$), unless $p\mid r_i$ for all entries in $\underline r=(r_i)_{i\in \bbN}$ in which case the $i$-th entry in $p\star \underline r$ is given by
\[
(p\star \underline r)_i \;=\;\begin{cases} \textstyle{{r_{{i/ p}}\over p}} & \text{if $p\mid i$;}\\ 0 & \text{else.}\end{cases}
\]
Since $p\geq 2$ the distribution $\delta_{(1,0,0,\ldots)}$ does not lie in the image of this map.
\xpf

\rema
In fact, the converse to Lemma \ref{last.surj.lemm} holds as well, i.e., for an integer $n\geq 1$ prime to $p$ the map $\bbG_{\on{m},k}\to \bbG_{\on{m},k}$, $\la \mapsto \la^n$ induces a surjection on spaces of distributions.
Indeed, one reduces to the case where $n\not= p$ is a prime number.
Then it follows from an explicit calculation --which we omit-- similarly as in the proof of Lemma \ref{last.surj.lemm}, or alternatively using affine Grassmannians as follows.
Consider the canonical map $\Gr_{\SL_n}\to \Gr_{\PGL_n}$ on affine Grassmannians, both of which are reduced by Proposition \ref{reduced.prop}.
Hence, as in \eqref{surj.dist.eq1} this induces a surjection on spaces of distributions. 
Following the proof of Proposition \ref{reduced.prop} further, we see that in \eqref{surj.dist.eq3} the elementary divisors $n_1\geq \ldots \geq n_r\geq 1$ (here $r=n$) are necessarily given by $n_1=n$ and $n_i=1$, $i\geq 2$ because $n$ is a prime number.
Now the surjectivity of the map in \eqref{surj.dist.eq4} gives the desired result.
\xrema

\subsection{Reducedness in the twisted case}
Here we give a different proof of nonreducedness of loop groups of tamely ramified semisimple groups $G$ such that $p$ divides the order of the fundamental group. The idea consists basically in observing that Weil restriction along purely inseparable extensions preserves loop groups and Grassmannians, but not flat and non-\'etale isogenies.

\prop
\label{reducedness_twisted}
Let $k$ be a perfect field of characteristic $p> 0$, and let $G$ be a tamely ramified reductive group over $F=k\rpot{t}$. 
For its loop group $LG$ to be reduced, it is necessary and sufficient that $G$ be semisimple and the order of $\pi_1(G)$ prime to $p$.
\xprop
\pf

By work of Pappas-Rapoport \cite[Thm.~6.1, Prop.~6.5]{PR08}, we only need to show that $LG$ is non-reduced whenever $G$ is semisimple and the order of the kernel $Z$ of its simply connected cover map $G_\scon \rightarrow G$ is divisible by $p>0$. 
Also, we may and do assume by étale descent that $k$ is algebraically closed.
As explained before the statement, we will consider the strictly smaller closed subgroup $\overline{G}:=\Res_{F/F^p}
G_\scon/\Res_{F/F^p}Z$ of $\Res_{F/F^p}G$, as observed in \cite[Exam.~A.7.9]{CGP15}. 
Note that Bruhat-Tits theory is available for $\overline{G}$ as well as for $\Res_{F/F^p}G$ by \cite{Lou21}, their buildings being isomorphic to the building of $G$ over $F$. 
We claim that the canonical morphism 
$$
\overline{\calG}_\bbf \rightarrow \Res_{\calO/\calO^p}{\calG}_\bbf
$$ 
between parahoric group schemes, which exists by having equivariantly identified their buildings and applying \cite[Prop.~1.7.6]{BT84}, is a locally closed immersion and its flat closure defines a normal smooth subgroup scheme whose quotient is representable by a quasi-affine group scheme, see \cite[Thm.~4A]{Ana73}.

Let $S$ be a maximal $F$-split torus of $G$, let $S_\scon$ be its unique lift to a maximal $F$-split torus of $G_\scon$ and let $\overline{S}$ be the image of $S_\scon$ in $\overline{G}$. 
We denote by $T$ (resp.~$T_\scon$, resp.~$\overline{T}$) the Cartan subgroups of $G$ (resp.~$G_\scon$, resp.~$\overline{G}$) obtained as centralizers of $S$ (resp.~$S_\scon$, resp.~$\overline{S}$). 
Arguing with big cells as in \cite[\S1.2.13, \S1.2.14]{BT84}, our claim about the canonical morphism of parahoric groups schemes will follow, once we establish that the map of connected Néron models
$$
\overline{\calT} \rightarrow \Res_{\calO/\calO^p}{\calT}
$$ 
of the Cartan subgroups is a locally closed immersion.

Assume first that $T$ is split. By using the elementary divisor theorem as in \eqref{surj.dist.eq3} above, we may assume $T=\bbG_{\on{m}}=T_\scon$ are $1$-dimensional and $Z=\Bmu_n$. 
If $n$ is prime to $p$, then $\overline{T}=\Res_{F/F^p} T$ and the claim is trivial. 
On the other hand, if $n$ is divisible by $p$, then $\overline{T}=\bbG_{\on{m}}\subseteq \Res_{F/F^p} \bbG_{\on{m}}=\Res_{F/F^p} T$ and the claim is clear as well. 

In general, let $K/F$ be a tamely ramified finite Galois extension with group $\Gamma$ splitting $T_\scon$ and $T$, and note that $K^p/F^p$ is a (pseudo-)splitting field for $\overline{T}$ with Galois group naturally isomorphic to $\Gamma$. 
Let $\calT$ (resp.\,$\overline{\calT}$, resp.\,$\calT_{\calO_K}$, resp.\,$\overline{\calT}_{\calO^p_K}$) denote the connected lft N\'{e}ron models for $T$ (resp.\,$\overline{T}$, resp.\,$T\otimes_{F}K$, resp.\,$\overline{T}\otimes_{F^p}{K^p}$).  (We warn the reader that in general $\calT_{\calO_K} \neq \calT \otimes_{\calO} \calO_K$, resp.\,$\overline{\calT}_{\calO^p_K} \neq \overline{\calT} \otimes_{\calO^p} \calO^p_K$.) 
Then we have locally closed immersions
$$
\calT \hookto \Res_{\calO_K/\calO}\calT_{\calO_K},
$$
$$
\overline{\calT} \hookto \Res_{\calO^p_{K}/\calO^p}\overline{\calT}_{\calO^p_{K}},
$$
extending the natural generic homomorphisms. 
Indeed, the maps exist by the universal property of connected N\'eron models. 
Moreover, their scheme-theoretic images are smooth by identifying them with the smooth $\Gamma$-invariants of the right hand sides, see \cite{Edi92} and compare also to \cite[Lem.~6.7]{PR08}, where we use the tameness hypothesis $p \nmid \# \Gamma $. 
Due to \cite[Prop.~10.1.4]{BLR12}, the resulting morphisms must be locally closed immersions. 
Taking restrictions of scalars along $\calO/\calO^p$ of the first map, and along $\calO^p_K/\calO^p$ of the split case morphism $\overline{\calT}_{\calO^p_K} \hookto \Res_{\calO_K/\calO^p_K} \calT_{\calO_K}$, we deduce the general claim. 

As a consequence of the group-theoretic facts just established, we derive that 
$$
\Gr_{\overline{\calG}_\bbf}^0 \rightarrow \Gr_{\Res_{\calO/\calO^p}{\calG}_\bbf}^0\cong \Gr_{\calG_\bbf}^0
$$
is a closed immersion. Indeed, if we let $\overline{\calG}_\bbf^1$ denote the flat closure of $\overline{\calG}_\bbf$ inside $\Res_{\calO/\calO^p}{\calG}_\bbf$ (see the definition of flat closure in $\S\ref{ind.flatness.sec}$ below), then the morphism
$$
\Gr_{\overline{\calG}_\bbf^1} \rightarrow  \Gr_{\calG_\bbf}
$$
is a quasi-compact immersion by \cite[Prop.~1.2.6]{Zhu16}, which must be closed, as source and target are ind-projective by \cite[Thm.~5.2]{Lou21}. 
Finally, we have to show that the Galois cover 
$$\Gr_{\overline{\calG}_\bbf} \rightarrow \Gr_{\overline{\calG}_\bbf^1} $$
with group $\overline{\calG}_\bbf^1(\calO^p)/\overline{\calG}_\bbf(\calO^p)$ induces an isomorphism between neutral components, which can be checked at the level of $k$-points. In other words, we must show that every element of $\overline{G}(\calO^p)$ stabilizing $\bbf$ which does not lie in the parahoric subgroup (i.e. ``the connected stabilizer'') maps to a connected component of the affine Grassmannian different from that of the identity. In the case of reductive groups, this is the main result of \cite{HR08}, stating that parahorics are the intersection of the stabilizers with the kernel of the Kottwitz map $\kappa_G$.
For pseudo-reductive groups, this was proved in \cite[Prop.~3.9, Thm.~5.2]{Lou21}.

If $LG$ were reduced, then the closed immersion $\Gr_{\overline{\calG}_\bbf}^0 \rightarrow  \Gr_{\calG_\bbf}^0$ would have to be an isomorphism, because $\Gr_{\calG_{\scon,\bbf}}^0\rightarrow \Gr_{\calG_\bbf}^0$ is a universal homeomorphism. 
In particular, their Lie algebras would be the same and via the uniformization $\Gr_{\calG_\bbf}=LG/L^+\calG_\bbf$ (similarly for $\overline{G}$), this would imply that the $F$-vector space $\text{Lie}\,G$ is the (non-direct) sum of the $F^p$-subspace $\text{Lie}\,\overline{G}$ and the $\calO$-lattice $\text{Lie}\,\calG_\bbf$. 
But, the dimension of $\overline{G}$ is strictly smaller than that of $\Res_{F/F^p}G$ by construction, so this is obviously a contradiction.
\xpf

\rema One would hope that a similar statement holds beyond the tamely ramified case, but one cannot control the N\'eron models with the same ease. On the other hand, if one tried to classify reducedness of the loop group for the more general class of pseudo-reductive groups, the above argument suggests this could be very difficult.
\xrema

\section{Ind-flatness}
In this section $G$ will denote a Chevalley group scheme over $\bbZ$.
Our aim is to prove in Proposition \ref{ind.flat.prop} that its affine Grassmannian $\Gr_{G,\bbZ}$ (equivalently, its loop group) is ind-flat over $\bbZ$ in the sense of Definition \ref{ind.flat.defi}. 
In Proposition \ref{ind.flat.prop.twisted} we explain how to generalize our proof to include the case of tamely ramified twisted groups.

\subsection{Preliminaries on ind-flatness}\label{ind.flatness.sec} 
Recall our conventions on ind-schemes, see Section \ref{conventions.sec}. 

\defi
\label{ind.flat.defi}
Let $S$ be a scheme.
An $S$-ind-scheme $X$ is called ind-flat if there exists a presentation $X=\on{colim}X_i$ where $X_i$ are flat $S$-schemes via the map $X_i\subset X\to S$.
\xdefi

Now let $R$ be a Dedekind ring with fraction field $K$. For an $R$-scheme $X$, the flat closure $X^{\on{fl}}$ is the scheme theoretic image of the inclusion $X_K\subset X$. Since $X_K\subset X$ is a quasi-compact map, the scheme theoretic image commutes with localization \cite[01R8]{StaProj}, and the closed immersion $X^{\on{fl}}\hookto X$ is an isomorphism on generic fibers. Then the scheme $X$ is flat over $R$ if and only if the map $X^{\on{fl}}\hookto X$ is an isomorphism if and only if $\calO_X$ is $R$-torsion-free. If $\varphi\co X\to Y$ is a map of $R$-schemes, then there is a map $\varphi^{\on{fl}}\co X^{\on{fl}}\to Y^{\on{fl}}$ with $\varphi_K=(\varphi^{\on{fl}})_K$.

\lemm \label{flat.closure.lemm}
Let $R$ be a Dedekind ring with fraction field $K$. For an $R$-ind-scheme $X$ the following conditions are equivalent:
\begin{enumerate}
	\item $X$ is ind-flat;
	\item for every presentation $X=\on{colim}X_i$, the map $\on{colim}X_i^{\on{fl}}\hookto \on{colim} X_i$ is an isomorphism;
	\item every ind-(closed immersion) $Y\hookto X$ which induces $Y_K\cong X_K$ is an isomorphism. 
\end{enumerate}
\xlemm
\pf The implications (3) $\Rightarrow$ (2) $\Rightarrow$ (1) are immediate, and we prove (1) $\Rightarrow$ (3). 

Let $X=\on{colim}X_i$ be a flat presentation.
Let $Y\hookto X$ be an ind-(closed immersion) which induces $Y_K\cong X_K$.
For each $i$, the induced map $Y\cap X_i\hookto X_i$ is an ind-(closed immersion) which induces $(Y\cap X_i)_K\cong (X_i)_K$.
We want to show that $Y\cap X_i \cong X_i$.
Replacing $X$ by $X_i$, we may assume that $X$ is a flat $R$-scheme. 
Covering $X$ by open affine schemes, we may further assume that $X$ is affine, hence quasi-compact.
Now let $Y=\on{colim}Y_j$ be any presentation. 
We will show that $Y_j\cong X$ for $j>\!\!>0$.
As $Y_K\cong X_K$ on generic fibers and $X_K$ is quasi-compact, there is a $j$ with $X_{K}\hookto Y_{j,K}$ so that $Y_{j,K}\cong X_K$ {(see Section \ref{conventions.sec})}. 
As $Y_j\hookto X$ is a closed immersion and $X$ is $R$-flat, we must have $Y_j\cong X$.
\xpf

\defi
\label{flat.closure.defi}
For an ind-scheme $X=\on{colim}_iX_i$, the {\it flat closure $X^{\on{fl}}$} is the ind-scheme $X^{\on{fl}}=\on{colim}_iX_i^{\on{fl}}$. 
\xdefi

In view of Lemma \ref{flat.closure.lemm}, the ind-(closed immersion) $X^{\on{fl}}\subset X$ is well-defined independently of the choice of a presentation.
Also a map of $R$-ind-schemes $X\to Y$ induces a map $X^{\on{fl}}\to Y^{\on{fl}}$ on the flat closures.

\subsection{Ind-flatness of affine Grassmannians}
The starting point is the following lemma.

\lemm
\label{elementary.flat.lemm}
Let $H$ be a smooth, affine group scheme over $\bbZ$.
\begin{enumerate}
\item The positive loop group $L^+H\to \Spec(\bbZ)$ is a flat, affine group scheme.
\item If $H$ is split unipotent or a split torus, then both the loop group $LH$ and the strictly negative loop group $L^{--}H$ are ind-flat over $\bbZ$. 
\end{enumerate}
\xlemm
\pf
For (1), let  $L^+_iH$, $i\geq 0$, be the smooth, affine $\bbZ$-group scheme defined by the functor $L^+_iH(R)=H\big(R[u]/(u^{i+1})\big)$ for a ring $R$.
Then $\{L^+_iH\}_{i\geq 0}$ naturally forms an inverse system, and the canonical map $L^+H\to \lim_{i\geq 0}L^+_iH$ is an isomorphism.
This implies (1).

For (2), observe that the map 
\begin{equation}\label{open.immersion.eq}
L^{--}H\x_\bbZ L^+H\to LH, \;\;(h^-,h^+)\mapsto h^-\cdot h^+,
\end{equation}
is representable by an open immersion.
Since $L^+H$ is faithfully flat by (1), the ind-flatness of $L^{--}H$ follows from the one for $LH$.
Now let $H$ be a split unipotent group scheme.
Then $H\simeq \bbA^n_\bbZ$ as schemes for some $n\geq 0$.
Since the formation of loop spaces commutes with products, it is enough to show that $L\bbA^1_\bbZ$ is ind-flat.
This is immediate from the identification, which is functorial in the ring $R$,
\[
L\bbA^1_\bbZ(R)\;=\; R\rpot{u}\;=\; \underset{{i>\!\!>-\infty}}{\on{colim}} \bigsqcap_{j\geq i}\bbA^1_\bbZ(R),
\]
given by mapping a Laurent series $\sum a_iu^i$ to the vector $(a_i)$.
Next let $H$ be a split torus. Then $H\simeq \bbG_{\on{m},\bbZ}^n$ as (group) schemes for some $n\geq 0$,
and we reduce to the case $T=\bbG_{\on{m},\bbZ}$.
In this case, the map \eqref{open.immersion.eq} is surjective and hence an isomorphism. 
We see that the ind-flatness of $L\bbG_{\on{m},\bbZ}$ is equivalent to the one of $L^{--}\bbG_{\on{m},\bbZ}$.
For the latter we note that for any ring $R$,
\[
L^{--}\bbG_{\on{m},\bbZ}(R)\;=\; \big(1+u^{-1}R[u^{-1}]\big)^\x\;=\;\big\{ 1+\textstyle\sum_{i\geq 1}a_iu^{-i}\;|\;\text{$a_i\in R$ nilpotent}\big\},
\]
so that $L^{--}\bbG_{\on{m},\bbZ}\simeq \on{colim}_{i\geq 1}\Spec\big(\bbZ[a_1,\ldots,a_i]/(a_1^i,\ldots,a_i^i)\big)$. 
This is clearly ind-flat.
\xpf

Recall that $G$ denotes a Chevalley group scheme over $\bbZ$.

\coro
\label{fake.open.flat.coro}
Let $T\subset G$ be a split, maximal torus over $\bbZ$, and let $B^\pm=T\ltimes U^\pm$ be Borel subgroups in $G$ such that $B^+\cap B^-=T$. 
Then the fake open cell (cf.~Lemma \ref{fake.open.lemm})
\[
L^{--}U^-\x_\bbZ L^{--}T\x_\bbZ L^{--}U^+\to \Spec(\bbZ)
\]
is ind-flat.
\xcoro
\pf
Since $U^\pm$ are split unipotent and $T$ is a split torus, this is immediate from Lemma \ref{elementary.flat.lemm} (2).
\xpf

The ind-flatness of the affine Grassmannian is deduced from Corollary \ref{fake.open.flat.coro} using the following observation due to Faltings \cite[Proof of Cor.~11]{Fal03}.

\begin{lem}
\label{what.faltings.knows}
Let $Y \hookto X$ be an ind-(closed immersion) of ind-(locally Noetherian) ind-schemes. 
If for every local Artinian ring $R$ the induced map 
\begin{equation}\label{bijective.faltings.eq}
Y(R)\rightarrow X(R)
\end{equation}
is bijective, then $Y\hookto X$ is an isomorphism.
\end{lem}

\begin{proof}
By Yoneda's lemma, we may view $X$, $Y$ as set-valued (contravariant) functors on the category of Noetherian, affine schemes.  
For any such $T$, the induced map $Y(T)\to X(T)$ is clearly injective, and we need to show the surjectivity.
This can be checked after base change $Y\x_XT\to T$ as our assumptions are stable under base change.
We reduce to the case where $X=T$ is a Noetherian (affine) scheme.
Write $Y=\on{colim}_iY_i$ as a filtered colimit of closed subschemes of $X$.
We claim that this sequence stabilizes with value $X$.

Let $\calI_i\subset \calO_X$ be the ideal sheaf defining $Y_i$.
Since the index set $I$ is filtered, it is enough to show the existence of an index $i$ with $\calI_{i}=0$.
For this we note that $\calI_i=0$ if and only if the annihilator ideal sheaf $\on{Ann}_{\calO_X}(\calI_i)$ is equal $\calO_X$, or equivalently if the closed subscheme $Z_i\subset X$ defined by $\on{Ann}_{\calO_X}(\calI_i)$ is empty.
For $i\leq j$, we have
\[
Y_i\subset Y_j \iff \calI_i\supset \calI_j \Longrightarrow \on{Ann}_{\calO_X}(\calI_i)\subset \on{Ann}_{\calO_X}(\calI_j) \iff Z_i\supset Z_j.
\]
Since $X$ is noetherian and the set $I$ is filtered, there is an $i_0 \in I$ such that $Z_{i_0}$ is a minimum, i.e. $Z_{i_0} \subseteq Z_i$ for all $i \in I$. Now suppose for the sake of contradiction that $Z_{i_0}\neq\varnothing$.

Now let $\eta\in Z_{i_0}$ be a generic point of an irreducible component.
It remains to find an index $i_1>i_0$ such that $\eta\notin Z_{i_1}$. 
The property defining $\eta$ means that the closed subscheme defined by $\on{Ann}_{\calO_{X}}(\calI_{i_0})_\eta=\on{Ann}_{\calO_{X,\eta}}(\calI_{i_0,\eta})$ in $\Spec(\calO_{X,\eta})$ is supported at the closed point $\frakm_\eta$, where $\frakm_\eta\subset\calO_{X,\eta}$ denotes the maximal ideal.
Since $\calO_{X,\eta}$ is Noetherian, this is equivalent to the existence of some integer $N>\!\!>0$ such that 
\[
\frakm_\eta\;\supset\; \on{Ann}_{\calO_{X,\eta}}(\calI_{i_0,\eta}) \;\supset\; \frakm_\eta^N.
\]
Hence, $\calI_{i_0,\eta}$ is a finitely generated module over the Artinian ring $\calO_{X,\eta}/\frakm_\eta^N$, and therefore Artinian itself.
Since our index set is filtered, we can choose $i_1>i_0$ such that $\calI_{i_1,\eta}$ is the minimum among the set $\{\calI_{j,\eta}\}_{j\geq i_0}$.
 Now  \eqref{bijective.faltings.eq} applied to the Artinian rings $\calO_{X,\eta}/\frakm_\eta^n$, $n\geq 1$, shows that $Y_{i_1}(\calO_{X,\eta}/\frakm^n_\eta) \overset{\sim}{\rightarrow} X(\calO_{X,\eta}/\frakm^n_\eta)$, that is, every homomorphism $\calO_{X, \eta} \rightarrow \mathcal O_{X, \eta}/\frakm^n_\eta$ factors uniquely through $\calO_{X,\eta} \rightarrow \calO_{X,\eta}/\calI_{i_1,\eta}$.  But then $\calI_{i_1, \eta}
\subseteq \cap_{n\geq 0}\frakm_\eta^n$, so that $\calI_{i_1, \eta} = 0 $ by Krull's intersection theorem. This implies that $\eta \notin Z_{i_1}$. This finishes the proof of the lemma.
\end{proof}

\rema 
\label{Faltings.Henselian.rema}
The proof of Lemma \ref{what.faltings.knows} shows that condition \eqref{bijective.faltings.eq} can be weakened. 
Namely, it is enough to use local Artinian rings which are strictly Henselian.
Indeed, in the last part of the proof it is enough to show $\calI_{i_1,\bar\eta}=0$ where $\bar\eta \to \eta$ is a geometric point and $\calI_{i_1,\bar\eta}$ denotes the stalk on the \'etale site.
\xrema

\prop
\label{ind.flat.prop}
The affine Grassmannian $\Gr_{G, \bbZ}$ is an ind-flat ind-scheme over $\bbZ$. In particular, the ind-scheme $\Gr_{G, \bbZ}$ is reduced if and only if $G$ is semisimple.
\xprop
\pf
Let $\Gr_{G,\bbZ}^{\on{fl}}\subset \Gr_{G, \bbZ}$ be the flat closure, cf.~Definition \ref{flat.closure.defi}.
By Lemma \ref{what.faltings.knows} and Remark \ref{Faltings.Henselian.rema}, it is enough to show that every local Artinian, strictly Henselian point $g\co \Spec(R)\to\Gr_{G, \bbZ}$ factors through $\Gr_{G,\bbZ}^{\on{fl}}$. 
Let $k$ be the residue field of $R$, and denote by $\bar g\co \Spec(k)\to \Gr_{G,\bbZ}$ the reduction of $g$.
Fix a split maximal torus $T\subset G$ over $\bbZ$, and $B^\pm=T\ltimes U^\pm$ as in Corollary \ref{fake.open.flat.coro}.
By the Cartan decomposition $\Gr_G(k)=\sqcup_{\mu\in X_*(T)}L^+G(k)\cdot u^\mu\cdot e$ (use that $k$ is separably closed)\footnote{Alternatively, the reader can note that for tamely ramified groups, with the help of our big cell results and with the rationality results of decompositions in \cite{HR08}, we know the Cartan decomposition holds for any field $k$ which is a ${\mathbb W}$-algebra, so Remark \ref{Faltings.Henselian.rema} is not really needed, even in the proof for the generalization proved in Proposition \ref{ind.flat.prop.twisted}.}, we can write $\bar g$ as a product $\bar{h}\cdot u^\mu\cdot e$ for some $\bar h \in L^+G(k)$, $\mu\in X_*(T)$.
By formal smoothness of $L^+G\to \Spec(\bbZ)$, we can lift $\bar{h}$ to an $R$-valued point $h\co \Spec(R)\to L^+G$.
Since $L^+G\to \Spec(\bbZ)$ is flat and $LT\to \bbZ$ is ind-flat by Lemma \ref{elementary.flat.lemm}, the inclusion $\Gr_{G,\bbZ}^{\on{fl}}\subset \Gr_{G, \bbZ}$ is invariant under the left action of $L^+G$ and $LT$.
Replacing $g$ by $u^{-\mu}\cdot h^{-1}\cdot g$, we may therefore assume that $g$ is supported at the base point.
Then $g$ factors through the fake open cell $L^{--}U^-\x L^{--}T\x L^{--}U^+$ by Lemma \ref{fake.open.lemm}.
Since this is ind-flat by Corollary \ref{fake.open.flat.coro}, the map \eqref{fake.open.eq} factors through $\Gr_{G,\bbZ}^{\on{fl}}$. This shows $\Gr_{G,\bbZ}^{\on{fl}}=\Gr_{G,\bbZ}$.

For the second assertion, note that $\Gr_{G,\bbQ}$ is reduced if and only if $G$ is semisimple, see \cite[Thm.~6.1, Prop.~6.5]{PR08}. 
Hence, $\Gr_{G,\bbZ}$ is not reduced whenever $G$ is not semisimple. 
Conversely, if $G$ is semisimple, then by taking the flat closure of any reduced presentation of $\Gr_{G,\bbQ}$ gives a reduced presentation of $\Gr_{G,\bbZ}^{\on{fl}}=\Gr_{G,\bbZ}$.

This finishes the proof of the proposition.
\xpf

We now briefly generalize this to the parahoric group schemes over ${\mathbb W}\pot{t}$ constructed previously.

\prop
\label{ind.flat.prop.twisted}
Let $G$ be a tamely ramified reductive $k\rpot{t}$-group.
The affine flag variety $\uFl_{G,{\bbf}}$ is an ind-flat ind-scheme over ${\mathbb W}$.
In particular, the ind-scheme $\uFl_{G,{\bbf}}$ is reduced if and only if $G$ is semisimple.
\xprop
\pf

 Without loss of generality, we may and do assume that $\bbf={\bf 0}$ is an absolutely special vertex. 
 We follow the reasoning of Proposition \ref{ind.flat.prop}:
After translating, it suffices to show that, for every Artinian $\bbW$-algebra $R$, every $R$-valued point of $L^{--}\ucG_{\bf 0}$ supported at the identity lies in the ind-flat closure, see Lemma \ref{what.faltings.knows}. 
Such points uniquely lift to the strictly negative loop group attached to the twisted open cell
$$\ucC_{\bf 0}=\ucU_{\bf 0}^-\times \ucT\times \ucU_{\bf 0}^+:=(\Res_{\bbW[u^{-1}]/\bbW[t^{-1}]} \,U^-_H\times T_H\times U_H^+)^{\sigma,\circ}$$
because $L^{--}\ucC_{\bf 0} \to  L^{--}\ucG_{\bf 0}$ is formally étale.
Hence, it suffices to see that $L^{--}\ucC_{\bf 0}\to \Spec(\bbZ)$ is ind-flat. 
The unipotent part of $\ucC_{\bf 0}$ is identified, as a scheme, with a product of restrictions of scalars of affine spaces, so its associated strictly negative loop group is ind-flat over $\bbZ$, compare with the argument in Lemma \ref{elementary.flat.lemm} (2).
As for the connected Néron model $\ucT$, we choose a smooth surjection $\ucT_1\to \ucT$ from the connected Néron model $\ucT_1$ of an induced torus $T_1$; see \cite[Lem.~6.7]{PR08} for the surjectivity assertion. 
We claim that $L^{--}\ucT_1(R) \to L^{--}\ucT(R)$ is surjective for every Artinian $\bbW$-algebra $R$. 
Since both functors are formally smooth, this reduces, by lifting across square zero nilpotent thickenings, to the surjectivity of Lie algebras $\on{Lie}L^{--}\ucT_1 \to \on{Lie}L^{--}\ucT$ viewed as $\bbW$-modules.
The latter is the map underlying the map of $\bbW[t^{-1}]$-modules $\on{Lie}\ucT_1\to \on{Lie}\ucT$ which is surjective by smoothness of $\ucT_1\to \ucT$.
As in Lemma \ref{elementary.flat.lemm} (2), it is easy to see that $L^{--}\ucT_1\to \Spec(\bbZ)$ is ind-flat because $T_1$ is an induced torus. 
This implies the ind-flatness of $L^{--}\ucT\to \Spec(\bbZ)$, again using Lemma \ref{what.faltings.knows}.
The final assertion on reducedness follows as in the proof of Proposition \ref{ind.flat.prop}.

\xpf

\section{Consequences for Pappas-Zhu local models}\label{Consequences.for.local.models}

We discuss some consequences of our findings for the theory of local models in cases where $p$ divides the order of $\pi_1(G_\der)$. 

In this final section, let $F$ be a discretely valued, complete field of characteristic $0$ with algebraically closed residue field $k$ of characteristic $p>0$. 
We fix a triple $(G,\{\mu\},\calG_\bbf)$ where $G$ is a tamely ramified reductive $F$-group, $\{\mu\}$ a conjugacy class of geometric cocharacters defined over a finite extension $E/F$, and $\calG_\bbf$ is a parahoric $\calO_F$-group scheme with generic fiber $G$. This notation seems to have first appeared in the survey article of Pappas-Rapoport-Smithling, see \cite{PRS13}, and first termed LM triple by He-Pappas-Rapoport in \cite[\S2.1]{HPR20} (but always under the assumption that $\{\mu\}$ is minuscule).
Pappas-Zhu \cite{PZ13} construct from the data $(G,\{\mu\},\calG_\bbf)$ the \emph{Pappas-Zhu local model}
\[
\bbM\;=\;\bbM{(G,\{\mu\},\calG_\bbf)},
\]
which is a flat, projective $\calO_E$-scheme equipped with a left action of a smooth affine group scheme. 
Recall that the construction of $\bbM$ requires the construction of a parahoric $\calO_F[t]$-group scheme $\ucG_\bbf$ in the sense of \cite[Thm.~4.1]{PZ13} which lifts $\calG_\bbf$ along the specialization $t\mapsto \varpi$ for some fixed uniformizer $\varpi\in \calO_F$. 
In particular, $\bbM$ depends a priori on
certain auxiliary choices, but it is shown in \cite[Thm.~2.7]{HPR20} that $\bbM$ actually depends, up to equivariant isomorphism, only on the data $(G,\{\mu\},\calG_\bbf)$.
The generic fiber $\bbM\otimes E$ is naturally the Schubert variety in the affine Grassmannian of $G/F$ associated with the class $\{\mu\}$. 
The special fiber is equidimensional, but not irreducible in general, and is equipped with a closed embedding
\begin{equation}\label{closed.embedding.special}
\bbM\otimes k \;\hookto\;\Fl_{G^\flat,\bbf^\flat}.
\end{equation}
The pair $(G^\flat,\bbf^\flat)$ is an equal characteristic analogue over a local function field $F^\flat=k\rpot{t}$ of the pair $(G,\bbf)$, see \cite{PZ13} (see also \cite{HR20a}).

More precisely, fix a maximal $F$-split torus $S$ whose apartment $\scrA(G,S,F)$ contains the facet $\bbf$.
Its centralizer $T$ is a maximal torus since $G$ is quasi-split by Steinberg's theorem. 
We also have the corresponding data $S^\flat \subset T^\flat$ inside the equal characteristic analogue $G^\flat$.
There is an identification of apartments $\scrA:=\scrA(G,S,F)=\scrA(G^\flat,S^\flat,F^\flat)$ compatible with the action of the Iwahori Weyl groups $W:=W(G,S,F)=W(G^\flat,S^\flat,F^\flat)$ under which $\bbf=\bbf^\flat$.
Then under \eqref{closed.embedding.special} the reduced locus of $\bbM\otimes k$ identifies by \cite[Thm.~6.12]{HR21} with the admissible locus
\[
\calA(G,\{\mu\},\calG_{\bbf})\;\defined\; \bigcup_{w} S_w(\bbf,\bbf),
\]
where the union is taken over the finitely many elements $w$ of the admissible set $W_\bbf\Adm(\{\mu\})W_\bbf$ in $W$, and where $S_w(\bbf,\bbf)$ denotes the Schubert variety arising from a $\calG_{\bbf}$-orbit in $\Fl_{G^\flat,\bbf}$ using $\bbf=\bbf^\flat$. Here ${\rm Adm}(\{\mu\})$ is the so-called $\mu$-{\em admissible set}, and can be given a purely combinatorial definition; see e.g.\,\cite{PRS13, HH17}.
The following result is an application of Proposition \ref{permanence}, which is used in \cite[Rem.~2.2]{HR22}. 

\begin{propo} \label{minuscule}
Let $\bba\subset \scrA$ be an alcove whose closure contains $\bbf$.
Suppose $\{\mu\} \subset X_*(T)$ has the property that $S_v(\bba, \bba)$ is normal for all \textup{(}equivalently, for the maximal elements\textup{)} $v \in {\rm Adm}(\{\mu\}) W_{\bbf}$. Then all $(\bbf, \bbf)$-Schubert varieties in $\mathcal A(G, \{ \mu\}, \calG_{\bbf})$ are normal. In particular, this last conclusion holds when $\bar{\mu} \in X_*(T)_I$ is minuscule for the \'{e}chelonnage roots and the closure of $\bbf$ contains a special vertex.
\end{propo}

\begin{proof}
Let $v \in W_{\bbf} {\rm Adm}(\{\mu\}) W_{\bbf}$; to prove the first assertion we need to show that $S_{v}(\bbf, \bbf)$ is normal. By Corollary \ref{normality.cor}, we may assume $v$ is left $\bbf$-maximal, so that $S_v(\bbf, \bbf) = S_v(\bba, \bbf)$ and so we only need to consider $({\bf a}, {\bf f})$-Schubert varieties. By \cite[Thm.~1.3]{HH17}, we have $W_{\bbf} {\rm Adm}(\{\mu\}) W_{\bbf} = {\rm Adm}(\{\mu\}) W_{\bbf}$, so we may reduce ourselves to proving that $S_v(\bba, \bbf)$ is normal for any $v \in {\rm Adm}(\{\mu\}) W_{\bbf}$. By Proposition \ref{permanence}, it is enough to prove that $S_{v \eta}(\bba, \bba)$ is normal for all $\eta \in W_{\bbf}$.  But these are normal by assumption.

For the second assertion, choose a special vertex ${\mathbf 0} \in \bar{\bbf}$. Since $\bar \mu$ is minuscule, the Schubert variety $S_{\bar{\mu}}(\bba, {\mathbf 0})$ is smooth, and hence so is its preimage under $\Fl_{\bba} \to \Fl_{{\mathbf 0}}$. This is itself a Schubert variety $S_{v_0}(\bba, \bba)$ indexed by the unique longest element $v_0 \in W_{{\mathbf 0}} {\rm Adm}(\{\mu\}) W_{{\mathbf 0}}$, a set which contains ${\rm Adm}(\{\mu\}) W_\bbf$; thanks again to Proposition \ref{permanence} the latter set indexes normal Schubert varieties, and then we are done by the first assertion.
\end{proof}

As an application we obtain Corollary \ref{local.models.corollary.intro} from the introduction:

\coro\label{local.models.corollary}
Assume $p$ divides the order of $\pi_1(G_{\on{der}})$. 
\begin{enumerate}
\item If every Schubert variety in the admissible locus $\calA(G, \{\mu\},\calG_{\bbf})$ is normal, then $\bbM$ is normal and its special fiber is reduced. This is the case when $\bar\mu$ is minuscule for the \'echelonnage roots and $\bbf$ contains a special vertex in its closure.
\item If any Schubert variety inside the admissible locus $\calA(G, \{\mu\},\calG_{\bbf})$ is not normal, then $\bbM$ is not normal and its special fiber is not reduced.
\end{enumerate}
\xcoro
\pf
Part (1) is immediate from Proposition \ref{minuscule} and \cite[Thm.~2.1]{HR22}. 
For (2), suppose one of the Schubert varieties inside $\calA(G, \{\mu\},\calG_{\bbf})$ is not normal.
Then the irreducible component containing this Schubert variety is not normal as well by Corollary \ref{normality.cor}.  
The normalization morphism 
\begin{equation}\label{normal.eq.map}
p\co \tilde\bbM\longto \bbM
\end{equation}
identifies by \cite[Cor.~2.3, Cor.~2.5]{HR22} with the map from the Pappas-Zhu local model of some $z$-extension of $G$.
In particular, \eqref{normal.eq.map} is a finite, birational, universal homeomorphism and an isomorphism on generic fibers; recall that $\bbM\otimes E$ is a Schubert variety in characteristic $0$.
This already shows that $\bbM$ is not normal, see \cite[Rem.~2.4]{HR22}.
It remains to show that the special fiber $\bbM \otimes k$ is not reduced.
Arguing by contradiction let us assume that $\bbM \otimes k$ is reduced.
For any line bundle $\calL$ on $\bbM $, this implies the injectivity of the canonical map $H^0(\bbM \otimes k,\calL)\to H^0(\tilde\bbM\otimes k,\calL)$.
Furthermore, if $\calL$ is ample and $N>0$ sufficiently large, then there is an equality 
$$
\dim_k H^0(\bbM\otimes k, \calL^N)= \dim_k H^0(\tilde\bbM\otimes k, \calL^N)
$$
by transporting the claim to the generic fiber using flatness. 
Hence, we get an isomorphism of vector spaces $H^0(\bbM\otimes k, \calL^N)\cong H^0(\tilde\bbM\otimes k, \calL^N)$, and thus \eqref{normal.eq.map} must be an isomorphism on special fibers.
We arrive at the desired contradiction since the irreducible components of $\tilde \bbM\otimes k$ are normal by \cite[Cor.~2.5]{HR22} for example.
\xpf

Let us give two concrete examples of badly behaved PZ local models. 
The obvious class of examples arises from Weil restrictions of scalars along ramified extensions:

\exam \label{Weil.restriction.example}
Let $F'/F$ be a totally ramified extension of $2$-adic fields of odd degree $e\geq 1$. 
Consider the Weil restriction of scalars $G=\Res_{F'/F}(\PGL_2)$, and let $\{\mu\}$ be the unique (nonzero) minuscule conjugacy class defined over $F$.
As parahoric subgroup we take the pointwise fixer of the standard lattice $\calO_{F'}^2$, that is, the associated parahoric group scheme is $\calG=\Res_{\calO_{F'}/\calO_F}(\PGL_2)$.
It corresponds to an absolutely special vertex $\bf 0$, and hence the special fiber of the PZ local model is irreducible.
Its underlying reduced subscheme is the unique $e$-dimensional $(\bf 0, \bf 0)$-Schubert variety in the affine Grassmannian for $\PGL_2$ in characteristic $2$. 
If $e\geq 2$ this Schubert variety is not normal by our classification in Corollary \ref{classification.pgl2.grassmannian}.
Hence, the special fiber of the PZ local model is not reduced in this case by Corollary \ref{local.models.corollary} (2).

We remark that if we drop the tameness assumption and take $e=2$, then $F'/F$ is wildly ramified and we can invoke \cite{Lev16} to define the local models.  
Again the special fiber of such a local model is not reduced because the corresponding Schubert variety is not normal, as follows immediately from Corollary \ref{qmin_in_Gr}.
Indeed, the corresponding Schubert variety is the quasi-minuscule one for the group $\PGL_2$ (although it is more natural to think of it as a Schubert variety attached to the standard pseudo-reductive group $\Res_{k\rpot{u}/k\rpot{t}}\PGL_2$ where $u^2 = t$, compare this to the approach of \cite{Lou19,Lou20}).
\xexam 

A less obvious example is given by ramified unitary groups. 
In this case, the underlying group is even absolutely simple:

\exam \label{unitary.example}
Let $F'/F$ be a totally ramified quadratic extension of $3$-adic fields. 
Let $G=\on{PU}_3(F'/F)$ be the adjoint, quasi-split unitary group associated with the Hermitian form $x_1\bar x_3+x_2\bar x_2 +x_3\bar x_1$ on $F'^3$.
Let $\{\mu\}$ be the minuscule conjugacy class corresponding to the coweight $(1,0,0)$. 
As parahoric subgroup we take the pointwise fixer of the standard lattice $\calO_{F'}^3$ which corresponds to an absolutely special vertex ${\bf 0}$ of the building, see \cite[\S7]{HR20a}. 
In this case, the special fiber is again irreducible and its underlying reduced locus is the unique $2$-dimensional Schubert variety in the twisted affine Grassmannian for $\on{PU}_3(F'/F)$ in characteristic $3$, that is, the quasi-minuscule one. 
This Schubert variety is not normal by Corollary \ref{qmin_in_Gr} so that the special fiber of the PZ local model is not reduced by Corollary \ref{local.models.corollary} (2).
\xexam

Next we comment on the behavior of PZ local models relatively to central extensions. The adjoint quotient $(G,\{\mu\}) \rightarrow (G_{\on{ad}},\{\mu_{\on{ad}}\})$ induces a natural morphism
$$\on{ad}_*\co \bbM \rightarrow \bbM_{\on{ad}} \otimes_{\calO_{E_{\on{ad}}}} \calO_E, $$
where $E_{\on{ad}}$ is the reflex field of $(G_{\on{ad}},\{\mu_{\on{ad}}\})$ and $\bbM_{\on{ad}}:=\bbM(G_{\on{ad}},\{\mu_{\on{ad}}\}, \calG_{\bbf,\on{ad}})$. This is a (fiberwise) birational universal homeomorphism, but not always an isomorphism. We are now going to look at the category consisting of all LM triples centrally lifting $(G_{\on{ad}},\{\mu_{\on{ad}}\}, \calG_{\bbf,\on{ad}})$, endowed with the obvious morphisms. It admits fiber products and we use this to study the variation of the PZ local models along central lifts.

\prop\label{local.model.central}
Let $(G_i,\{\mu_i\}, \calG_{i,\bbf})$, $i=1,2$, be two LM central lifts of $(G_{\on{ad}},\{\mu_{\on{ad}}\}, \calG_{\bbf,\on{ad}})$ and denote by $(G_3,\{\mu_3\}, \calG_{3,\bbf})$ their fiber product. 
For $i=1,2,3$, let $\bbM_i$ be the PZ local model for $(G_i,\{\mu_i\}, \calG_{i,\bbf})$, and denote by $\on{ad}_{i,*}\colon \bbM_i\to \bbM_\ad\otimes_{\calO_{E_{\on{ad}}}}\calO_{E_i}$ the induced map.
If $p \nmid \frac {\#\pi_1(G_{1,\on{der}}) }{\#\pi_1(G_{3,\on{der}})}$, then the rational map $ \on{ad}_{2,*}\circ \on{ad}_{1,*}^{-1}$ extends to an actual morphism of schemes over $\calO_{E_3}$.
\xprop
\pf
Recall that, by construction of parahoric $\calO_F[t]$-group schemes extensions, the natural maps $\calG_{3,\bbf}\to \calG_{i, \bbf}$ can be extended to $\ucG_{3,\bbf}\to \ucG_{i, \bbf}$ inducing morphisms
\begin{equation}\label{local.model.central.eq}
 \bbM_3 \rightarrow \bbM_i\otimes_{\calO_{E_i}} \calO_{E_3}
 \end{equation}
for $i=1,2$ over $\bbM_\ad\otimes_{\calO_{E_{\on{ad}}}}\calO_{E_3}$. 
Each map \eqref{local.model.central.eq} is finite birational between integral schemes.
We claim that the map \eqref{local.model.central.eq} for $i=1$ is an isomorphism, which implies the proposition.
By the argument of Proposition \ref{normality.prop}, it suffices to show that the induced maps on tangent spaces at closed points of the special fibers are injective.
Now, observe that the reduced neutral component of the Beilinson--Drinfeld Grassmannian attached to $\ucG_{i,\bbf}$ is contained in the Beilinson--Drinfeld Grassmannian of $\ucG_{i, \bbf, \der}$, since this is true at the level of the generic fibers (which are classical affine Grassmannians), compatibly with the map \eqref{local.model.central.eq}. 
In particular, by translating closed points in the special fibers of the local models, it suffices to check that
$$ \Fl_{G^\flat_{3,\der} \bbf} \to \Fl_{G^\flat_{1,\der} \bbf}$$
is formally étale. This is true under the assumption $p \nmid \frac {\#\pi_1(G_{1,\on{der}}) }{\#\pi_1(G_{3,\on{der}})}$, because $G^\flat_{3,\der}\to G^\flat_{1,\der}$ has an étale kernel, comp.~\cite[6.a]{PR08}.
\xpf

For defining canonical integral models of Shimura varieties of abelian type, under the assumptions of tame ramification and $p>2$, Kisin-Pappas need to approach local models in terms of embeddings, see \cite[\S2.3, \S3.2]{KP18}. 

\defi
The tuple $(G,\{\mu\})$ is called {\em of abelian type} if there is a central lift $(G_1,\{\mu_1\})$ of $(G_{\on{ad}},\{\mu_{\on{ad}}\})$ endowed with a closed embedding $\rho_1\co G_1\to \GL_n$, $n\geq 1$ such that $\{\rho_1\circ \mu_1\}=\{\varpi_d^\vee\}$, where $\varpi_d^\vee$ denotes the $d$-th minuscule coweight of $\text{GL}_n$ for some $1\leq d \leq n-1$.  
The central lift $(G_1,\{\mu_1\})$ is also called {\em of Hodge type}.
\xdefi

First note that every Shimura datum of abelian type in the sense of \cite{Del79} (see also \cite{Kis10}) gives rise to a tuple of abelian type as above. 
Further, our definition coincides with that given in \cite[II, \S3.11]{Lou20}, but appears to differ from that of \cite[\S2.7]{HPR20} in the following way. 
Let us, for simplicity, assume that $G_{\on{ad}}$ contains no $F$-simple factor over which $\{\mu_{\on{ad}}\}$ becomes trivial. 
Then the classification of Hodge embeddings due to Deligne, see \cite[1.3.8, table 1.3.9]{Del79} (compare with \cite[Lem.~3.4.13]{Kis10} and \cite[Lem.~4.6.22]{KP18}) implies that $\rho_1$ is minuscule (as required in \cite[\S2.7]{HPR20}). 
The main difference here is not requiring the existence of an isogeny $G_{1,\on{der}}\rightarrow G_{\on{der}}$ because we want our class to be stable under central lifts.

Next we note that since $G$ is assumed to be tamely ramified, we can arrange for $G_1$ to be tamely ramified as well. 
So the canonical maps $G\to G_\ad\gets G_1$ extend to maps of $\calO_F[t]$-group schemes $\ucG_{\bbf}\to \ucG_{\bbf,\on{ad}}\gets \ucG_{\bbf,1}$, and hence to maps of PZ local models $\bbM\to \bbM_\ad\gets \bbM_1$ defined over the ring of integers of the compositum $E\cdot E_1$ where $E_1$ is the reflex field of $\{\mu_1\}$. 
Also note that adding the center of $\text{GL}_n$ to our central lift $G_1$ changes neither the condition on $\rho_1$ nor the PZ local model $\bbM_1$ in view of Proposition \ref{local.model.central}. 

The importance of our central lift of Hodge type is that $\rho_1$ extends by \cite[Prop.~1.3.3]{KP18} to a (not necessarily closed) immersion of parahoric group schemes 
$$\calG_{\bbf,1} \rightarrow \calG\!\calL_n,$$
which is heavily based on work of Landvogt (beware that \cite{KP18} use the notation $\calG_x$ for the fixer group scheme of $x$ and reserve $\calG_x^0$ for its parahoric neutral component). 
Then the same authors construct in \cite[Prop.~2.3.7]{KP18} a uniquely determined closed embedding:
$$\bbM_1 \hookrightarrow \bbM_{\on{lat}}$$
where we set $\bbM_{\on{lat}}:=\bbM(\on{GL}_n, \{\varpi_d^\vee\}, \calG\!\calL_n)\otimes \calO_{E_1}$.
We remark that the symplectic embeddings used in the given reference and the hypothesis $p \nmid \# \pi_1(G_{1,\on{der}})$ are unnecessary, the former pertaining to later applications to Shimura varieties and the latter to ensure normality of $\bbM_1$. Here we will give a closer look at the possibilities for the geometry of this scheme, analyzing all possible cases.

\prop\label{Hodge.classification}
Let $(G_1,\{\mu_1\})$ be a central lift of Hodge type as above and let $\bbM_1$ be the PZ local model attached to the corresponding LM triple.
Then the following properties hold:
\begin{enumerate}
\item If $p>2$ or $G_{\on{ad}}$ has no $D$-factors, then $\bbM_1$ is always normal and only depends on $(G,\{\mu\},\calG_\bbf)$ up to extending scalars.

\item If $p=2$ and $(G_{\on{ad}},\{\mu_{\on{ad}}\})$ is $\breve{F}$-simple of type $D_n^\bbH$, $n\geq 5$, then $\bbM_1$ only depends on $(G,\{\mu\},\calG_\bbf)$ up to base change, but will be non-normal for sufficiently large $\bar{ \mu}$.
\item If $p=2$ and $(G_{\on{ad}},\{\mu_{\on{ad}}\})$ is $\breve{F}$-simple of type $D_{2m+1}^\bbR$, $m\geq 2$, then $\bbM_1$ is always normal and only depends on $(G,\{\mu\},\calG_\bbf)$ up to base change.

\item If $p=2$ and $(G_{\on{ad}},\{\mu_{\on{ad}}\})$ is $\breve{F}$-simple of type $D_{2m}^\bbR$, $m\geq 2$, then we can always choose $(G_1,\{\mu_1\},\calG_{\bbf,1})$ and $\rho_1$ such that $\bbM_1$ is normal. For sufficiently large $\bar{ \mu}$, we can simultaneously choose $(G_2,\{\mu_2\},\calG_{\bbf,2})$ and $\rho_2$, such that $\bbM_2$ is non-normal.
\end{enumerate}
\xprop

\pf
By Corollary \ref{local.models.corollary} it suffices to examine the normality of Schubert varieties in the special fiber of $\bbM_1$ and for this we need to understand when $p$ divides the order of $Z$, the kernel of $G_{\on{sc}} \rightarrow G_{1,\on{der}}$. Inspecting Deligne's table, see \cite[table 1.3.9]{Del79}, we see that $Z$ is always a multiplicative $2$-group and trivial if $G_{\on{sc}}$ has no $D$-factors. This gives (1). For simple orthogonal adjoint groups, the pullback of $\rho_1$ to $G_{\on{sc}}=\Res_{F'/F}\on{Spin}_{2n}$ is the restriction of scalars of one of the three minuscule representations of $\on{Spin}_{2n}$, that is, the two half-spin irreducible factors of the faithful spin representation of $\on{Spin}_{2n}$ and the pulled back vector representation of $\on{SO}_{2n}$. Moreover, our inspection of \cite[table 1.3.9]{Del79} reveals that we can only use the vector representation for $D_n^\bbH$ and the half-spins for $D_n^\bbR$. In (2), the kernel $Z$ is a certain $2$-group independent of the choice of a Hodge lift. For (3), $Z$ is always trivial because the half-spin representations are faithful if $n=2m+1$ is odd. For (4), we can choose our Hodge lift such that $\rho_1$ restricts to the faithful spin representation (= sum of the two half-spins), but we can also choose some other Hodge lift with $\rho_2$ restricting to an half-spin representation, whose kernel is a non-trivial $2$-group.
\xpf

\begin{rmk}
Concerning integral models of Shimura varieties of abelian type, it seems that the hypothesis $p \nmid \# \pi_1(G_{\on{der}})$ in \cite[Thm.~0.4]{KP18} can be removed, as long as one replaces the PZ local model in the statement by its (weak) normalization. For the $p=2$, $D^\bbH$ case with $\bar{\mu}$ large, the proposition seems to indicate some additional work might be needed in order to construct such integral models, so as to circumvent the fact that the Hodge embedding defines a non-normal orbit closure.
\end{rmk}

\begin{rmk}
Let us comment on the relation with the Scholze--Weinstein conjecture \cite[Conj.~21.4.1]{SW20}. 
The $B_{\on{dR}}$-affine Grassmannian $\Gr^{\on{dR}}_{\calG_{\bbf}}\to \calO_F^\diamondsuit$ is an ind-proper v-sheaf.
The conjecture states that for $\{\mu\}$ minuscule, the $\{\mu\}$-bounded sub-v-sheaf $\Gr^{\on{dR}}_{\calG_{\bbf}, \{\mu\}}\to \calO_E^\diamondsuit$ is representable by a unique flat projective $\calO_E$-model with reduced special fiber of the variety of type $\{\mu\}$-parabolics, called {\it the local model}.
This conjecture is proven in \cite{HPR20, Lou20} for many cases with $(G,\{\mu\})$ of abelian type, and in \cite[Theorem 1.1]{AGLR22} for general reductive groups $G$, relying on \cite[Theorem 1.3]{GL22} for the reducedness of special fibers in some cases of wildly ramified groups $G$. 
The local model singled out by the conjecture is the weak normalization $\tilde{\bbM}$ of $\bbM$ and coincides with the PZ local model of some $z$-extension, simply the local model in the sense of \cite[\S2.6]{HPR20}.
\end{rmk}

\begin{appendix}

\section{Frobenius ind-splitting} \label{Frob_ind-split}
Fix a field $k$ of characteristic $p>0$.
Here we revisit the notion of Frobenius splittings and prove several basic lemmas regarding this technique in the realm of ind-schemes over $k$.

\defi
A $k$-scheme $X$ is said to be (Frobenius) split if the morphism of $\calO_X$-modules $\calO_X \rightarrow F_*\calO_X$ admits a section $s$, where $F$ denotes the absolute Frobenius morphism. 
A closed subscheme $Y$ of $X$ is said to be compatibly split if the splitting of $X$ descends to that of $Y$. 
Finally, we say that an ind-scheme $X$ is ind-split (resp.~compatibly ind-split with an ind-closed sub-ind-scheme) if it admits a presentation $X=\on{colim}X_i$ (resp.~as well as $Y=\on{colim}Y_i$) by simultaneously compatibly split schemes.
\xdefi

\lemm
Given a collection $X_i$ of simultaneously compatibly ind-split ind-closed sub-ind-schemes of $X$, finite intersections and finite unions are also simultaneously compatibly split.
\xlemm
\pf
This is known in the case of schemes (\cite[Prop.~1.2.1]{BK07}), and it generalizes to that of ind-schemes by taking appropriate presentations.
\xpf

Thanks to \cite[(3.7)]{BL94} or \cite[7.11.3]{BD91}, we have a good notion of sheaves of modules on ind-schemes $X=\on{colim}X_i$, namely obtained as a family of compatible $\calO_{X_i}$-modules in the obvious way. 
For an ind-proper ind-scheme $X$ over $k$ equipped with a coherent $\calO_X$-module $\calM$ arising from coherent $\calO_{X_i}$-modules $\calM_i$, we define $H^n(X,\calM) := \on{lim}H^{n}(X_i,\calM_i) $ for $n\geq 0$.
This definition is sensible as the cohomology groups are finite-dimensional and thus $R^n\!\on{lim}$ vanishes for $n>0$.

\lemm
Let $Y\subset X$ be a closed immersion of compatibly ind-split ind-proper ind-schemes over $k$. 
If $\calL$ is an ample line bundle on $X$, then $H^0(X, \calL)\rightarrow H^0(Y, \calL)$ is surjective and $H^{>0}(X,\calL)=H^{>0}(Y,\calL)=0$.
Additionally, if $Y\subset X$ is a closed immersion, then $H^1(X,\calI_Y\otimes \calL)=0$ where $\calI_Y$ denotes the ideal sheaf defining $Y$.
\xlemm

\pf
At finite level, this is just \cite[Thm.~1.2.8]{BK07} and it follows in general by taking projective limits.
\xpf

The next results study the implications of certain splittings for the graded algebra $H^0(X, \calL^{\bullet})$, going back to Ramanathan \cite{Ram87}, but we mostly follow the treatment of \cite[\S1.5]{BK07}.

\prop\label{ind-split-diagonal}
Let $Y\subset X$ be a closed immersion of compatibly ind-split ind-proper ind-schemes over $k$, and assume that the diagonal $\Delta_X$ is compatibly ind-split with $X\times X$. 
Given an ample line bundle $\calL$ on $X$, the graded $k$-algebra $H^0(Y, \calL^{\bullet})$ is generated by its degree one elements, and $\calL$ defines an ind-closed immersion of $Y$ into $\bbP(H^0(Y, \calL)^{\vee})$ where $H^0(Y, \calL)^{\vee}:=\on{colim}H^0(Y_i, \calL)^{\vee}$. 
\xprop

\pf
Observe that due to the previous lemma and our hypothesis, we get surjectivity of the map $H^0(X, \calL^n)\otimes H^0(X,\calL)\rightarrow H^0(X, \calL^{n+1})$ which implies the claim as long as $Y=X$, by induction on $n$. 
If $Y$ is not necessarily equal to $X$, we still have an epimorphism $H^0(X,\calL^n)\rightarrow H^0(Y, \calL^n)$. 
The projective embedding is given by taking the colimit of the resulting closed immersions for a compatibly split presentation.
\xpf

Let us recall some terminology regarding commutative graded algebras and modules (compare with \cite[Def.~1.5.5]{BK07}).

\defi\label{quadratic-graded-algebras-and-modules}
 A commutative $\bbZ_{\geq 0}$-graded $k$-algebra $A_{\bullet}$ is called quadratic if $A_0=k$, if it is generated by $A_1$ and if the kernel $K_{\bullet}$ of the induced surjection $S^{\bullet}A_1\rightarrow A_{\bullet}$ is generated by $K_2$. 
 An $A_{\bullet}$-graded module $M_{\bullet}$ is said to be quadratic if it is generated by $M_0$ and the kernel $K_{\bullet}$ of $M_0\otimes A_{\bullet} \rightarrow M_{\bullet}$ is generated by $K_1$.
\xdefi

The next result subsumes \cite[Prop.~2.7, Prop.~2.19]{Ram87} and is the ind-scheme version of \cite[Prop.~1.5.8]{BK07}.

\prop\label{ind-split-triple-diagonal}
Let $Z \subset Y \subset X$ be closed immersions of simultaneously compatible ind-split ind-proper ind-schemes and $\calL$ be an ample line bundle on $X$. Suppose moreover that $\Delta_{X^2}\times {X}$, $X\times \Delta_{X^2}$, $Y\times \Delta_{X^2}$ and $Z\times \Delta_{X^2}$ are simultaneously compatibly ind-split in $X^3$. Then $H^0(Y, \calL^{\bullet})$ is a quadratic graded algebra and $H^0(Z, \calL^{\bullet})$ is a quadratic graded module over $H^0(Y, \calL^{\bullet})$.
\xprop

In geometric terms, this tells us that the projective embeddings of $Z$, $Y$ and $X$ determined by $\calL$ are given by quadratic homogeneous polynomials and the transition morphisms are defined by linear ones.

\pf
If we intersect the given ind-schemes with $X\times \Delta_{X^2}$, the conditions of Proposition \ref{ind-split-diagonal} are satisfied and hence the graded algebras in the statement are generated by its degree $1$ elements. To show that $A_{\bullet}:=H^0(X, \calL^{\bullet})$ is quadratic, we consider the Mayer-Vietoris short exact sequence $0 \rightarrow \calI_{\Delta_{X^2}\times X \cup X\times \Delta_{X^2}} \rightarrow \calI_{\Delta_{X^2}\times X } \oplus \calI_{ X\times \Delta_{X^2}} \rightarrow \calI_{\Delta_{X^3}} \rightarrow 0$, where the ideals of definition are with respect to $X^3$, and we tensor it with $\calL^{n_1}\boxtimes\calL^{n_2}\boxtimes \calL^{n_3}$ for some integers $n_1,n_2,n_3\geq 1$. 
Let $K_{n_1,n_2,n_3}$ be the kernel of $A_{n_1}\otimes A_{n_2} \otimes A_{n_3} \to A_{n_1+n_2+ n_3}$, and analogously for $K_{n_1,n_2}$ and $K_{n_2,n_3}$, respectively. 
By taking cohomology, we arrive at surjectivity of the map $K_{n_1,n_2}\otimes A_{n_3}\, \oplus\, A_{n_1}\otimes K_{n_2,n_3}\rightarrow K_{n_1,n_2,n_3}$ by the proof of \cite[Prop.~1.5.8]{BK07}, which implies that $A_{\bullet}$ is quadratic by \cite[Lem.~1.5.7]{BK07}. In order to show that $B_{\bullet}:=H^0(Y, \calL^{\bullet})$ and $C_{\bullet}:=H^0(Z, \calL^{\bullet})$ are quadratic algebras, we repeat the same strategy with the couple $(\Delta_{X^2}\times X , Y\times \Delta_{X^2})$ which intersects in $\Delta_{Y^3}$ and use surjectivity of the transition maps $A_{\bullet}\rightarrow B_{\bullet}\rightarrow C_{\bullet}$ to derive the same formulae for $B_{\bullet}$ and $C_{\bullet}$. By \cite[Rmk.~1.5.6 (iii)]{BK07} every transition map defines a graded quadratic module structure.
\xpf

\section{The quasi-minuscule Schubert scheme for $\PGL_2$}\label{appendix_PGL2}
Let $S_\scon$ (resp.~$S_\ad$) be the quasi-minuscule Schubert variety in the affine Grassmannian $\Gr_{\SL_2}$ (resp.~$\Gr_{\PGL_2}$) over $\bbZ$.
Let $L^{--}\SL_2$ (resp.~$L^{--}\PGL_2$) be the strictly negative loop group (see Section \ref{Neg_Loop_Grp_Sec}) which defines an ind-affine open neighborhood of the base point in $\Gr_{\SL_2}$ (resp.~$\Gr_{\PGL_2}$).
The canonical map $\Gr_{\SL_2}\to \Gr_{\PGL_2}$ induces a scheme theoretically surjective morphism $S_\scon\to S_\ad$, and hence a morphism
\[
\Spec\,A:=L^{--}\SL_2\cap S_\scon \to L^{--}\PGL_2\cap S_\ad =:\Spec\,B,
\]
which identifies $B$ with an integral subdomain of $A$.
The aim of this section is to prove the following result.

\prop\label{quasi_min_PGL2}
There is an isomorphism $A\cong \bbZ[x,y,z]/(z^2+xy)$ under which $B$ is the subring generated by the elements $x, y, 2z, xz, yz$.
\xprop

\coro\label{quasi_min_PGL2_special}
The ring $B\otimes \bbF_2$ is not reduced.
Its reduction identifies with the subring of $A\otimes \bbF_2\cong \bbF_2[x,y,z]/(z^2+xy)$ generated by $x, y, xz, yz$. 
\xcoro
\pf
The element $u=2z$ is not $0$ in $B\otimes \bbF_2$ because $z\not \in B$. 
But its square $u^2=4z^2=-4xy$ is $0$ in $B\otimes \bbF_2$ because $x,y\in B$.
This shows that $B\otimes \bbF_2$ is not reduced. 
Clearly, the image of $B\otimes \bbF_2\to A\otimes \bbF_2\cong \bbF_2[x,y,z]/(z^2+xy)$ is the subring generated by $x,y,xz,yz$.
The kernel of this map is nilpotent because the spectra of all rings are irreducible of Krull dimension $2$.
Hence, the ring $(B\otimes \bbF_2)_\red$ identifies with the desired subring of the integral domain $\bbF_2[x,y,z]/(z^2+xy)$.
\xpf

This corollary shows that the special fiber $S_\ad\otimes \bbF_2$ is not reduced. 
More precisely, the reduction $(S_\ad\otimes \bbF_2)_\red$ is the quasi-minuscule Schubert variety for $\PGL_2$ over $\bbF_2$, but the inclusion $(S_\ad\otimes \bbF_2)_\red\subset S_\ad\otimes \bbF_2$ is strict.

To prove Proposition \ref{quasi_min_PGL2}, we first calculate the ring $A$. 
For this, we consider the Lie algebra $\fraks\frakl_2$ of $\SL_2$ over $\bbZ$.
The nilpotent cone $\frakn$ in $\fraks\frakl_2$ is the closed subscheme of matrices whose determinant is zero.
We choose the isomorphism $\bbA^3_\bbZ\cong \fraks\frakl_2$ given by the map
\[
(x,y,z)\mapsto \begin{pmatrix} z & x \\ y & -z \\ \end{pmatrix},
\]
so that $\{z^2+xy=0\}\cong \frakn$ as schemes over $\bbZ$.

\lemm\label{quasi_min_SL2}
Let $e\in \Gr_{\SL_2}(\bbZ)$ denote the base point. 
The map $\frakn \to \Gr_{\SL_2}$, $X\mapsto (1+t^{-1}X)\cdot e$ induces an isomorphism $\frakn \cong L^{--}\SL_2\cap S_{\scon}$, that is, an isomorphism $\bbZ[x,y,z]/(z^2+xy)\cong A$ on coordinate rings.
\xlemm
\pf
The map $\frakn\to L^{--}\SL_2$, $X\mapsto 1+t^{-1}X$ is well-defined and a closed immersion.  
It induces an isomorphism onto the closed subscheme $(L^{--}\SL_2)_{[-1,1]}$ of $L^{--}\SL_2$ of all matrices $M=1+t^{-1}M_1+t^{-2}M_2+\ldots$ such that $M_i=0$ and $(M^{-1})_i=0$ for $i\geq 2$.
We now regard $L^{--}\SL_2$ via the map $g\mapsto g\cdot e$ as an open sub-ind-scheme of $\Gr_{\SL_2}$, see \cite[Lem.~2]{Fal03} (cf.~Lemma \ref{Open_Lem_Split}).
It remains to show $L^{--}\SL_2\cap S_\scon=(L^{--}\SL_2)_{[-1,1]}$ as subschemes of $\Gr_{\SL_2}$.

Recall the lattice interpretation of the affine Grassmannian, see \cite[p.~42]{Fal03} (cf.~\cite[p.~697]{Gor01}).
For any ring $R$, the $R$-valued points of $\Gr_{\SL_2}$ are given by $R\pot{t}$-lattices $\La\subset R\rpot{t}^2$ such that $\det \La=R\pot{t}$ in $R\rpot{t}$. 
We denote by $\La_{0,R}=R\pot{t}^2$ the standard lattice which corresponds to the base point.
Let $\Gr_{\SL_2, [-1,1]}$ denote the closed subscheme of $\Gr_{\SL_2}$ of $R\pot{t}$-lattices $\La$ such that $t\La_{0,R}\subset \La\subset t^{-1}\La_{0,R}$.
A direct computation on $R$-valued points shows $(L^{--}\SL_2)_{[-1,1]}=L^{--}\SL_2\cap \Gr_{\SL_2, [-1,1]}$.
Recall that $S_\scon$ is defined as the scheme theoretic closure of the orbit map $L^+\SL_2\to \Gr_{\SL_2}$, $g\mapsto g\cdot \left (\begin{smallmatrix} t & 0 \\ 0 & t^{-1} \\ \end{smallmatrix} \right )\cdot e$. 
We see that $S_\scon \subset \Gr_{\SL_2, [-1,1]}$, and hence that $L^{--}\SL_2\cap S_\scon$ is a closed subscheme of $(L^{--}\SL_2)_{[-1,1]}$.
Since both are integral of Krull dimension $2+\dim \bbZ=3$, they must be equal.   
\xpf
 
In order to calculate the subring $B$ of $A\cong \bbZ[x,y,z]/(z^2+xy)$, we consider the adjoint representation of $\SL_2$:
The map $g\mapsto (x\mapsto gxg^{-1})$ induces a morphism of $\bbZ$-group schemes $ \SL_2\to {\rm GL}(\fraks\frakl_2)=\GL_3$ given by
\begin{equation}\label{adjoint_rep}
\Ad\co 
\begin{pmatrix} 
a & b \\ 
c & d \\
\end{pmatrix}
\mapsto 
\begin{pmatrix} 
1 + 2bc & -ac & bd \\ 
-2ab & a^2 & -b^2 \\
2cd & -c^2 & d^2 \\
\end{pmatrix}
,
\end{equation}
where we use the ordered basis $\left (\begin{smallmatrix} 1 & 0 \\ 0 & -1 \\ \end{smallmatrix} \right )$, $\left (\begin{smallmatrix} 0 & 1 \\ 0 & 0 \\ \end{smallmatrix} \right )$, $\left (\begin{smallmatrix} 0 & 0 \\ 1 & 0 \\ \end{smallmatrix} \right )$ of $\fraks\frakl_2$. 
This map induces a closed immersion $\PGL_2\hookto \SL_3$ of reductive $\bbZ$-group schemes, and hence a closed immersion $\Gr_{\PGL_2}\hookto \Gr_{\SL_3}$ of affine Grassmannians over $\bbZ$.
Therefore, the image of $S_\ad$ in $\Gr_{\SL_3}$ identifies with the scheme theoretic image of $S_\scon$ under $\Ad\co \Gr_{\SL_2}\to \Gr_{\SL_3}$.

\pf[Proof of Proposition \ref{quasi_min_PGL2}]
We identify $A=\bbZ[x,y,z]/(z^2+xy)$ under the isomorphism of Lemma \ref{quasi_min_SL2}.
Combining this with \eqref{adjoint_rep} gives
\[
\Ad\co 
\begin{pmatrix} 
1+t^{-1}z & t^{-1}x \\ 
t^{-1}y & 1-t^{-1}z \\
\end{pmatrix}
\mapsto 
\begin{pmatrix} 
1 + 2t^{-2}xy & -(1+t^{-1}z)t^{-1}y & t^{-1}x(1-t^{-1}z) \\ 
-2(1+t^{-1}z)t^{-1}x & (1+t^{-1}z)^2 & -t^{-2}x^2 \\
2t^{-1}y(1-t^{-1}z) & -t^{-2}y^2 & (1-t^{-1}z)^2 \\
\end{pmatrix}
.
\]
As this formula holds on $R$-valued points, the ring $B$ is precisely the subring of $A$ generated by the monomials in $x,y,z$ appearing as coefficients of $t^{i}$ for $i=-1,-2$.
An inspection of this matrix using $z^2+xy=0$ and $\bbZ^\x={\pm1}$ shows that $B=\bbZ[x,y,2z,xz,yz]$ as a subring of $A$.
\xpf

\section{Minimal nilpotent orbits in twisted affine Grassmannians}\label{app-minimal-nilpotent-orbits}
Fix an algebraically closed field $k$ of characteristic $0$. Let $H$ be a simply connected or adjoint simple $k$-group of type $A_n$ ($n\geq 2$), $D_n$ ($n \geq 4$), or $E_6$ endowed with a pinning$(H,T_H,B_H,X_H)$. Let $\sigma_0$ be the canonical involution of $\text{Aut}(H,T_H,B_H,X_H)$ induced by the non-trivial involution of the Dynkin diagram of $\Phi(H,T_H)$\footnote{For $H$ of type $D_4$, there are 3 possible choices of involution, and we pick the one fixing $\alpha_1$. Note that since these involutions are all conjugate, other choices lead to isomorphic group-theoretic data.}. Let $M:=H^{\sigma_0}$ be the affine $k$-group deduced from $H$ by taking $\sigma_0$-fixed points. It is smooth of finite type,  connected, reductive, simple, simply connected or adjoint by \cite[Prop.~4.1]{Hai15} (see also \cite[Prop.~A.1]{HR20a}). Moreover, it carries a natural pinning $(M,T_M,B_M,X_M) $ where the middle entries are given by fixed points under the involution and $X_M=X_H$. The root system $\Phi(M,T_M)$ is the set of non-divisible elements in the image of $\Phi(H,T_H)$ under the natural restriction morphism:
$$X^*(T_H)\otimes \bbR \rightarrow X^*(T_M)\otimes \bbR \cong (X^*(T_H)\otimes \bbR)_{\sigma_0} $$
the latter of which will often be identified with $(X^*(T_H)\otimes \bbR)^{\sigma_0}$ via the obvious averaging map.

Let $k\rpot{t}$ be the Laurent series field over $k$ and consider its quadratic Galois extension $k\rpot{u}$ with $u=t^{1/2}$. The restriction of scalars $\Res_{k\rpot{u}/k\rpot{t}} { H_{k\rpot{u}}}$ admits the involution $\sigma:=\sigma_0 \otimes \iota$ where $\iota$ stands for the Galois involution of $k\rpot{u}/k\rpot{t}$. Its fixed points $G:=(\Res_{k\rpot{u}/k\rpot{t}} { H_{k\rpot{u}}})^\sigma$ form a reductive, quasi-split group equipped with a natural pinning $(G,T_G,B_G,X_G)$, see \cite[\S2]{PZ13} or \cite[\S2.1-2.2]{Lou19}. We have an obvious absolutely special parahoric model of $G$ given by the same formula after replacing $k\rpot{u}/k\rpot{t}$ by $k\pot{u}/k\pot{t}$. We will still denote this $k\pot{t}$-group by $G$. It is important to note as well that at a combinatorial level, the groups $G$ and $M$ are not so far from one another, in the sense that $\Phi(M,T_M)$ is the set of non-divisible roots of the relative root system $\Phi(G,S_G)$ via the obvious identification $X^*(S_G)\otimes \bbR=X^*(T_M)\otimes \bbR$.

Our goal is to establish a link between certain nilpotent orbits of $M$ (not necessarily for the adjoint representation) and certain Schubert varieties of $\Gr_G$.  Note that the classification of simply connected tamely ramified reductive groups would force us to consider the case where $\sig_0$ is either the identity or has order $3$ (for the $D_4$ root system and associated triality). However, the material of this section has already been treated in \cite[\S2.10]{MOV05} and \cite[\S8]{HR20a} in those additional cases.

\subsection{Minimal nilpotent orbits of the $M$-module $\frakg_{-1}$}\label{minimal-nilpotent-orbits}

Let $\frakh\otimes k[u,u^{-1}]$ be the algebraic loop algebra of $\frakh$ with the obvious action of $\sigma$ by $\sigma_0$ on the left and Galois conjugation on the right. We let $\frakg$ denote the $\sigma$-invariants of this Kac-Moody algebra - this is a graded version of $\text{Lie}\,G$. The action preserves moreover the obvious $u$-grading and we write $\frakg_{-1}:=\frakh[u^{-1}]^\sigma$. This is acted upon by $M$ in the evident manner and we are going to analyze the structure of this representation as well as some of its nilpotent orbits. In the following, we denote by $\Phi_{M,<}$ the short roots of $\Phi_M:=\Phi(M, T_M)$ and by $\theta_{M,<}$ the unique dominant short root of $\Phi_M$.

\prop\label{repn.nilp.orbit.reduced}
Suppose that $\Phi_G:=\Phi(G,S_G)$ is reduced or, equivalently, that $\Phi_H:=\Phi(H,T_H)$ is not of type $A_{2m}$. The following properties hold:
\begin{enumerate}
\item The $M$-module $\frakg_{-1}$ is irreducible and quasi-minuscule, that is, its highest weight equals $\theta_{M,<}$.
\item Let $v\in\frakg_{-1}$ be any non-zero weight vector. Then the closed orbit $\scrO_{\on{min}}:= \overline{M\cdot v}$ inside the affine space $\frakg_{-1}$ is independent of $v$ and contains the origin. It satisfies the following dimension formula
$$ \dim\scrO_{\on{min}}=2+\# \{ a\in \Phi_M: a+\theta_{M,<} \in \Phi_{M,<} \} $$
and its tangent space $T_0\scrO_{\on{min}}$ at the origin is identified with $\frakg_{-1}$.
\end{enumerate}
\xprop

\pf
The null weight space of our representation equals $\frakt_{G,-1}:=\frakt_H[u^{-1}]^\sigma$ which has dimension equal to the cardinality of $\Delta_{M,<}$, that is, the subset of short positive simple roots. This can be seen by writing down its basis
\begin{equation} \label{averaged_coroot}
u^{-1}h_{\al}-u^{-1}h_{\sigma_0(\al)},
\end{equation}
where $\al \in \Delta_H$ is not $\sigma_0$-invariant and $h_{\al}=[e_{\al},e_{-\al}]$ is the canonical coroot element induced by the choice of the non-zero root vectors $e_{\al}$ belonging to a Chevalley-Steinberg basis of $\frakh$ extending the components of $X_{H}$.
Similarly, we see that the only nonzero weights are (short) roots of the form $\frac{\al +\sigma_0(\al)}{2}$ with $\al \in \Phi_H$. 
Indeed, their weight spaces are $1$-dimensional spanned by
$$ 
u^{-1}e_{\al}-u^{-1}e_{\sigma_0(\al)},
$$
see \cite[4.1.3]{BT84} and compare with \cite[\S2.1]{Lou19} for more explanations and references.
The reducedness hypothesis is crucial here to ensure that $\sigma_0(e_{\al})=e_{\sigma_0(\al)}$ for all roots $\al \in \Phi_H$. 
Since all the roots in $\Phi_H$ have the same length, it follows that the short roots $\Phi_{M,<}$ are those of the form $\frac{\al + \sigma_0(\al)}{2}$ for non-$\sigma_0$-invariant roots $\al \in \Phi_H$, and the unique dominant short root $\theta_{M,<}$ is thus the highest weight of the representation $\frakg_{-1}$. 
Since all the weight-spaces with non-zero weight are 1-dimensional, $\frakg_{-1}$ is the sum of the quasi-minuscule representation of $M$, plus possibly a trivial representation with some multiplicity $m$. 
But it is known that the weight-zero space in the quasi-minuscule representation of $M$ has dimension $\#\Delta_{M,<}$, and thus it follows that $m = 0$. This completes the proof of (1).

Now we consider the minimal\footnote{To actually know that this is the smallest nilpotent orbit of $\frakg_{-1}$ as happens for the adjoint representation, we would need an analogue of the Jacobson-Morozov theorem.} nilpotent orbit $\scrO_{\on{min}}=\overline{M\cdot v}$. 
Since all non-zero weight vectors are extremal by (1) and these are conjugate under the $M$-action, the orbit closure $\scrO_{\on{min}}$ is independent of the choice of the non-zero weight vector $v$.
Further, it is called nilpotent because $v$ belongs to the nullcone of $\frakg_{-1}$. 
In other words, $v$ is an unstable point in the sense of geometric invariant theory, as one can find a cocharacter $\lambda$ of $M$ such that
$$\lim_{t \to 0}\lambda(t)\cdot v =0$$
by the Hilbert-Mumford criterion. 
 This also proves $0 \in \scrO_{\on{min}}$.
 So the tangent space $T_0\scrO_{\on{min}}$, being an $M$-submodule of $\frakg_{-1}$, must be the entire space by irreducibility.

As for computing the dimension, we need to subtract from $\dim M$ the dimension of the stabilizer $Z_M(v)$ of $v$ which is preserved under conjugation by $T_M$. 
This can be done at the level of Lie algebras and then $\frakz_\frakm(v)$ actually decomposes into its intersection with weight spaces for the $T_M$-action. 
Obviously, $\frakz_\frakm(v) \cap \frak{t}_M $ is a hyperplane in $\frak{t}_M$ and hence its cokernel contributes once to the dimension of the minimal nilpotent orbit. 
Now we need to count roots $a \in \Phi_M$ such that $e_a$ annihilates $v$. 
Choosing $v$ to be a highest weight vector, it certainly suffices to have $a  \not \in -\theta_{M,<}+(\Phi_{M,<} \cup \{0 \})$. 

Suppose, on the other hand, that $a +\theta_{M,<} \in \Phi_{M,<} \cup \{0 \}$. If $a =-\theta_{M,<}$ and if we write $\theta_{M,<} = \frac{\psi + \sigma_0(\psi)}{2}$, then $v = u^{-1}e_{\psi} - u^{-1}e_{\sigma_0(\psi)}$ and $e_{-\theta_{M,<}} \cdot v$ is a non-zero multiple of the averaged coroot element (\ref{averaged_coroot}) for $\alpha = \psi$, using that $\{ \psi, \sig_0(\psi)\}$ form a perpendicular orbit pair. If $a +\theta_{M,<}$ is a short root of $M$, then we can write $a=\frac{\al + \sigma_0(\al)}{2}$ without necessarily having $\al\neq \sigma_0(\al)$ and we claim that we can arrange $\alpha+\psi\in \Phi_H$ up to replacing $\alpha$ by its $\sig_0$-conjugate. Otherwise, the bracket $[e_a,e_{\theta_{M,<}}]$ would have to vanish while simultaneously generating the root space of $a+\theta_{M,<}$.
Now if $\alpha \neq \sigma_0(\alpha)$ then $e_a := e_\alpha + e_{\sigma_0(\alpha)}$ and we see that $e_a\cdot v \neq 0$ because after expanding we get a non-zero multiple of $e_{\alpha+\psi}\in\bbC^\times[e_{\alpha}, e_{\psi}]$, which cannot be canceled out since $a+\theta_{M,<}$ is short and hence $\sig_0(\al +\psi)\neq \alpha+\psi$ and also since $\alpha \neq \sigma_0(\alpha)$. If $\alpha = \sigma_0(\alpha)$ then $e_a := e_\alpha$ and similarly $[e_a, u^{-1}(e_\psi - e_{\sigma_0(\psi)})] \neq 0$.
This yields the dimension formula.\footnote{Alternatively, we could have used that the Kac-Moody roots $a$ and $\theta_{M,<}-\delta$ of the Kac-Moody algebra $\frakg$, where $\delta$ is the minimal positive imaginary root, constitute a prenilpotent pair of real roots in the sense of Tits \cite[\S3.2]{Tit87}, so their bracket is non-trivial.}
\xpf

We treat separately the case when $\Phi_G$ is non-reduced, thus of type $BC_n$, for reasons that will become clear to the reader in a moment. We let $\theta_G$ be the highest root of $\Phi_G$ (notice that it must always be divisible). Recall that $\sigma_0$-invariant roots of $\Phi_H$ do not induce roots in $\Phi_M$, but only divisible roots of $\Phi_G$. Under the non-reducedness assumption on $\Phi_G$, the short roots of $\Phi_M$ consist of averages of non-orthogonal $\sigma_0$-orbit pairs of roots, whereas long roots are the averages of the orthogonal pairs.

\prop\label{repn.nilp.orbit.non.reduced}
Suppose that $\Phi_G$ is non-reduced or, equivalently, that $\Phi_H$ is of type $A_{2n}$. The following properties hold:
\begin{enumerate}
\item The $M$-module $\frakg_{-1}$ is irreducible of highest weight $\theta_G$.
\item Let $v\in\frakg_{-1}$ be any extremal weight vector. Then the closed orbit $\scrO_{\on{min}}:= \overline{M\cdot v}$ inside the affine space $\frakg_{-1}$ is independent of $v$ and contains the origin. 
It satisfies the following dimension formula
$$ \dim\scrO_{\on{min}}=1+\# \{ a\in \Phi_M: a+\theta_G \in \Phi_M \} $$
and its tangent space $T_0\scrO_{\on{min}}$ at the origin is identified with $\frakg_{-1}$.
\end{enumerate}
\xprop

\pf
We start by producing a reasonable basis of $\frakg_{-1}$, in the very same spirit of the previous proposition. 
The null weight space is still $\frakt_{G,-1}:=\frakt_H[u^{-1}]^{\sigma}$, it has dimension $\# \Delta_M $ spanned by the basis
$$ u^{-1}(h_{\al}-h_{\sigma_0(\al)}) $$
for any orbit pair $\{\al, \sigma_0(\al) \}$ regardless of their orthogonality behavior. 
All roots $a=\frac{\al+\sig_0(\al)}{2}$ of $\Phi_M$ are multiplicity one weights with weight vectors given by
$$ v_a=u^{-1}(e_\al- e_{\sigma_0(\al)}),$$
where we use a Chevalley-Steinberg basis which must necessarily satisfy the property $\sig_0(e_{\alpha})=\varepsilon_{\alpha}e_{\sig_0(\al)}$, with $\varepsilon_{\alpha}\in \{\pm 1\}$ being a fixed sign. Here $\varepsilon_{\alpha}=1$ if $\alpha\neq \sig_0(\alpha)$ and $\varepsilon_\alpha =-1$ otherwise.

This shows already that $\frakg_{-1}$ is not quasi-minuscule, as $\Phi_M$ is not simply-laced, but we also have extremal vectors of weight $a=\al+\sigma_0(\al)\in \Phi_G \backslash \Phi_M$ for all non-orthogonal non-singleton orbit pairs $\{\al, \sigma_0(\al) \}$, equivalently, all $\sigma_0$-invariant roots $a$ of $\Phi_H$. 
Indeed, these extremal weight spaces are spanned by $v_{a}=u^{-1}e_{a}$, which are fixed by $\sigma$, because $\sigma_0(e_a)=-e_a$. 

Therefore we conclude that $\frakg_{-1}$ contains the highest weight module attached to $\theta_G$. Moreover, since every non-zero weight has multiplicity one, belonging to the highest weight module by saturatedness, the only other possible summand would be the trivial representation. However, it is easy to see that for each $a \in \Delta_M$, $[e_{-\al}+ e_{-\sigma_0(\al)},v_a]=-u^{-1}(h_\al-h_{\sigma_0(\al)})$, whence irreducibility of $\frakg_{-1}$. Indeed, this shows that the entire zero weight space $\frakt_H[u^{-1}]^{\sigma}$ is contained in the module with highest weight $\theta_G$.

As for the remaining assertions on $\scrO_{\on{min}}:=\overline{M\cdot v_{\theta_G}}$, we can argue in the same manner as in the reduced case. Let us take care of the combinatorics. We need to study roots $a$ in $\Phi_M$ such that $a +\theta_G \in \Phi_G \cup \{0\}$ and examine whether $e_av_{\theta_G}\neq 0$. Since $\theta_G\notin \Phi_M$ and $a+\theta_G$ cannot be divisible, we can replace $\Phi_G \cup \{0\}$ by $\Phi_M$. Write $a=\frac{\al+\sig_0(\al)}{2}$ and note that we can arrange $\alpha+\theta_G \in \Phi_H$ just as in the proof of Proposition \ref{repn.nilp.orbit.reduced}. Then $[e_{\al}+e_{\sigma_0(\al)}, u^{-1}e_{\theta_G}]$ is a non-zero multiple of $u^{-1}(e_{\al+\theta_G})$ plus a non-zero multiple of $u^{-1}(e_{\sigma_0(\al) + \theta_G})$, and cancellation cannot occur since $\al \neq \sigma_0(\al)$.
\xpf

\subsection{Quasi-minuscule Schubert variety of $\on{Gr}_G$}
Recall that \cite[Thm.\,6.1]{Hai18} describes the \'echelonnage root and coroot systems $\Phi_{\breve{\Sigma}}$, resp.~$\Phi^\vee_{\breve{\Sigma}}$ for $G$, in terms of the $\sigma_0$-action on $\Phi_H$, resp.~on $\Phi_H^\vee$. We obtain
$$
\Phi_{\breve{\Sigma}}=N'_{\sigma_0}\big(\Phi_H\big),
$$
where the modified norm is defined as in \cite[\S3]{Hai18}. Dually, we get
$$
\Phi^\vee_{\breve{\Sigma}}=\text{res}_{\sigma_0}\big(\Phi_H^\vee\big),
$$
which is given by taking $\sigma_0$-averages and excluding the resulting divisible coroots. We note that parallel to the above, $\Phi_M = {\rm res}_{\sigma_0}(\Phi_H)$.

We are particularly interested in the unique quasi-minuscule coweight $\overline{\psi^\vee}$ of $\Phi^\vee_{\breve{\Sigma}}$. This is obtained from the highest orbit pair $\{\psi^\vee, \sigma_0(\psi)^\vee\}$, the highest orbit pair of $\Phi_H^\vee$ with average $\overline{\psi^\vee} $ of the shortest possible length. In other words, $\{\psi^\vee, \sigma_0(\psi)^\vee\}$ is the set of coroots of $\{\psi, \sigma_0(\psi)\}$, the highest non-singleton orbit pair in $\Phi_H$ if $\Phi_G$ is reduced and the highest non-orthogonal non-singleton orbit pair otherwise (so $\psi$ carries the same meaning as in the previous section). Indeed, in the reduced case, $\psi + \sigma_0(\psi)$ is the highest long root, hence the highest root, of $N'_{\sigma_0}(\Phi_H) = \Phi_{\breve{\Sigma}}$. So $(\psi + \sigma_0(\psi))^\vee = \frac{\psi^\vee +\sigma_0(\psi)^\vee}{2} = \overline{\psi^\vee}$ is the quasi-minuscule coroot for $\Phi_{\breve{\Sigma}}$ (recall $\psi$ and $\sigma_0(\psi)$ are perpendicular). In the non-reduced case, $2(\psi + \sigma_0(\psi))$ is the highest long root, hence the highest root, of $\Phi_{\breve{\Sigma}}$, and a calculation shows the quasi-minuscule coroot is again expressed as $\overline{\psi^\vee}$, compare with \cite[Lem.~3.2]{Hai18}.

We have the following important lemma:

\begin{lem}\label{equality.dim.nilp.orb.qmin.var}
The quasi-minuscule Schubert variety $S_{G,\overline{\psi^\vee}}$ and the minimal nilpotent orbit $\scrO_{\on{min}}$ for $\frakg_{-1}$ have the same dimension.
\end{lem}

\pf
This amounts to establishing the combinatorial identity
$$\langle 2\rho_H, \psi^\vee\rangle=
\begin{cases}
2+\# \{ a\in \Phi_M: a+\theta_{M,<} \in \Phi_{M,<} \}  \, \, \on{if}\, \Phi_G \on{is}\, \on{reduced};\\
1+\# \{ a\in \Phi_M: a+\theta_G \in \Phi_M \}  \, \, \on{else.}
\end{cases}$$

Let us first assume $\Phi _G$ is reduced. Consider the two types of roots in $\Phi_{H, +}$: $\beta$ with $\beta \perp \sigma_0(\beta)$ and $\gamma$ such that $\gamma=\sigma_0(\gamma)$. Write $b := \frac{\beta + \sigma_0(\beta)}{2}$ and $c := \gamma$ for the corresponding positive roots in $\Phi_M$; note that $b \in \Phi_{M,<}$ is a short root and $c \in \Phi_{M, >}$ a long root. Since $\theta_M^\vee := (\theta_{M,<})^\vee=\psi^\vee+\sigma_0(\psi)^\vee$ we have identities
\begin{align}
\langle b, \theta_M^\vee \rangle &= \langle \beta, \psi^\vee \rangle + \langle \sigma_0(\beta), \psi^\vee \rangle \label{b_eq}\\
\langle c, \theta_M^\vee \rangle &= 2\langle \gamma, \psi^\vee \rangle. \label{c_eq}
\end{align}
We claim that (\ref{b_eq}) (resp.~(\ref{c_eq})) takes values in $\{0,1\}$, if $\beta \notin \{\psi, \sigma_0(\psi)\}$ (resp.~$\{0,2\}$). To see this recall that the root $\beta + \sigma_0(\beta)$ (resp.~$\gamma$) of $N'_{\sigma_0}(\Phi_H) = \Phi_{\breve{\Sigma}}$ is not proportional to the highest root $\psi + \sigma_0(\psi)$ of $\Phi_{\breve{\Sigma}}$, so by \cite[VI.1.8, Prop.~25]{Bou68}, we obtain $\langle \beta + \sigma_0(\beta),\overline{\psi^\vee} \rangle$ (resp.~$\langle \gamma, \overline{\psi^\vee} \rangle$) belongs to $\{0,1\}$.

Next we observe that $\langle b, \theta_M^\vee \rangle = 1$ if and only $-b + \theta_{M,<} \in \Phi_{M,<}$, when $\beta \notin \{\psi, \sigma_0(\psi)\}$ (resp.~$\langle c, \theta_M^\vee \rangle = 2$ if and only if $-c + \theta_{M,<} \in \Phi_{M,<}$).  (Note that all roots $a=-b$ (resp.~$a = -c$) appearing in the desired formula are necessarily negative.) If $\langle b, \theta^\vee_M \rangle = 1$, then $s_{\theta_{M,<}}(b) = b - \theta_{M,<} \in \Phi_{M,<}$. If  $\langle c, \theta_M^\vee \rangle = 2$, then $c$ and $s_{\theta_{M, <}}(c) = c - 2 \theta_{M,<}$ are both long roots, so $c-\theta_{M, <}$ is a short root.  
Conversely, if $b$ (resp.~$c$) and $\theta_{M,<}$ are perpendicular, their difference is longer than $\theta_{M,<}$ so in particular is not a short root.

Finally, note that, under the bijection between $\sigma$-orbits in $\Phi_{H}$ and elements of $\Phi_G$, our previous considerations imply by a counting argument that $\langle 2\rho_{ H}, \psi^\vee \rangle$ equals the right side of the combinatorial identity. Indeed, the missing case $b=\theta_{M,<}$ provides the extra summand $2$ in the right side of the identity, whereas positive roots (beware the sign changes) cannot sum with $\theta_{M,<}$ to a short root by maximality.

Now consider the case where $\Phi_G$ is not reduced of type $BC_n$ so that $\Phi_H$ is of type $A_{2n}$. 
As much as we could probably give a combinatorial proof, it is quite simple to verify that the right side equals $2n$ by inspecting \cite[Ch.~VI, Planches~II-III]{Bou68}, whereas a calculation reveals that $\langle2 \rho_H, \psi^\vee\rangle=2n $ as well.
\xpf

We have a natural morphism of reduced (ind)-schemes $\on{exp}\co\frakn_{H,-1} \rightarrow L^{--}H$ induced by the exponential map of Lie algebras, where $\frakn_{H,-1}$ is the set of nilpotent matrices in $\frakh_{-1}$. For $\text{SL}_n$, this can be written as the usual exponential $u^{-1}X\mapsto \sum_{i=0}^\infty \frac{(u^{-1}X)^i}{i!}$, and it follows that the above morphism is a closed immersion. Moreover, it is $\sigma$-equivariant, so we also obtain a closed immersion on fixed points
$$\on{exp}\co\frakn_{G,-1}\longto L^{--}G,$$
where $\frakn_{G,-1}:=\frakn_{H,-1}^\sigma$.
We have the following generalization of \cite[Thm.~8.1, Prop.~8.6]{HR20a}:

\prop\label{exponential.iso}
The morphism $\on{exp}\co \frakn_{G,-1} \rightarrow L^{--}G$ restricts to an isomorphism $\scrO_{\on{min}} \cong L^{--}G\cap S_{G,\overline{\psi^\vee}}$.
\xprop

\pf
Once we show that the image of any extremal weight vector lies in the quasi-minuscule Schubert cell, the result follows immediately from Lemma \ref{equality.dim.nilp.orb.qmin.var}. 
Indeed, we would have a closed immersion between two varieties of the same dimension, so it has to be an isomorphism.

Now for the factorization claim, we must once again divide our approach depending on the reducedness of $\Phi_G$. 
Let us first treat the reduced case. 
We observe that the exponential of $v_{\theta_{M,<}}$ is by definition $x_\psi(u^{-1})x_{\sigma_0(\psi)}(-u^{-1})$. 
But this product of commuting elements comes from an isogeny $\text{SL}_2 \times \text{SL}_2 \rightarrow H$ onto the root group attached to the orbit $\{\psi, \sigma_0(\psi)\}$ of commuting roots. 
Notice that the element $x_a(\pm u^{-1})$ of $\text{Gr}_{\text{SL}_2}$ belongs to the $a^\vee$-Schubert cell. 
Hence, by naturality, we get that $\exp(v_{\theta_{M,<}})$ is sent to the Schubert cell of $\text{Gr}_H$ associated with $\psi^\vee+\sigma_0(\psi)^\vee$. 
But this is exactly the image of $t^{\overline{\psi^\vee}} \in T_G(k\rpot{t})$ in $T_H(k\rpot{u})/T_H(\pot{u})$ under the Kottwitz map, see \cite[(7.3.2)]{Kot97}.

Finally, suppose that $\Phi_G$ is non-reduced. We have the extremal weight vector $v_{\theta_G}$ of $\frakg_{-1}$, whose exponential equals $x_{\psi+\sigma_0(\psi)}(u^{-1})$. 
This element lies again in the Schubert cell of $\text{Gr}_H$ attached to $(\psi + \sigma_0(\psi))^\vee = \psi^\vee+\sigma_0(\psi)^\vee$, so we are done again by \cite[(7.3.2)]{Kot97}.
\xpf

\section{Equivalence of geometric properties} \label{geom_props}

\begin{propo} An $({\bf a}, {\bf f})$-Schubert variety $S_w$, $w \in W/W_{\bf f}$ is normal if and only if it is weakly normal \textup{(}resp.\,\textup{(}S2\textup{)}, resp.\,Cohen-Macaulay, resp.\,Frobenius split if ${\rm char}(k) >0$, resp.,\,it has rational singularities\textup{)}.
\end{propo}

\begin{proof}
Recall from \cite[Prop.~9.7]{PR08} that the normalization $\tilde S_w\to S_w$ is a universal homeomorphism, and that $\tilde S_w$ is Cohen-Macaulay, Frobenius split if $\on{char}(k)>0$, and has rational singularities. Therefore `normal' implies each of those properties. This also handles the equivalence of `normal' with `weakly normal'. To show that `Frobenius split' implies `normal', we invoke \cite[Lem.\,1]{MS87}:  if $Y \to X$ is a proper surjective morphism of irreducible $k$-schemes with connected fibers such that $Y$ is normal and $X$ is Frobenius split, then $X$ is normal. 
We apply this to the Demazure resolution $D(\tilde{w}) \rightarrow S_w$ attached to any reduced decomposition $\tilde w$ of $w$.

If $S_w$ has rational singularities in the sense of \cite[Prop.~9.7]{PR08}, then we know thanks to \cite[Rem.\,9.2]{Kov17} that $S_w$ is Cohen-Macaulay and normal. Also, in general Cohen-Macaulay (resp.\,$S2$) and regularity in codimension $1$ imply normality.

Thus, for the equivalence of `normal' with the remaining properties it suffices to prove the following result.
\end{proof}

\begin{lem} \label{reg_codim_1} 
Every Schubert variety $S_w$ is regular in codimension $1$.
\end{lem}

\pf
We may assume ${\bf f} = {\bf a}$ and $w \in W_{\rm aff}$.  Let $e$ be the base point in $\Fl_{G,\bba}$.
We abbreviate by letting $U^+ = L^+\calG_{\bba}$ and $U^- = L^{--}\calG_{\bba}$. Setting $^xy := xyx^{-1}$, we define for any $x \in W_{\rm aff}$ the groups $U^+_x = \,^xU^- \cap U^+$ and $U^-_x = \,U^- \cap \,^{x^{-1}}U^+$, so that $^xU^-_x = U^+_x$.

Since the open Schubert cell in $S_w$ is a smooth orbit under $U^+ = L^+\calG_\bba$, we are reduced to checking regularity in an open neighbourhood of the point $w' e$ associated to Schubert varieties $S_{w'} \subset S_w$ of codimension $1$. 
Write $w=usv$ as a partial reduced word with $s$ being a simple reflection, such that $w'=uv$ is still a partial reduced word.

For any $x$, the Schubert variety $S_x$ has an open neighborhood of $xe$ of the form $U^+_x xe$.  
This follows from properties of the negative loop groups (Corollary \ref{Open_Cor}; see also \cite[Eqn. (4.2.24)]{Lou19}).
Let $\pi\co D(\tilde w)\to S_w$ be the partial Demazure resolution attached with the partial reduced decomposition $\tilde w= (u,s,v)$.
Therefore, we get an open neighborhood in $D(\tilde{w})$ around $(u,1,v)$ which is isomorphic to 
$$
O := U^+_u u \times U^-_s \times U^+_v v \cong U^+_u \times \,^u U^-_s \times \,^uU^+_v \, uv.
$$
Moving $uv$ to the left, we identify this with
$$
uv \cdot ^{v^{-1}} U^-_u \times \,^{v^{-1}}U^-_s \times \, U^-_v.
$$
Using that $uv$ and $sv$ are partial reduced words, we see that each group factor lies in $U^-$, and thus under the product morphism $\pi$, this maps into $uv U^- \cap S_w$.



We claim that $\pi|_O$ is a monomorphism. 
By the root group decomposition lemma \cite[Prop.~4.2.6]{Lou19}, it is enough to show that the terms of the above product share no affine roots in common.  Comparing the first two terms, this follows from $u < us$.  Comparing the third term with either of the others, it follows because $^vU_v^- = U^+_v$ involves only positive affine roots.


Now the morphism $\pi|_O\co O \rightarrow uvU^- \cap S_w$ is a finite type monomorphism, hence is unramified \cite[Tag 06ND]{StaProj}. 
This implies that $\pi|_O$ induces an isomorphism $\mathcal O_{D(\tilde{w}), (u,1,v)} \cong \mathcal O_{S_w, uv}$, and hence that $S_w$ is regular at $uv$, as follows:
the map $\pi$ is birational by general properties of Demazure resolutions. 
Its restriction $\pi|_O$ is quasi-finite, and so factors as an open immersion followed by a finite morphism by Zariski's Main Theorem.
Hence, the map on local rings is finite, birational and unramified. 
An application of Nakayama's lemma as in the proof of (3)$\Rightarrow$(2) of Proposition \ref{normality.prop} shows that it must be an isomorphism. 
\xpf

\rema 
For Schubert varieties attached to finite Weyl groups, a direct proof of Lemma \ref{reg_codim_1} (not relying on normality) was already known by e.g.\,\cite[Thm.\,A.12.1.10]{LR07}, whose proof proceeds by descending induction on $l(w)$ and does not extend to the affine case. 
While this paper was undergoing revision, a referee pointed us to \cite[Cor.\,3.3]{CFL22}, which similarly handles arbitrary Kac--Moody groups and does not imply Lemma \ref{reg_codim_1} when $\on{char}(k)$ divides $\#\pi_1(G_\der)$.  
\xrema

\end{appendix}

\bibliography{biblio}
\bibliographystyle{alpha}
\end{document}